\magnification=1200
\input amssym.def
\input amssym.tex

\hsize=13.5truecm
\baselineskip=16truept
\baselineskip=16truept
\font\secbf=cmb10 scaled 1200
\font\eightrm=cmr8
\font\sixrm=cmr6

\font\eighti=cmmi8

\font\sixi=cmmi6
\skewchar\eighti='177 \skewchar\sixi='177

\font\eightsy=cmsy8
\font\sixsy=cmsy6
\skewchar\eightsy='60 \skewchar\sixsy='60

\font\eightit=cmti8

\font\eightbf=cmbx8
\font\sixbf=cmbx6

\let\sc=\tensc

\font\eightsc=cmcsc10 scaled 800
\font\secbf=cmb10 scaled 1200
\font\subsecfont=cmb10 scaled \magstephalf
\font\amb=cmmib10

\font\ambi=cmmib10 scaled 700

\newfam\mbfam 

\textfont\mbfam\amb \scriptfont\mbfam\ambi


\def\aa{\def\rm{\fam0\eightrm}%
  \textfont0=\eightrm \scriptfont0=\sixrm \scriptscriptfont0=\fiverm
  \textfont1=\eighti \scriptfont1=\sixi \scriptscriptfont1=\fivei
  \textfont2=\eightsy \scriptfont2=\sixsy \scriptscriptfont2=\fivesy
  \textfont3=\tenex \scriptfont3=\tenex \scriptscriptfont3=\tenex
  \def\sc{\eightsc}
  \def\it{\fam\itfam\eightit}%
  \textfont\itfam=\eightit
  \def\bf{\fam\bffam\eightbf}%
  \textfont\bffam=\eightbf \scriptfont\bffam=\sixbf
   \scriptscriptfont\bffam=\fivebf
  \normalbaselineskip=9.7pt
  \setbox\strutbox=\hbox{\vrule height7pt depth2.6pt width0pt}%
  \normalbaselines\rm}

\def\Proof{\vskip12pt\noindent{\bf Proof.} }

\def\Def#1{\vskip12pt\noindent{\bf Definition #1}}

\def\m@th{\mathsurround=0pt}

\def\cc#1{\hbox to .89\hsize{$\displaystyle\hfil{#1}\hfil$}\cr}
\def\lc#1{\hbox to .89\hsize{$\displaystyle{#1}\hfill$}\cr}
\def\rc#1{\hbox to .89\hsize{$\displaystyle\hfill{#1}$}\cr}

\def\eqal#1{\null\,\vcenter{\openup\jot\m@th
  \ialign{\strut\hfil$\displaystyle{##}$&&$\displaystyle{{}##}$\hfil
      \crcr#1\crcr}}\,}

\def\section#1{\vskip 22pt plus6pt minus2pt\penalty-400
        {{\secbf
        \noindent#1\rightskip=0pt plus 1fill\par}}
        \par\vskip 12pt plus5pt minus 2pt
        \penalty 1000}

\def\subsection#1{\vskip 20pt plus6pt minus2pt\penalty-400
        {{\subsecfont
        \noindent#1\rightskip=0pt plus 1fill\par}}
        \par\vskip 8pt plus5pt minus 2pt
        \penalty 1000}

\def\subsubsection#1{\vskip 18pt plus6pt minus2pt\penalty-400
        {{\subsecfont
        \noindent#1}}
        \par\vskip 7pt plus5pt minus 2pt
        \penalty 1000}

\def\center#1{{\begingroup \leftskip=0pt plus 1fil\rightskip=\leftskip
\parfillskip=0pt \spaceskip=.3333em \xspaceskip=.5em \pretolerance 9999
\tolerance 9999 \parindent 0pt \hyphenpenalty 9999 \exhyphenpenalty 9999
\par #1\par\endgroup}}

\def\\{\hfill\break}

\def\mida#1{{{\null\kern-4.2pt\left\bracevert\vbox to 6pt{}\!\hbox{$#1$}\!\right\bracevert\!\!}}}
\def\midy#1{{{\null\kern-4.2pt\left\bracevert\!\!\hbox{$\scriptstyle{#1}$}\!\!\right\bracevert\!\!}}}

\def\diagint{{\raise1.5pt\hbox{$\scriptscriptstyle\diagup$}\hskip-8.7pt\intop}}

\def\divv{{\rm div}\,}

\def\supp{{\rm supp}\,}

\def\today{${\scriptscriptstyle\number\day-\number\month-\number\year}$}
\footline={{\hfil\rm\the\pageno\hfil${\scriptscriptstyle\rm\jobname}$\ \ \today}}
\def\ifnextchar#1#2#3{\bgroup
  \def\reserveda{\ifx\reservedc #1 \aftergroup\firstoftwo
    \else \aftergroup\secondoftwo\fi\egroup{#2}{#3}}%
  \futurelet\reservedc\ifnch
  }
\def\ifnch{\ifx \reservedc \sptoken \expandafter\xifnch
      \else \expandafter\reserveda
      \fi}
\def\firstoftwo#1#2{#1}
\def\secondoftwo#1#2{#2}
\def\tempswafalse{\let\iftempswa\iffalse}
\def\tempswatrue{\let\iftempswa\iftrue}

\def\cite{\ifnextchar [{\tempswatrue\citea}{\tempswafalse\citeb}}
\def\citea[#1]#2{[#2, #1]}
\def\citeb#1{[#1]}

\def\D{{\Bbb D}}

\def\R{{\Bbb R}}
\def\N{{\Bbb N}}
\def\T{{\Bbb T}}

\def\meas{{\rm meas\,}}

\center{\secbf Global regular axially symmetric solutions to the
Navier-Stokes equations}
\vskip1.5cm

\centerline{\bf Wojciech M. Zaj\c aczkowski}

\vskip1cm
\noindent
Institute of Mathematics, Polish Academy of Sciences,\\
\'Sniadeckich 8, 00-956 Warsaw, Poland\\
E-mail:wz@impan.pl;\\
Institute of Mathematics and Cryptology, Cybernetics Faculty,\\
Military University of Technology, Kaliskiego 2,\\
00-908 Warsaw, Poland
\vskip0.8cm


\vskip1.5cm

\noindent
{\bf Abstract.}
Global regular axially symmetric solutions to the Navier-Stokes equations is
proved. The solution is such that velocity belongs to
$W_2^{2,1}(\Omega\times\R_+)$ and gradient of pressure to
$L_2(\Omega\times\R_+)$, where $\Omega$ is a finite axially symmetric
cylinder in $\R^3$ and the slip boundary conditions are imposed on its
boundary. First we prove the existence of local solutions in the mentioned
spaces for time less or equal $T$. Having such solutions we are able to show
that swirl $(rv_\varphi)$  is the H\"older continuous. This gives
a possibility to show that $\|v'\|_{H^1(\Omega)}$, $v'=(v_r,v_z)$, is bounded
independently
on time near the axis of symmetry. Similar estimate for $v'$ is obtained in
a neighborhood located in a positive distance of the axis of symmetry.
Separately we show that $\|v_\varphi\|_{H^1(\Omega)}$ is bounded by a quantity
independent on time. This implies that the local solution can be prolonged on
intervals $(kT,(k+1)T)$, $k\in\N$, $T>0$. Employing the decay estimates
appropriate for the Navier-Stokes equations we show that however the external
force does not decrease with time the norm $\|v(t)\|_{H^1(\Omega)}$ does not
increase.
\vfil\eject

\section{1. Introduction}

In this paper we prove the existence of global regular axially symmetric
solutions with large swirl to the Navier-Stokes equations. We consider the
motion of an incompressible fluid in an axially symmetric cylinder with the
slip boundary conditions. We generalize the results from \cite{Z1, Z2, Z3},
where the periodic cylinder with respect to the variable along its axis and
the case without any external force are considered. The slip boundary
conditions are necessary because the
main step in a proof of global estimate is the energy type estimate for
vorticity (see \cite{L1, Z4, Z5, Z6}). To get such estimate we need
integration by parts so appropriate boundary conditions for
vorticity are necessary.

In this paper we consider the axially symmetric solutions (see Definition 1.1
below) to the following problem
$$\eqal{
&v_t+v\cdot\nabla v-\nu\Delta v+\nabla p=f\quad &{\rm in}\ \
\Omega^T=\Omega\times(0,T),\cr
&\divv v=0\quad &{\rm in}\ \ \Omega^T,\cr
&v\cdot\bar n=0\quad &{\rm on}\ \ S^T=S\times(0,T),\cr
&\bar n\cdot\D(v)\cdot\bar\tau_\alpha=0,\ \ \alpha=1,2,\quad &{\rm on}\ \
S^T,\cr
&v|_{t=0}=v_0\quad &{\rm in}\ \ \Omega,\cr}
\leqno(1.1)
$$
where $\Omega$ is an axially symmetric cylinder with boundary $S=S_1\cup S_2$,
$x=(x_1,x_2,x_3)$ is the Cartesian system of coordinates in $\R^3$ such that
$x_3$ is the axis of the cylinder $\Omega$.
By $S_1$ we denote the part of the boundary of the cylinder parallel to the
$x_3$-axis and $S_2$ is perpendicular to it. Next,
$v=v(x,t)=(v_1(x,t),v_2(x,t),v_3(x,t))\in\R^3$ is the velocity of the
considered fluid, $p=p(x,t)\in\R$ the pressure,
$f=f(x,t)=(f_1(x,t),f_2(x,t),f_3(x,t))\in\R^3$ the external force field,
$\nu>0$ is the constant viscosity coefficient, $\D(v)=\nabla v+\nabla v^T$ the
dilatation tensor which is the double symmetric part of $\nabla v$, $\bar n$
is the unit outward normal vector to $S$ and $\bar\tau_\alpha$,
$\alpha=1,2$, is a tangent one.

\noindent
To examine axially symmetric solutions to (1.1) we introduce the cylindrical
coordinates $r$, $\varphi$, $z$ by the relations $x_1=r\cos\varphi$,
$x_2=r\sin\varphi$, $x_3=z$. Moreover, we introduce the vectors
$\bar e_r=(\cos\varphi,\sin\varphi,0)$,
$\bar e_\varphi=(-\sin\varphi,\cos\varphi,0)$, $\bar e_z=(0,0,1)$ connected
with the cylindrical coordinates.

\noindent
Then, the cylindrical coordinates of $v$ and $f$ are defined by the relations
$$
v_r=v\cdot\bar e_r,\ \ \ v_\varphi=v\cdot\bar e_\varphi,\ \ \
v_z=v\cdot\bar e_z,\ \ \ f_r=f\cdot\bar e_r,\ \ \
f_\varphi=f\cdot\bar e_\varphi,\ \ \ f_z=f\cdot\bar e_z,
$$
where the dot denotes the scalar product in $\R^3$. Finally, by
$u=v_\varphi r$ we denote a swirl.

To describe the domain $\Omega$ and its boundary in greater details we
introduce the notation
$$\eqal{
&\Omega=\{x\in\R^3:\ r<R,\ |z|<a\},\cr
&S_1=\{x\in\R^3:\ r=R,\ |z|<a\},\cr
&S_2=\{x\in\R^3:\ r<R,\ z\in\{-a,a\}\},\cr}
$$
where $R$ and $a$ are given positive numbers.

\Def{1.1.}
By the axially symmetric solutions we mean such solutions to problem (1.1) that
$$\eqal{
&v_{r,\varphi}=v_{\varphi,\varphi}=v_{z,\varphi}=p_{,\varphi}=0,\cr
&f_{r,\varphi}=f_{\varphi,\varphi}=f_{z,\varphi}=0.\cr}
\leqno(1.2)
$$
In the cylindrical coordinates equations (1.1) for the axially symmetric
solutions can be expressed in the form (see \cite{LL, K})
$$
v_{r,t}+v\cdot\nabla v_r-{v_\varphi^2\over r}-\nu\Delta v_r+\nu{v_r\over r^2}=
-p_{,r}+f_r,
\leqno(1.3)
$$
$$
v_{\varphi,t}+v\cdot\nabla v_\varphi+{v_r\over r}v_\varphi-\nu\Delta v_\varphi+
\nu{v_\varphi\over r^2}=f_\varphi,
\leqno(1.4)
$$
$$
v_{z,t}+v\cdot\nabla v_z-\nu\Delta v_z=-p_{,z}+f_z,
\leqno(1.5)
$$
$$
v_{r,r}+v_{z,z}=-{v_r\over r},
\leqno(1.6)
$$
where $v\cdot\nabla=v_r\partial_r+v_z\partial_z$,
$\Delta u={1\over r}(ru_{,r})_{,r}+u_{,zz}$.

\noindent
Expressing the boundary conditions $(1.1)_4$ in the cylindrical coordinates
yields (see \cite[Ch. 4, Lemma 2.1]{Z4})
$$\eqal{
&v_r=0,\ \ v_{z,r}=0,\ \ v_{\varphi,r}={1\over r}v_\varphi\quad &{\rm on}\ \
S_1,\cr
&v_{r,z}=0,\ \ v_{\varphi,z}=0,\ \ v_z=0\quad &{\rm on}\ \ S_2.\cr}
\leqno(1.7)
$$
Finally, initial conditions assume the form
$$
v_r|_{t=0}=v_r(0),\quad v_\varphi|_{t=0}=v_\varphi(0),\quad
v_z|_{t=0}=v_z(0).
\leqno(1.8)
$$
We prove the existence of global axially symmetric solutions to (1.1) such
that $v\in W_2^{2,1}(\Omega\times\R_+)$, $\nabla p\in L_2(\Omega\times\R_+)$.
The proof is divided into the following steps. First we prove the existence
of a local solution $v\in W_2^{2,1}(\Omega\times(0,T_0))$,
$\nabla p\in L_2(\Omega\times(0,T_0))$, where $T_0$ is sufficiently small.
To show such existence we need that $v(0)\in H^1(\Omega)$ and
$f\in L_2(\Omega\times(0,T_0))$.

\noindent
The existence of the local solutions is proved by the Leray-Schauder fixed
point theorem, where the necessary a priori estimate is possible thanks to
the restriction on $T_0$ (see Lemma 2.5). We are not able to extend the local
solution on the interval $(T_0,2T_0)$ without an additional estimate on
$\|v(T_0)\|_{H^1(\Omega)}$. However, to prove a global existence we have to
show that
$$
\|v(kT_0)\|_{H^1(\Omega)}\le\|v(0)\|_{H^1(\Omega)},\quad k\in\N.
\leqno(1.9)
$$
To show (1.9) we separate cylinder $\Omega$ into two parts: a cylinder of
radius $r_0$, $\Omega_{r_0}=\{x\in\Omega:\ r<2r_0\}$ (a neighborhood of the
axis of symmetry) and the domain $\Omega_{\bar r_0}=\{x\in\Omega:\ r>r_0/2\}$
(a neighborhood located in a positive distance from the $x_3$-axis).

\noindent
Since $\Omega=\Omega_{r_0}\cup\Omega_{\bar r_0}$ we subordinate a partition
of unity corresponding to this division. The most difficult part is to prove
the estimate
$$
\|v'\|_{V_2^1(\Omega_{r_0}^T)}\le a_1,
\leqno(1.10)
$$
where $v'=(v_r,v_z)$, $V_2^1(\Omega^T)$ is defined in Section 2, $T\in\R_+$
and $a_1$ is a constant depending on data.

\noindent
However, (1.10) is proved for the local solution the bound $a_1$ does not
depend on $T$. To prove (1.10) we need the H\"older continuity of swirl $u$
with its vanishing on the axis of symmetry and the restriction
$$
\max_t\max_{\Omega_{r_0}}|u|\le\root{4}\of{5\over8}\nu.
\leqno(1.11)
$$
The results are proved in Section 3.

Hence (1.11) holds for $r_0$ sufficiently small and also for
$v\in L_{10}(\Omega^T)$, so for local solution. This needs step by step in
time approach. In view of (1.11) we have
(1.10) with $a_1=a_1(1/r_0)$, where $a_1$ is an increasing function. Estimate
(1.10) is proved in a series of lemmas (see Lemmas 4.1, 4.2, 4.3, 4.5, 4.6).

\noindent
The next step is to prove the estimate (see Lemmas 6.4, 6.5)
$$
\|v'\|_{V_2^1(\Omega_{\bar r_0}^T)}\le a_2,\quad T\in\R_+,
\leqno(1.12)
$$
where $a_2$ is a constant depending on data. Finally, by Lemma 7.4 we have
$$
\|v_\varphi(t)\|_{H^1(\Omega)}\le a_3,\quad t\in\R_+,
\leqno(1.13)
$$
where $a_3$ depends also on data.

\noindent
However, estimates (4.10), (4.12), (4.13) are proved for the local solution
quantities $a_1,a_2,a_3$ do not depend explicitly on time. The dependence on
time is only by time integral norms of the external force.

\noindent
Combining (1.10), (1.12) and (1.13) we derive the a priori estimate
$$
\|v(t)\|_{H^1(\Omega)}\le a_4,\quad t\in\R_+,
\leqno(1.14)
$$
where $a_4$ depends on data but not on $t$.

\noindent
A priori estimate (1.14) becomes a real estimate for the local solution in the
interval $(0,T_0)$. Since (1.14) is a global type estimate the local solution
can be extended step by step on $\R_+$ (see the Main Theorem below).

\proclaim Main Theorem.
Let $T>0$ be given. Let $k\in\N_0=\N\cup\{0\}$. Assume that
\item{1.} $v(0)\in H^1(\Omega)$,
$f\in L_2(\Omega\times(kT,(k+1)T))\cap L_\infty(\R_+;L_{6/5}(\Omega))$
(see Lemmas 2.2, 2.3)
\item{2.} $u_0=rv_\varphi(0)\in L_\infty(\Omega)$,
$g=rf_\varphi\in L_\infty(\Omega\times\R_+)$ (see Lemma 3.1)
\item{3.} $u_0\in C^\alpha(\Omega)$,
$g\in C^{\alpha,\alpha/2}(\Omega\times(kT,(k+1)T))$,
$\alpha\in(0,1/2]$ (see Lemma 3.2)
\item{4.} $\chi=v_{r,z}-v_{z,r}$, ${\chi(0)\over r}\in L_2(\Omega)$,
$F=f_{r,z}-f_{z,r}$,
${F\over r}\in L_2(\Omega\times(kT,(k+1)T))$,
${v_\varphi(0)\over r}\in L_4(\Omega)$,
${f_\varphi\over\sqrt{r}}\in L_{20/11}(\Omega\times(kT,(k+1)T))$
\item{5.} there exists $r_0>0$ such that
$$
\|u\|_{L_\infty(\Omega_{r_0}\times(kT,(k+1)T))}\le
\root{4}\of{a\over3}\nu,\quad a<1
$$
where $\Omega_{r_0}=\{x\in\Omega:\ r<2r_0\}$.\\
Then there exists a solution to problem (1.1) such that
$v\in W_2^{2,1}(\Omega\times(kT,(k+1)))$,
$\nabla p\in L_2(\Omega\times(kT,(k+1)T))$ and the estimate holds
$$\eqal{
&\|v\|_{W_2^{2,1}(\Omega\times(kT,(k+1)T))}+
\|\nabla p\|_{L_2(\Omega\times(kT,(k+1)T))}\cr
&\le c(B_0(1/r_0,d_4,d_5,d_6,d_7,d_8,d_9)+\|f\|_{L_2(\Omega\times(kT,(k+1)T))}
\cr
&\quad+\|v(kT)\|_{H^1(\Omega)})\cr}
\leqno(1.15)
$$
where (see Lemma 8.2)
$$\eqal{
&\|v(kT)\|_{H^1(\Omega)}\le c(B_0+X(0)e^{-\nu_0kT}),\cr
&\chi^2(0)={1\over r^2}\bigg\|{v_\varphi^2(0)\over r}\bigg\|_{L_2(\Omega)}^2+
\bigg\|{\chi(0)\over r}\bigg\|_{L_2(\Omega)}^2\cr}
$$
and $B_0$ is an increasing positive function of its arguments, where $d_4$ is
introduced in (2.1), $d_5$ in (2.2), $d_6$ in (3.1), $d_7$ in (2.15), $d_8$
and $d_9$ in Lemma 8.2.

\noindent
In the famous paper of Caffarelli, Kohn, Nirenberg \cite{CKN} there is shown
that the singular set for solutions to the Navier-Stokes equations might
have at most an one-dimensional Hausdorff measure. Therefore is can be
expected that for axially symmetric solutions such singularity might appear
on the axis of symmetry. It is shown in this paper that $v_r$, $v_\varphi$,
$u$ vanish on the axis of symmetry, so the considered solution behaves there
very regularly. In \cite{Z8} existence of global regular axially symmetric
solutions with prescribed sufficiently small initial swirl in some
neighborhood of the axis of symmetry is proved. Then the property is
preserved. In this paper the property is shown without any smallness
restrictions on the initial data. Otherwise the necessary a priori estimate
can not be derived.

In this paper we generalize the result from \cite{Z1--Z3} to the problem
with an external force and the slip boundary conditions on whole boundary.
The result is appropriate for examining stability of the axially symmetric
solutions.

\noindent
The generalization is not trivial because the external force has a strong
influence on any solution to the Navier-Stokes equations. The external force
has an opposite influence to the dissipation. Therefore a global existence
can be proved if the dissipation prevails the influence of the external
force. To escape restrictions that velocity and the external force vanish as
time converges to infinity we prove existence step by step on each finite
time interval $[kT,(k+1)T]$, $k\in\N_0$. This approach needs that the
following inequality must be shown
$$
\|v((k+1)T)\|_{H^1(\Omega)}\le\|v(kT)\|_{H^1(\Omega)},\quad k\in\N_0,
\leqno(1.16)
$$
which follows from the decay properties of the Navier-Stokes equations (see
Sections 7 and 8). We have to emphasize that to derive (1.16) $T$ must be
sufficiently large. Moreover, the approach needs that existence must be
proved in each time  interval $[kT,(k+1)T]$, $k\in\N_0$, separately. This
follows from Lemma 2.5 and Theorem 8.1.

\noindent
The Main Theorem says that however $v_r$, $v_\varphi$ vanish on the axis
of symmetry $v_z$ remains large and is only restricted by estimate (1.15).
Hence, considering an inflow-outflow conditions on $S_2$ seems that $v_z$
can be made as much as we want.

\noindent
This paper is organized in the following way. In Section 2 there are
formulated energy type estimates for solutions to problem (1.1) without
showing existence (see Lemmas 2.2, 2.3). Existence is very well presented in
\cite{CKN, L2, T}. However, in our case, we do not need existence of weak
solutions because existence in each time interval $(kT,(k+1)T)$ is proved by
the Leray-Schauder fixed point theorem. For this we need only an appropriate
a priori estimate. Moreover, in Section 2 there are formulated problems for
$u$ (see (2.23)) and $\chi$ (see (2.25)). Finally local existence of solutions
to (1.1) $(v\in W_2^{2,1}(\Omega\times(kT,(k+1)T))$,
$\nabla p\in L_2(\Omega\times(kT,(k+1)T))$, $k\in\N_0$) is proved in
Lemma 2.5.
In Section 3 boundedness and the H\"older continuity of swirl is proved (see
Lemmas 3.1, 3.2).
In Section 4 a priori estimate in a neighborhood of the axis of symmetry is
proved. A crucial point of getting the estimate is restriction (1.11) (see
(4.30)) which implies a different treatment near the axis
of symmetry and far of it.

\noindent
In Section 6 in view of estimate $\|\chi/r\|_{V_2^0(\Omega_{r_0}^t)}\le A$
(see (6.1)) there is proved that
$$
\|v'(t)\|_{H^1(\Omega)}\le cA,
\leqno(1.17)
$$
where $A$ is defined by (6.2). Finally, in Section 7 the estimate is found
$$
\|v_\varphi(t)\|_{H^1(\Omega)}\le\varphi(A_*),
\leqno(1.18)
$$
where $A_*$ is introduced in (6.5).

\noindent
Thanks to estimates (1.17) and (1.18) global existence of regular solutions
such $v\in W_2^{2,1}(\Omega\times(kT,(k+1)T))$,
$\nabla p\in L_2(\Omega\times(kT,(k+1)T))$ is proved step by step in Section 8.

\section{2. Notation and auxiliary results}

By $c$ we denote a generic constant which changes its value from line to line.
A constant $c_k$ with index $k$ is defined by the first formula, where it
appears. By $\phi$ we denote the generic functions which changes its form
from formula to formula and is always positive and increasing function. By
$c(\sigma)$ we denote a generic constant increasing with $\sigma$.

\noindent
We use also the notation
$$\eqal{
&\Omega_\varepsilon=\{x\in\Omega:\ \varepsilon<r\},\cr
&\Omega_\varepsilon\cap\supp\zeta=\Omega_{\varepsilon,\zeta},\quad
\Omega_\varepsilon\cap\supp\zeta_{,r}=\Omega_{\varepsilon,\zeta_{,r}},\quad
\Omega\cap\supp\zeta_{,r}=\Omega_{\zeta,r}.\cr}
$$

\Def{2.1.}
By $V_2^k(\Omega^T)$, $k\in\N\cup\{0\}\equiv\N_0$, we denote a space of
functions with the finite norm
$$
\|u\|_{V_2^k(\Omega^T)}=\|u\|_{L_\infty(0,T;H^k(\Omega))}+
\|\nabla u\|_{L_2(0,T;H^k(\Omega))},
$$
where $H^0(\Omega)=L_2(\Omega)$ and $H^k(\Omega)=W_2^k(\Omega)$ is a Sobolev
space with the finite norm $\|u\|_{H^k(\Omega)}=
\big(\sum_{|\alpha|\le k}\intop_\Omega|D_x^\alpha u|^2dx\big)^{1/2}$, where
$D_x^\alpha=\partial_{x_1}^{\alpha_1}\partial_{x_2}^{\alpha_2}
\partial_{x_3}^{\alpha_3}$, $|\alpha|=\alpha_1+\alpha_2+\alpha_3$,
$\alpha_i\in\N_0$, $i=1,2,3$.

\proclaim Lemma 2.2.
Assume that $v(0)\in L_2(\Omega)$, $f\in L_\infty(\R_+;L_{6/5}(\Omega))$,
$g=f\cdot\eta$, $u=v\cdot\eta$, $\eta=(-x_2,x_1,0)$,
$u(0)\in L_\infty(\Omega)$.
Assume that there exist positive constants $d_0$, $d_1$, $d_2$ such that
$$
\|v(0)\|_{L_2(\Omega)}\le d_0,\quad
\|f\|_{L_\infty(\R_+;L_{6/5}(\Omega))}\le d_2,
$$
$$
\sup_{k\in\N_0}\sup_{t\in(kT,(k+1)T]}\bigg|\intop_{kT}^t
\intop_\Omega gdxdt'+\intop_\Omega u(kT)dx\bigg|\le d_1.
$$
Assume that $\nu_*={\nu\over c_k}$, where $c_k$ is the constant from the
Korn inequality (see (2.4)). Assume that $\nu_*=\nu_1+\nu_2$, $\nu_i>0$,
$i=1,2$.\\
Assume that $T>0$ is fixed and $k\in\N_0$. Then for the weak solutions to
problem (1.1) we have the estimates
$$
\|v(t)\|_{L_2(\Omega)}^2\le{d_3^2\over1-e^{-\nu_1T}}+e^{-\nu_1t}d_0^2\equiv
d_4^2,\quad t\in\R_+,
\leqno(2.1)
$$
where $d_3={1\over\nu_1}\big({1\over\nu_*}d_2^2+2\nu d_1^2\big)$ and
$$\eqal{
&\|v\|_{V_2^0(\Omega\times(kT,t))}^2\le{d_3^2\over\nu_2}e^{\nu_1T}+
\bigg(1+{1\over\nu_2}\bigg)\bigg({2d_3^2\over1-e^{-\nu_1T}}+
e^{-\nu_1kT}d_0^2\bigg)\cr
&\le\bigg(1+{1\over\nu_2}\bigg)\bigg[{2d_3^2\over1-e^{-\nu_1T}}
e^{\nu_1T}+d_0^2\bigg]\equiv d_5^2,\cr}
\leqno(2.2)
$$
where $t\in(kT,(k+1)T)$, $k\in\N_0$.

\Proof
Multiplying (1.1) by $v$, integrating over $\Omega$, using $(1.1)_2$ and the
boundary conditions we obtain
$$
{1\over2}{d\over dt}\|v\|_{L_2(\Omega)}^2+\nu E_\Omega(v)=\intop_\Omega
f\cdot vdx,
\leqno(2.3)
$$
where
$$
E_\Omega(v)=\|\D(v)\|_{L_2(\Omega)}^2.
$$
From \cite[Ch. 4, Lemma 2.4]{Z4} we have the Korn inequality
$$
\|v\|_{H^1(\Omega)}^2\le c_k\bigg(E_\Omega(v)+
\bigg|\intop_\Omega v\cdot\eta dx\bigg|^2\bigg),
\leqno(2.4)
$$
where $\eta=(-x_2,x_1,0)$, $v\cdot\eta=rv_\varphi=u$.

\noindent
Now we calculate the last term on the r.h.s. of (2.4). Multiplying $(1.1)_1$
by $\eta$, integrating over $\Omega$ and using $(1.1)_2$, $(1.1)_3$ we obtain
$$\eqal{
&{d\over dt}\intop_\Omega v\cdot\eta dx-\intop_\Omega v_iv_j\nabla_i\eta_jdx+
\intop_\Omega\T_{ij}\nabla_i\eta_jdx\cr
&=\intop_\Omega f\cdot\eta dx,\cr}
\leqno(2.5)
$$
where $\T(v,p)=\{\T_{ij}\}_{i,j=1,2,3}$ is the stress tensor of the form
$$
\T(v,p)=\nu\D(v)-pI,
$$
$I$ is the unit matrix and the summation convention over the repeated indices
is assumed. Since $\nabla\eta$ is an antisymmetric tensor equation (2.5)
implies
$$
{d\over dt}\intop_\Omega v\cdot\eta dx=\intop_\Omega f\cdot\eta dx.
\leqno(2.6)
$$
Integrating (2.6) with respect to time from $kT$ to $t\in(kT,(k+1)T]$ yields
$$
\intop_\Omega u(t)dx=\intop_{kT}^t\intop_\Omega gdxdt'+\intop_\Omega u(kT)dx.
\leqno(2.7)
$$
Using (2.4) and (2.7) in (2.3) implies
$$\eqal{
&{1\over2}{d\over dt}\|v\|_{L_2(\Omega)}^2+\nu_*\|v\|_{H^1(\Omega)}^2\le
\intop_\Omega f\cdot vdx+\nu\bigg|\intop_{kT}^t\intop_\Omega gdxdt'\cr
&\quad+\intop_\Omega u(kT)dx\bigg|^2,\cr}
\leqno(2.8)
$$
where the inequality is considered in the time interval $[kT,(k+1)T]$.

\noindent
Applying the H\"older and the Young inequalities to the first term on the
r.h.s. of (2.8) and multiplying the result by 2 we derive
$$
{d\over dt}\|v\|_{L_2(\Omega)}^2+\nu_*\|v\|_{H^1(\Omega)}^2\le{1\over\nu_*}
\|f\|_{L_{6/5}(\Omega)}^2+2\nu d_1^2
\leqno(2.9)
$$
Employing the decomposition $\nu_*=\nu_1+\nu_2$ inequality (2.9) takes the form
$$\eqal{
&{d\over dt}(\\v\|_{L_2(\Omega)}^2e^{\nu_1t})+\nu_2\|v\|_{H^1(\Omega)}^2
e^{\nu_1t}\cr
&\le\bigg({1\over\nu_*}d_2^2+2\nu d_1^2\bigg)e^{\nu_1t}\equiv\nu_1d_3^2
e^{\nu_1t}.\cr}
\leqno(2.10)
$$
Omitting the second term on the l.h.s. of (2.10) and integrating the result
with respect to time from $kT$ to $t\in(kT,(k+1)T]$ yields
$$
\|v(t)\|_{L_2(\Omega)}^2\le d_3^2+\|v(kT)\|_{L_2(\Omega)}^2e^{-\nu_1(t-kT)}.
$$
By iteration we obtain
$$
\|v(kT)\|_{L_2(\Omega)}^2\le{d_3^2\over1-e^{-\nu_1T}}+
\|v(0)\|_{L_2(\Omega)}^2e^{-\nu_1kT}.
$$
Hence for $t\in(kT,(k+1)T]$ we get
$$
\|v(t)\|_{L_2(\Omega)}^2\le d_3^2{2-e^{-\nu_1T}\over1-e^{-\nu_1T}}+
\|v(0)\|_{L_2(\Omega)}^2e^{-\nu_1T}
\leqno(2.11)
$$
so (2.1) holds. Integrating (2.10) with respect to time from $kT$ to
$t\in(kT,(k+1)T]$ we have
$$\eqal{
&\|v(t)\|_{L_2(\Omega)}^2+\nu_2e^{-\nu_1t}\intop_{kT}^t
\|v(t')\|_{H^1(\Omega)}^2e^{\nu_1t'}dt'\cr
&\le d_3^2+\|v(kT)\|_{L_2(\Omega)}^2e^{-\nu_1(t-kT)}.\cr}
\leqno(2.12)
$$
Continuing, we obtain
$$\eqal{
&\|v(t)\|_{L_2(\Omega)}^2+\nu_2e^{-\nu_1(t-kT)}\intop_{kT}^t
\|v(t')\|_{H^1(\Omega)}^2dt'\cr
&\le d_3^2+\|v(kT)\|_{L_2(\Omega)}^2e^{-\nu_1(t-kT)}.\cr}
\leqno(2.13)
$$
Finally, (2.13) implies
$$
\intop_{kT}^t\|v(t')\|_{H^1(\Omega)}^2dt'\le{d_3^2\over\nu_2}e^{\nu_1T}+
{1\over\nu_2}\bigg({d_3^2\over1-e^{-\nu_1T}}+e^{-\nu_1kT}d_0^2\bigg).
\leqno(2.14)
$$
Combining (2.11) and (2.14) gives (2.2). This concludes the proof.

\proclaim Lemma 2.3.
Let the assumptions of Lemma 2.2 hold. Let $f\in L_2(\Omega\times(kT,t))$,
$t\in(kT,(k+1)T]$. For weak solutions to problem (1.3)--(1.8) the following
estimate holds
$$\eqal{
&\|v\|_{V_2^0(\Omega\times(kT,t))}^2+
\bigg\|{v_r\over r}\bigg\|_{L_2(\Omega\times(kT,t))}^2\cr
&\quad+\bigg\|{v_\varphi\over r}\bigg\|_{L_2(\Omega\times(kT,t))}^2\le c(T+1)
d_4^2+c\|f\|_{L_2(\Omega\times(kT,t))}^2\equiv d_7^2,\cr}
\leqno(2.15)
$$
where $t\in(kT,(k+1)T]$, $k\in\N_0$.

\Proof
Multiplying (1.3) by $v_r$, integrating over $\Omega$ and using boundary
conditions (1.7) yields
$$\eqal{
&{1\over2}{d\over dt}\intop_\Omega v_r^2dx-\intop_\Omega{v_\varphi^2\over r}
v_rdx+\nu\intop_\Omega(v_{r,r}^2+v_{r,z}^2)dx\cr
&\quad+\nu\intop_\Omega{v_r^2\over r^2}dx=-\intop_\Omega p_{,r}v_rdx+
\intop_\Omega f_rv_rdx.\cr}
\leqno(2.16)
$$
Multiplying (1.4) by $v_\varphi$, integrating over $\Omega$ and using boundary
conditions (1.7) implies
$$\eqal{
&{1\over2}{d\over dt}\intop_\Omega v_\varphi^2dx+\intop_\Omega{v_r\over r}
v_\varphi^2dx+\nu\intop_\Omega(v_{\varphi,r}^2+v_{\varphi,z}^2)dx-\nu
\intop_{-a}^av_{\varphi}^2dz\cr
&\quad+\nu\intop_\Omega{v_\varphi^2\over r^2}dx=\intop_\Omega f_\varphi
v_\varphi dx.\cr}
\leqno(2.17)
$$
Multiplying (1.5) by $v_z$, integrating over $\Omega$ and using the boundary
conditions (1.7) we obtain
$$
{1\over2}{d\over dt}\intop_\Omega v_z^2dx+\nu\intop_\Omega(v_{z,r}^2+v_{z,z}^2)
dx=-\intop_\Omega p_{,z}v_zdx+\intop_\Omega f_zv_zdx.
\leqno(2.18)
$$
Adding the above equations, using (1.6), the inequalities
$$\eqal{
&\|v_\varphi\|_{L_2(S_1)}^2\le\varepsilon\|\nabla v_\varphi\|_{L_2(\Omega)}^2
+c(1/\varepsilon)\|v_\varphi\|_{L_2(\Omega)}^2,\cr
&\intop_\Omega(v_r^2+v_\varphi^2)dx\le R^2\intop_\Omega\bigg({v_r^2\over r^2}+
{v_\varphi^2\over r^2}\bigg)dx,\cr}
$$
and the Poincare inequality for $v_z$ we arrive to the inequality
$$\eqal{
&{d\over dt}\intop_\Omega v^2dx+\nu\intop_\Omega|\nabla v|^2dx+\nu\intop_\Omega
\bigg({v_r^2\over r^2}+{v_\varphi^2\over r^2}
\bigg)dx\cr
&\le c\bigg(\intop_\Omega v_\varphi^2dx+\intop_\Omega f^2dx\bigg).\cr}
\leqno(2.19)
$$
Integrating (2.19) with respect to time from $kT$ to $t$ and using (2.1) we
obtain
$$\eqal{
&\|v\|_{V_2^0(\Omega\times(kT,t))}^2+\nu\intop_{kT}^t\intop_\Omega
\bigg({v_r^2\over r^2}+{v_\varphi^2\over r^2}\bigg)dxdt'\cr
&\le cTd_4^2+c\|f\|_{L_2(\Omega\times(kT,t))}^2+\|v(kT)\|_{L_2(\Omega)}^2.\cr}
\leqno(2.20)
$$
Using again (2.1) we obtain (2.15). This concludes the proof.

Let us consider the Stokes problem
$$\eqal{
&v_t-\divv\T(v,p)=f\quad &{\rm in}\ \ \Omega^T,\cr
&\divv v=0\quad &{\rm in}\ \ \Omega^T,\cr
&v\cdot\bar n=0,\ \ \bar n\cdot\D(v)\cdot\bar\tau_\alpha=0,\ \ \alpha=1,2,\quad
&{\rm on}\ \ S^T,\cr
&v|_{t=0}=v_0\quad &{\rm in}\ \ \Omega.\cr}
\leqno(2.21)
$$

\proclaim Lemma 2.4. (see \cite{S, Z7, ZZ})
Assume that $f\in L_s(\Omega^T)$, $v_0\in W_s^{2-2/s}(\Omega)$,
$s\in(1,\infty)$, $S_1\in C^2$. Then there exists a solution to problem (2.21)
such that $v\in W_s^{2,1}(\Omega^T)$, $\nabla p\in L_s(\Omega^T)$ and
there exists a constant $c_0=c(\Omega,s)$ such that
$$
\|v\|_{W_s^{2,1}(\Omega^T)}+\|\nabla p\|_{L_s(\Omega^T)}\le c(\Omega,s)
(\|f\|_{L_s(\Omega^T)}+\|v_0\|_{W_s^{2-2/s}(\Omega)}).
\leqno(2.22)
$$

\noindent
From (1.4) and (1.7) we obtain the following problem for swirl $u=rv_\varphi$
$$\eqal{
&u_t+v\cdot\nabla u-\nu\Delta u+2\nu{u_{,r}\over r}=rf_\varphi\equiv g\quad
&{\rm in}\ \ \Omega^T,\cr
&u_{,r}={2\over r}u\quad &{\rm on}\ \ S_1^T,\cr
&u_{,z}=0\quad &{\rm on}\ \ S_2^T,\cr
&u|_{t=0}=u_0\quad &{\rm in}\ \ \Omega,\cr}
\leqno(2.23)
$$
where $v$ is the divergence free vector.

Let us introduce the $\varphi$-component of vorticity by
$$
\chi=v_{r,z}-v_{z,r}.
\leqno(2.24)
$$
In view of (1.3), (1.5) and (1.7), $\chi$ is a solution to the problem (see
\cite[Ch. 4, Sect. 3 (3.1) and Lemma 2.2 (2.13)]{Z4})
$$\eqal{
&\chi_t+v\cdot\nabla\chi-{v_r\over r}\chi-\nu\bigg[\bigg(r\bigg(
{\chi\over r}\bigg)_{,r}\bigg)_{,r}+\chi_{,zz}+2\bigg({\chi\over r}\bigg)_{,r}
\bigg]\quad\cr
&\quad={2v_\varphi v_{\varphi,z}\over r}+F\quad &{\rm in}\ \ \Omega^T,\cr
&\chi|_S=0\quad &{\rm on}\ \ S^T,\cr
&\chi|_{t=0}=\chi_0\quad &{\rm in}\ \ \Omega,\cr}
\leqno(2.25)
$$
where $F=f_{r,z}-f_{z,r}$.
To prove global regular solutions to problem (1.1) we need a priori estimate
with weights which are singular on the axis of symmetry (see \cite{Z1}).
Therefore it is convenient to consider instead of problem (1.1) the following
approximated problem
$$\eqal{
&v_t+v\cdot\nabla v-\nu\Delta v+\nabla p=f\quad &{\rm in}\ \
\Omega_\varepsilon^T,\cr
&\divv v=0\quad &{\rm in}\ \ \Omega_\varepsilon^T,\cr
&\bar n\cdot v=0,\ \ \bar n\cdot\D(v)\cdot\bar\tau_\alpha=0,\ \ \alpha=1,2,
\quad &{\rm on}\ \ S_1^T\cup S_\varepsilon^T\cup S_{2\varepsilon}^T,\cr
&v|_{t=0}=v_0\quad &{\rm in}\ \ \Omega_\varepsilon,\cr}
\leqno(2.26)
$$
where $\Omega_\varepsilon=\{x\in\Omega:\ r>\varepsilon\}$,
$S_\varepsilon=\{x\in\R^3:\ r=\varepsilon,\ |z|<a\}$,
$S_{2\varepsilon}=\{x\in S_2:\ r>\varepsilon\}$.
Setting $\varepsilon=0$ we obtain problem (1.1).

\proclaim Lemma 2.5.
Assume that $v(0)\in H^1(\Omega_\varepsilon)$,
$f\in L_2(\Omega_\varepsilon^T)$, $\varepsilon\ge0$, $T>0$. Assume that $T$
is so small that
$$
\varphi(c_0,c_1,c_2)T^{1/2}d_4[\|f\|_{L_2(\Omega_\varepsilon^T)}^3+
\|v(0)\|_{H^1(\Omega_\varepsilon)}^3+T^{3/2}d_4^3+1]\le1/2,
\leqno(2.27)
$$
where $\varphi$ is some positive increasing function of its arguments and
$d_4$ is introduced in (2.1).
Then there exists a solution to problem (2.26) such that
$v\in W_2^{2,1}(\Omega_\varepsilon^T)$,
$\nabla p\in L_2(\Omega_\varepsilon^T)$ and
$$
\|v\|_{W_2^{2,1}(\Omega_\varepsilon^T)}+
\|\nabla p\|_{L_2(\Omega_\varepsilon^T)}\le8c_0
(\|f\|_{L_2(\Omega_\varepsilon^T)}+\|v(0)\|_{H^1(\Omega_\varepsilon)}),
\leqno(2.28)
$$
where $c_0$ appears in (2.22), $c_1$ in (2.30) and $c_2$ in (2.31).
The existence and the estimate hold also for $\varepsilon=0$.

\Proof
Since we are going to prove the existence of solutions to problem (2.26) by
the Leray-Schauder fixed point theorem we restrict the proof to show a priori
estimate (2.28) only because other steps of it are clear. In the proof we
omit the index $\varepsilon$ for simplicity.

\noindent
Applying Lemma 2.4 with $s=2$ to problem (2.26) yields
$$\eqal{
&\|v\|_{W_2^{2,1}(\Omega^T)}+\|\nabla p\|_{L_2(\Omega^T)}\le c_0
(\|v\cdot\nabla v\|_{L_2(\Omega^T)}\cr
&\quad+\|f\|_{L_2(\Omega^T)}+\|v(0)\|_{H^1(\Omega)}).\cr}
\leqno(2.29)
$$
Now we examine the first term on the r.h.s. of (2.29). We estimate it by
$$\eqal{
&\bigg(\intop_0^Tdt\intop_\Omega|v\cdot\nabla v|^2dx\bigg)^{1/2}\le
\bigg(\intop_0^T\|v(t)\|_{L_\infty(\Omega)}^2
\|\nabla v(t)\|_{L_2(\Omega)}^2dt\bigg)^{1/2}\cr
&\le\sup_t\|\nabla v(t)\|_{L_2(\Omega)}\bigg(\intop_0^T
\|v(t)\|_{L_\infty(\Omega)}^2dt\bigg)^{1/2}\equiv I_1.\cr}
$$
From \cite[Ch. 3, Sect. 15]{BIN} we have the interpolation
$$
\|v\|_{L_\infty(\Omega)}\le c_1\|v_{xx}\|_{L_2(\Omega)}^{3/4}
\|v\|_{L_2(\Omega)}^{1/4}+c_1\|v\|_{L_2(\Omega)}.
\leqno(2.30)
$$
Employing (2.30) in $I_1$ yields
$$\eqal{
I_1&\le c_1\sup_t\|\nabla v(t)\|_{L_2(\Omega)}(T^{1/8}
\|v_{xx}\|_{L_2(\Omega^T)}^{3/4}\sup_t\|v(t)\|_{L_2(\Omega)}^{1/4}\cr
&\quad+T^{1/2}\sup_t\|v(t)\|_{L_2(\Omega)})\equiv I_2.\cr}
$$
Using the estimate
$$\eqal{
&\sup_t\|\nabla  v(t)\|_{L_2(\Omega)}\le\sup_t\|v(t)\|_{H^1(\Omega)}\cr
&\le c_2(\|v\|_{W_2^{2,1}|(\Omega^T)}+\|v(0)\|_{H^1(\Omega)})\cr}
\leqno(2.31)
$$
and the energy estimate (2.1) in $I_2$ we obtain from (2.29) the inequality
$$\eqal{
&\|v\|_{W_2^{2,1}(\Omega^T)}+\|\nabla p\|_{L_2(\Omega^T)}\le c_0c_1c_2
(\|v\|_{W_2^{2,1}(\Omega^T)}\cr
&\quad+\|v(0)\|_{H^1(\Omega)})(T^{1/8}\|v\|_{W_2^{2,1}(\Omega^T)}^{3/4}
d_4^{1/4}+T^{1/2}d_4)\cr
&\quad+c_0(\|f\|_{L_2(\Omega^T)}+\|v(0)\|_{H^1(\Omega)}).\cr}
\leqno(2.32)
$$
Assuming that $T$ is so small that
$$
c_0c_1c_2(T^{1/8}\|v\|_{W_2^{2,1}(\Omega^T)}^{3/4}d_4^{1/4}+T^{1/2}d_4)\le
{1\over2}
\leqno(2.33)
$$
we derive from (2.32) the inequality
$$\eqal{
&\|v\|_{W_2^{2,1}(\Omega^T)}+\|\nabla p\|_{L_2(\Omega^T)}\le2c_0c_1c_2\cdot\cr
&\quad\cdot(T^{1/8}\|v\|_{W_2^{2,1}(\Omega^T)}^{3/4}d_4^{1/4}+T^{1/2}d_4)
\|v(0)\|_{H^1(\Omega)}\cr
&\quad+2c_0(\|f\|_{L_2(\Omega)}+\|v(0)\|_{H^1(\Omega)}).\cr}
\leqno(2.34)
$$
The condition (2.33) is not written in a final form because it contains an
unknown norm $\|v\|_{W_2^{2,1}(\Omega^T)}$. Applying the Young inequality to
the first expression on the r.h.s. of (2.34) yields
$$\eqal{
&\|v\|_{W_2^{2,1}(\Omega^T)}+\|\nabla p\|_{L_2(\Omega^T)}\le
{\varepsilon^{4/3}\over4/3}\|v\|_{W_2^{2,1}(\Omega^T)}\cr
&\quad+{1\over4\varepsilon^4}(2c_0c_1c_2T^{1/8}d_4\|v(0)\|_{H^1(\Omega)})^4+
2c_0c_1c_2T^{1/2}d_4\|v(0)\|_{H^1(\Omega)}\cr
&\quad+2c_0(\|f\|_{L_2(\Omega^T)}+\|v(0)\|_{H^1(\Omega)}).\cr}
$$
Setting ${\varepsilon^{4/3}\over4/3}={1\over2}$ we obtain that
$\varepsilon=\big({2\over3}\big)^{3/4}$ so the above inequality yields
$$\eqal{
&\|v\|_{W_2^{2,1}(\Omega^T)}+\|\nabla p\|_{L_2(\Omega^T)}\le27(c_0c_1c_2)^4
T^{1/2}d_4\|v(0)\|_{H^1(\Omega)}^4\cr
&\quad+4c_0c_1c_2T^{1/2}d_4\|v(0)\|_{H^1(\Omega)}+4c_0
(\|f\|_{L_2(\Omega^T)}+\|v(0)\|_{H^1(\Omega)}).\cr}
\leqno(2.35)
$$
Assuming that $T$ is so small that
$$
27(c_0c_1c_2)^4T^{1/2}d_4\|v(0)\|_{H^1(\Omega)}^3+4c_0c_1c_2T^{1/2}d_4\le4c_0
\leqno(2.36)
$$
we obtain that (2.35) implies (2.28). Using (2.28) we express condition (2.33)
in the form
$$\eqal{
&8(c_0c_1c_2)^4T^{1/2}d_4[(8c_0)^3(\|f\|_{L_2(\Omega^T)}^3+
\|v(0)\|_{H^1(\Omega)}^3)\cr
&\quad+(T^{1/2}d_4)^3]\le{1\over2}\cr}
\leqno(2.37)
$$
There exists a function $\varphi(c_0,c_1,c_2)$ such that (2.36) and (2.37)
imply (2.27). This concludes the proof.

\section{3. Regularity of swirl $u$}

Let us recall that $u$ is a solution to problem (2.23). To show $L_\infty$
estimate of $u$ we need some notation
$$\eqal{
&A_k(t)=\{x\in\Omega:\ u(x,t)>k\},\quad u^{(k)}=\max\{u-k,0\}\cr
&\mu(k)=\intop_0^T\meas^{r\over q}A_k(t)dt,\cr}
$$
where ${3\over q}+{2\over r}={3\over2}$.

\proclaim Lemma 3.1.
Assume that $u_0\in L_\infty(\Omega)$, $g\in L_\infty(\Omega^T)$. Then there
exists a constant $d_6$ depending on $\|u_0\|_{L_\infty(\Omega)}$,
$\|g\|_{L_\infty(\Omega^T)}$ such that
$$
\|u\|_{L_\infty(\Omega^T)}\le d_6.
\leqno(3.1)
$$

\Proof
Multiplying $(2.23)_1$ by $u^{(k)}$ and integrating the result over $\Omega$
yields
$$\eqal{
&{1\over2}{d\over dt}\intop_\Omega|u^{(k)}|^2dx+\nu\intop_\Omega
|\nabla u^{(k)}|^2dx-\nu\intop_{-a}^a|u^{(k)}|_{r=R}|^2dz\cr
&\quad+2\nu\intop_\Omega u_{,r}u^{(k)}drdz=\intop_\Omega gu^{(k)}dx.\cr}
\leqno(3.2)
$$
The last term on the l.h.s. of (3.2) equals
$$
\nu\intop_\Omega(|u^{(k)}|^2)_{,r}drdz=\nu\intop_{-a}^a|u^{(k)}|_{r=R}|^2dz
$$
because $u^{(k)}|_{r=0}=0$. Otherwise the condition
$\intop_0^T\intop_\Omega{u^2\over r^4}dxdt<\infty$ (see (2.15)) implies
a contradition for $u>k$. Hence, (3.2) takes the form
$$
{1\over2}{d\over dt}\intop_\Omega|u^{(k)}|^2dx+\nu\intop_\Omega
|\nabla u^{(k)}|^2dx=\intop_\Omega gu^{(k)}dx.
\leqno(3.3)
$$
Integrating (3.3) with respect to time and using that
$k>\|u_0\|_{L_\infty(\Omega)}$ we obtain
$$\eqal{
&\|u^{(k)}\|_{V_2^0(\Omega^T)}^2\le c\intop_{\Omega^T}|gu^{(k)}|dxdt\cr
&\le c\|g\|_{L_\infty(\Omega^T)}\|u^{(k)}\|_{L_{10/3}(\Omega^T)}
(\mu(k))^{7/10}\cr}
\leqno(3.4)
$$
In view of imbedding $V_2^0(\Omega^T)\subset L_{10/3}(\Omega^T)$ we have
$$
\|u^{(k)}\|_{V_2^0(\Omega^T)}\le c\|g\|_{L_\infty(\Omega^T)}
|\mu(k)|^{{3\over10}(1+\varkappa)},
\leqno(3.5)
$$
where ${7\over10}={3\over10}(1+\varkappa)$ so ${7\over3}=1+\varkappa$,
$\varkappa={4\over3}$ and $\mu(k)$ is calculated for $r=q={10\over3}$.
From \cite[Ch. 2, Sect. 6, Theorem 6.1]{LSU} we conclude the proof.

\noindent
Repeating the considerations from \cite{Z3} and \cite[Ch. 2, Sect. 8]{LSU}
we have

\proclaim Lemma 3.2. (H\"older continuity of swirl).
Let $g\in C^{\alpha,\alpha/2}(\Omega^T)$, $u_0\in C^\alpha(\Omega)$,
$\alpha\in(0,1/2]$. Let $v\in L_{10}(\Omega^T)$.
Then $u\in C^{\alpha,\alpha/2}(\Omega^T)$.

\noindent
For more details see \cite{B}.

\proclaim Lemma 3.3.
Let the assumptions of Lemmas 3.1 and 2.3 hold. Then
$$
\|v_\varphi\|_{L_4(\Omega^t)}^4\le d_6^2d_7^2,\quad t\le T,
\leqno(3.6)
$$
where $d_6$ is introduced by (3.1) and $d_7$ by (2.15).

\section{4. A priori estimates for $\chi$ in a neighborhood of the $x_3$-axis}

First we examine the elliptic problem
$$\eqal{
&v_{r,z}-v_{z,r}=\chi\quad &{\rm in}\ \ \Omega,\cr
&v_{r,r}+v_{z,z}+{v_r\over r}=0\quad &{\rm in}\ \ \Omega,\cr
&v_r|_{S_1}=0,\ \ v_z|_{S_2}=0.\cr}
\leqno(4.1)
$$
Expressing $(4.1)_2$ in the form
$$
(rv_r)_{,r}+(rv_z)_{,z}=0
\leqno(4.2)
$$
we have existence of a potential $\psi=\psi(r,z,t)$ such that
$$
v_r={\psi_{,z}\over r},\quad v_z=-{\psi_{,r}\over r}.
\leqno(4.3)
$$
Then the boundary conditions $(4.1)_3$ are satisfied if
$$
\psi|_{S_1}=0\ \ {\rm so}\ \psi(R,z,t)=0\ \ {\rm and}\ \
\psi|_{S_2}=0\ \ {\rm so}\ \ \psi(r,-a,t)=\psi(r,a,t)=0.
\leqno(4.4)
$$
Therefore, problem (4.1) takes the form
$$\eqal{
&\bigg({\psi_{,z}\over r}\bigg)_{,z}+\bigg({\psi_{,r}\over r}\bigg)_{,r}=\chi
\quad &{\rm in}\ \ \Omega,\cr
&\psi|_S=0\quad &{\rm on}\ \ S.\cr}
\leqno(4.5)
$$
To obtain an a priori estimate for $v'=(v_r,v_z)$ near the $x_3$-axis we have
to work with problems with singular coefficients on it. Therefore, we can work
either with functions vanishing sufficiently fast near the axis of symmetry
or examining approximate solutions described by (2.26). We shall restrict
our considerations to the second case because it seems to be more precise.
Then problem (2.25) for $\chi$ takes the form
$$\eqal{
&\chi_{,t}+v\cdot\nabla\chi-{v_r\over r}\chi-\nu\bigg[\bigg(r
\bigg({\chi\over r}\bigg)_{,r}\bigg)_{,r}+\chi_{,zz}+
2\bigg({\chi\over r}\bigg)_{,r}\bigg]\hskip-20pt\cr
&\quad={2v_\varphi v_{\varphi,z}\over r}+F\quad &{\rm in}\ \
\Omega_\varepsilon^T,\cr
&\chi=0\quad &{\rm on}\ \ S_1^T\cup S_\varepsilon^T\cup S_{2\varepsilon}^T,\cr
&\chi|_{t=0}=\chi_0\quad &{\rm in}\ \ \Omega_\varepsilon,\cr}
\leqno(4.6)
$$
where $\Omega_\varepsilon$, $S_\varepsilon$, $_{2\varepsilon}$ are defined
below problem (2.26).

\noindent
The condition $\chi|_{S_\varepsilon}=0$ follows from \cite[Lemma 2.2]{Z4}.
Hence we assume that
$$
\chi=0\quad {\rm for}\quad r\le\varepsilon.
\leqno(4.7)
$$
Solutions of (2.26) and (4.6) should be labeled with index $\varepsilon$
which we omit for simplicity.
To examine a behavior of solutions in a neighborhood of the axis of symmetry
we introduce a smooth cut-off function $\zeta_1=\zeta_1(r)$ such that
$\zeta_1(r)=1$ for $r\le r_0$ and $\zeta_1(r)=0$ for $r\ge2r_0$, where
$2r_0<R$.

\noindent
Let us introduce the notation
$$\eqal{
&\tilde\chi=\chi\zeta_1^2,\quad \tilde v_\varphi=v_\varphi\zeta_1,\quad
\tilde v'=v'\zeta_1^2,\quad \tilde F=F\zeta_1^2,\cr
&\tilde f'=f'\zeta_1^2,\quad \tilde f_\varphi=f_\varphi\zeta_1.\cr}
\leqno(4.8)
$$
Then $\tilde\chi$ is a solution to the problem
$$\eqal{
&\tilde\chi_{,t}+v\cdot\nabla\tilde\chi-{v_r\over r}\tilde\chi-\nu
\bigg[\bigg(r\bigg({\tilde\chi\over r}\bigg)_{,r}\bigg)_{,r}+\tilde\chi_{,zz}+
2\bigg({\tilde\chi\over r}\bigg)_{,r}\bigg]\quad\cr
&=v\cdot\nabla\zeta_1^2\chi-\nu(\chi\zeta_{1,r}^2)_{,r}-\nu r
\bigg({\chi\over r}\bigg)_{,r}\zeta_{1,r}^2-2\nu\bigg({\chi\over r}
\zeta_{1,r}^2\bigg)\cr
&\quad+{2\tilde v_\varphi\tilde v_{\varphi,z}\over r}+\tilde F\quad
&{\rm in}\ \ \Omega_\varepsilon^T,\cr
&\tilde\chi|_{r=2r_0}=0,\ \ \tilde\chi|_{r=\varepsilon}=0,\ \
\tilde\chi|_{S_2}=0,\cr
&\tilde\chi|_{t=0}=\tilde\chi_0\quad &{\rm in}\ \ \Omega_\varepsilon.\cr}
\leqno(4.9)
$$

\proclaim Lemma 4.1.
Let $\Omega_\varepsilon=\{x\in\Omega:\ 0<\varepsilon<r\}$. Assume that there
exists a weak solution to problem (1.1) described by Lemma 2.2. Assume that
${\tilde v_\varphi\over r}\in L_4(\Omega_\varepsilon^T)$,
$\chi\in L_{20\over7}(\Omega\cap\supp\zeta_{1,r}\times(0,T))$, where
$\supp\zeta_{1,r}=\{x\in\Omega_\varepsilon:\ r_0\le r\le2r_0\}$ and
${\tilde F\over r}\in L_2(\Omega_\varepsilon^T)$. Then for sufficiently
smooth solutions to problem (1.1) we have
$$\eqal{
&\bigg\|{\tilde\chi\over r}\bigg\|_{L_\infty(0,t;L_2(\Omega_\varepsilon))}^2+
\nu\bigg\|\nabla{\tilde\chi\over r}\bigg\|_{L_2(\Omega_\varepsilon^t)}^2\cr
&\le{2\over\nu}
\bigg\|{\tilde v_\varphi\over r}\bigg\|_{L_4(\Omega_\varepsilon^t)}^4+
c(1/r_0)d_5^2+c(1/r_0)d_5
\|\chi\|_{L_{20\over7}(\Omega_{\varepsilon,\zeta_{1,r}}^t)}^2\cr
&\quad+\bigg\|{\tilde\chi(0)\over r}\bigg\|_{L_2(\Omega_\varepsilon)}^2+c
\bigg\|{\tilde F\over r}\bigg\|_{L_2(\Omega_\varepsilon^t)},\quad t\le T,\cr}
\leqno(4.10)
$$
where $\varepsilon<r_0$ and $d_5$ is introduced in (2.2).

\Proof
Multiplying $(4.9)_1$ by ${\tilde\chi\over r^2}$ and integrating the result
over $\Omega_\varepsilon$ we obtain
$$\eqal{
&{1\over2}{d\over dt}
\bigg\|{\tilde\chi\over r}\bigg\|_{L_2(\Omega_\varepsilon)}^2+\nu
\bigg\|\nabla{\tilde\chi\over r}\bigg\|_{L_2(\Omega_\varepsilon)}^2\cr
&=\intop_{\Omega_\varepsilon}\bigg[v\cdot\nabla\zeta_1^2\chi
{\tilde\chi\over r^2}-\nu(\chi\zeta_{1,r}^2)_{,r}{\tilde\chi\over r^2}\cr
&\quad-\nu r\bigg({\chi\over r}\bigg)_{,r}\zeta_{1,r}^2{\tilde\chi\over r^2}-
2\nu{\chi\over r}\zeta_{1,r}^2{\tilde\chi\over r^2}\bigg]dx\cr
&\quad+2\intop_{\Omega_\varepsilon}
{\tilde v_\varphi\tilde v_{\varphi,z}\over r}{\tilde\chi\over r^2}dx+
\intop_{\Omega_\varepsilon}\tilde F{\tilde\chi\over r^2}dx.\cr}
\leqno(4.11)
$$
Now we estimate the particular terms on the r.h.s. of (4.11). The second term
implies
$$\eqal{
&\intop_{\Omega_\varepsilon}{(\tilde v_\varphi^2)_{,z}\over r^2}
{\tilde\chi\over r}dx=-\intop_{\Omega_\varepsilon}
{\tilde v_\varphi^2\over r^2}\bigg({\tilde\chi\over r}\bigg)_{,z}dx\cr
&\le{\varepsilon_1\over2}\intop_{\Omega_\varepsilon}
\bigg({\tilde\chi\over r}\bigg)_{,z}^2dx+{1\over2\varepsilon_1}
\intop_{\Omega_\varepsilon}\bigg({\tilde v_\varphi\over r}\bigg)^4dx.\cr}
$$
To estimate the first term on the r.h.s. of (4.11) we use properties of the
cut-off function $\zeta_1=\zeta_1(r)$.
The first expression under the square bracket is bounded by
$$
 c(1/r_0)\intop_{\Omega_{\varepsilon,\zeta_{1,r}}}|v_r|\chi^2dx.
$$
The fourth term under the square bracket is estimated by
$$
c(1/r_0)\intop_{\Omega_{\varepsilon,\zeta_{1,r}}}\chi^2dx.
$$
The second term under the square bracket equals
$$
-\nu\intop_{\Omega_\varepsilon}\chi_{,r}\zeta_{1,r}^2\zeta_1^2{\chi\over r^2}dx
-\nu\intop_{\Omega_\varepsilon}\chi\zeta_{1,rr}^2\zeta_1^2{\chi\over r^2}dx
\equiv I_1.
$$
Integrating by parts in the first integral in $I_1$ it takes the form
$$
-{\nu\over2}\intop_{\Omega_\varepsilon}(\chi^2)_{,r}{1\over r^2}\zeta_{1,r}^2
\zeta_1^2dx={\nu\over2}\intop_{\Omega_\varepsilon}\chi^2
\bigg(\zeta_{1,r}^2\zeta_1^2{1\over r}\bigg)_{,r}drdz.
$$
Hence
$$
|I_1|\le c(1/r_0)\intop_{\Omega_{\varepsilon,\zeta_{1,r}}}\chi^2dx
$$
Finally, the third term under the square bracket yields
$$
-{\nu\over2}\intop_{\Omega_\varepsilon}
\bigg[\bigg({\chi\over r}\bigg)^2\bigg]_{,r}\zeta_{1,r}^2\zeta_1^2dx=
{\nu\over2}\intop_{\Omega_\varepsilon}\bigg({\chi\over r}\bigg)^2
(\zeta_{1,r}^2\zeta_1^2r)_{,r}drdz\equiv I_2.
$$
Then,
$$
|I_2|\le c(1/r_0)\intop_{\Omega_{\varepsilon,\zeta_{1,r}}}\chi^2dx.
$$
Applying the Cauchy inequalities and the Poincare inequality to the last term
on the r.h.s. of (4.11) we estimate it by
$$
{\varepsilon_2\over 2}\intop_{\Omega_\varepsilon}
\bigg|\nabla{\tilde\chi\over r}\bigg|^2dx+{1\over2\varepsilon_2}
\intop_{\Omega_\varepsilon}\bigg|{\tilde F\over r}\bigg|^2dx.
$$
Using the above estimates in (4.11) and assuming that
$\varepsilon_1=\varepsilon_2={\nu\over2}$ we obtain from (4.11) the inequality
$$\eqal{
&{1\over2}{d\over dt}
\bigg\|{\tilde\chi\over r}\bigg\|_{L_2(\Omega_\varepsilon)}^2+{\nu\over2}
\bigg\|\nabla{\tilde\chi\over r}\bigg\|_{L_2(\Omega_\varepsilon)}^2\cr
&\le{1\over\nu}
\bigg\|{\tilde v_\varphi\over r}\bigg\|_{L_4(\Omega_\varepsilon)}^4+
c(1/r_0)\intop_{\Omega_{\varepsilon,\zeta_{1,r}}}|v_r|\chi^2dx\cr
&\quad+c(1/r_0)\intop_{\Omega_{\varepsilon,\zeta_{1,r}}}\chi^2dx+c
\intop_{\Omega_\varepsilon}\bigg|{\tilde F\over r}\bigg|^2dx.\cr}
\leqno(4.12)
$$
Integrating (4.12) with respect to time and using estimate (2.2) to the third
term on the r.h.s. of (4.12) we have
$$\eqal{
&\bigg\|{\tilde\chi\over r}\bigg\|_{L_\infty(0,t;L_2(\Omega_\varepsilon))}^2+
\nu\bigg\|\nabla{\tilde\chi\over r}\bigg\|_{L_2(\Omega_\varepsilon^t)}^2\cr
&\le{2\over\nu}
\bigg\|{\tilde v_\varphi\over r}\bigg\|_{L_4(\Omega_\varepsilon^t)}^4+
c(1/r_0)\intop_{\Omega_{\varepsilon,\zeta_{1,r}}^t}|v_r|\chi^2dxdt'\cr
&\quad+c(1/r_0)d_5^2+
\bigg\|{\tilde F\over r}\bigg\|_{L_2(\Omega_\varepsilon^t)}^2+
\bigg\|{\tilde\chi(0)\over r}\bigg\|_{L_2(\Omega_\varepsilon)}^2,\quad t\le T.
\cr}
\leqno(4.13)
$$
Applying the H\"older inequality and using estimate (2.2) in the second term
on the r.h.s. of (4.13) we estimate it by
$$\eqal{
&c(1/r_0)\|v\|_{L_{10/3}(\Omega_\varepsilon^t)}
\|\chi\|_{L_{20/7}(\Omega_{\varepsilon,\zeta_{1,r}}^t)}^2\cr
&\le c(1/r_0)d_5
\|\chi\|_{L_{20/7}(\Omega_{\varepsilon,\zeta_{1,r}}^t)}.\cr}
$$
Using the estimate in (4.13) we obtain (4.10). This concludes the proof.

A crucial step in the proof of an a priori estimate in a neighborhood of the
axis of symmetry is an estimate for $v_\varphi$ (see \cite{Z1}). Therefore
we localize equation (1.4) with boundary and initial conditions included in
(1.7) and (1.8). Hence we have the following problem
$$\eqal{
&\tilde v_{\varphi,t}+v\cdot\nabla\tilde v_\varphi+{v_r\over r}
\tilde v_\varphi-\nu\Delta\tilde v_\varphi+\nu{\tilde v_\varphi\over r^2}=
v\cdot\nabla\zeta_1v_\varphi\cr
&\quad-2\nu\nabla v_\varphi\nabla\zeta_1-\nu v_\varphi\Delta\zeta_1+
\tilde f_\varphi\quad {\rm in}\ \ \Omega_\varepsilon^T,\cr
&\tilde v_\varphi|_{t=0}=\tilde v_\varphi(0),\cr
&\tilde v_\varphi|_{r=2r_0}=0,\ \ \tilde v_\varphi|_{r=\varepsilon}=0,\ \
\tilde v_{\varphi,z}|_{S_2}=0.\cr}
\leqno(4.14)
$$

\proclaim Lemma 4.2.
Assume that $\tilde f_\varphi/\sqrt{r}\in L_{20/11}(\Omega_\varepsilon^t)$,
${\tilde v_\varphi(0)\over\sqrt{r}}\in L_4(\Omega_\varepsilon)$,
$\varepsilon>0$. Assume that estimates (2.2) and (3.1) hold, $d_5$ appears
in (2.2) and $d_6$ in (3.1). Assume that $v$ is a solution of (1.1). Then
for solutions of (4.14) the inequality is valid
$$\eqal{
&{1\over8}\intop_{\Omega_\varepsilon}{\tilde v_\varphi^4(t)\over r^2}dx+
{3\over8}\nu\intop_{\Omega_\varepsilon^t}
\bigg|\nabla{\tilde v_\varphi^2\over r}\bigg|^2dxdt'+{3\over4}\nu
\intop_{\Omega_\varepsilon^t}{\tilde v_\varphi^4\over r^4}dxdt'\cr
&\le{3\over2}\intop_{\Omega_\varepsilon^t}{|v_r|\over r}
{\tilde v_\varphi^4\over r^2}dxdt'+c(1/r_0)d_5^2(1+d_6)d_6^2+c\bigg\|
{\tilde f_\varphi\over r^{1/2}}\bigg\|_{L_{20/11}(\Omega_\varepsilon^t)}^4\cr
&\quad+{1\over4}\intop_{\Omega_\varepsilon}{\tilde v_\varphi^4(0)\over r^2}dx,
\quad t\le T.\cr}
\leqno(4.15)
$$

\Proof
Multiplying $(4.14)_1$ by ${\tilde v_\varphi|\tilde v_\varphi|^2\over r^2}$
and integrating over $\Omega_\varepsilon$ we obtain
$$\eqal{
&\intop_{\Omega_\varepsilon}\partial_t\tilde v_\varphi
{\tilde v_\varphi|\tilde v_\varphi|^2\over r^2}dx+\intop_{\Omega_\varepsilon}
v\cdot\nabla\tilde v_\varphi
{\tilde v_\varphi|\tilde v_\varphi|^2\over r^2}dx+
\intop_{\Omega_\varepsilon}{v_r\over r}{\tilde v_\varphi^4\over r^2}dx\cr
&\quad-\nu\intop_{\Omega_\varepsilon}\Delta\tilde v_\varphi
{\tilde v_\varphi|\tilde v_\varphi|^2\over r^2}dx+\nu
\intop_{\Omega_\varepsilon}{|\tilde v_\varphi|^4\over r^4}dx\cr
&=\intop_{\Omega_\varepsilon}v\cdot\nabla\zeta_1v_\varphi
{\tilde v_\varphi|\tilde v_\varphi|^2\over r^2}dx-\nu
\intop_{\Omega_\varepsilon}[2\nabla v_\varphi\nabla\zeta_1+v_\varphi
\Delta\zeta_1]{\tilde v_\varphi|\tilde v_\varphi|^2\over r^2}dx\cr
&\quad+\intop_{\Omega_\varepsilon}\tilde f_\varphi\tilde v_\varphi
{|\tilde v_\varphi|^2\over r^2}dx.\cr}
\leqno(4.16)
$$
Now we examine the particular terms in (4.16).
The first term on the l.h.s. of (4.16) equals
$$
{1\over4}{d\over dt}\intop_{\Omega_\varepsilon}{\tilde v_\varphi^4\over r^2}dx.
$$
The second term on the l.h.s. of (4.16) assumes the form
$$\eqal{
&{1\over4}\intop_{\Omega_\varepsilon}v\cdot\nabla(\tilde v_\varphi^4)
{1\over r^2}dx={1\over4}\intop_{\Omega_\varepsilon}\bigg[v\cdot\nabla
{\tilde v_\varphi^4\over r^2}-v\cdot\nabla{1\over r^2}\tilde v_\varphi^4\bigg]
dx
={1\over2}\intop_{\Omega_\varepsilon}{v_r\over r}
{\tilde v_\varphi^4\over r^2}dx,\cr}
$$
where the first integral vanishes after integration by parts because
$\divv v=0$ and $v\cdot\bar n|_{S_{2\varepsilon}}=0$, $\tilde v_\varphi=0$
on $r=2r_0$ and $r=\varepsilon$.
Hence, the sum of the second and the third terms on the l.h.s. of (4.16) equals
$$
{3\over2}\intop_{\Omega_\varepsilon}{v_r\over r}
{\tilde v_\varphi^4\over r^2}dx.
$$
Integrating by parts and using the boundary conditions in the fourth term on
the l.h.s. of (4.16) yields
$$\eqal{
I_1&=\nu\intop_{\Omega_\varepsilon}\nabla\tilde v_\varphi\cdot\nabla
\bigg(\tilde v_\varphi{\tilde v_\varphi^2\over r^2}\bigg)dx=3\nu
\intop_{\Omega_\varepsilon}|\nabla\tilde v_\varphi|^2
{\tilde v_\varphi^2\over r^2}dx
-2\nu\intop_{\Omega_\varepsilon}\tilde v_{\varphi,r}
{\tilde v_\varphi^3\over r^2}drdz\cr
&=3\nu\intop_{\Omega_\varepsilon}\bigg|\nabla\tilde v_\varphi
{\tilde v_\varphi\over r}\bigg|^2dx-{\nu\over2}\intop_{\Omega_\varepsilon}
{(\tilde v_\varphi^4)_{,r}\over r^2}drdz\cr
&={3\over4}\nu\intop_{\Omega_\varepsilon}
\bigg|{\nabla\tilde v_\varphi^2\over r}\bigg|^2dx-{\nu\over2}
\intop_{\Omega_\varepsilon}\bigg({\tilde v_\varphi^4\over r^2}\bigg)_{,r}drdz-
\nu\intop_{\Omega_\varepsilon}{\tilde v_\varphi^4\over r^4}dx\cr
&={3\over4}\nu\intop_{\Omega_\varepsilon}\bigg|\nabla
{\tilde v_\varphi^2\over r}+{\tilde v_\varphi^2\over r^2}\nabla r\bigg|^2dx+
{\nu\over2}\intop_{-a}^a{\tilde v_\varphi^4\over r^2}\bigg|_{r=\varepsilon}
dz-\nu\intop_{\Omega_\varepsilon}{\tilde v_\varphi^4\over r^4}dx\cr
&={3\over4}\nu\intop_{\Omega_\varepsilon}\bigg[\bigg|\nabla
{\tilde v_\varphi^2\over r}\bigg|^2+{\tilde v_\varphi^4\over r^4}+2\partial_r
{\tilde v_\varphi^2\over r}
{\tilde v_\varphi^2\over r^2}\bigg]dx-\nu\intop_{\Omega_\varepsilon}
{\tilde v_\varphi^4\over r^4}dx\cr
&\quad+{\nu\over2}\intop_{-a}^a
{\tilde v_\varphi^4\over r^2}\bigg|_{r=\varepsilon}dz
={3\over4}\nu\intop_{\Omega_\varepsilon}\bigg|\nabla
{\tilde v_\varphi^2\over r}\bigg|^2dx-{\nu\over4}\intop_{\Omega_\varepsilon}
{\tilde v_\varphi^4\over r^4}dx+{\nu\over2}\intop_{-a}^a
{\tilde v_\varphi^4\over r^2}\bigg|_{r=\varepsilon}dz\cr
&\quad+{3\over2}\nu
\intop_{\Omega_\varepsilon}\partial_r{\tilde v_\varphi^2\over r}
{\tilde v_\varphi^2\over r}drdz
={3\over4}\nu\intop_{\Omega_\varepsilon}\bigg|\nabla
{\tilde v_\varphi^2\over r}\bigg|^2dx-{\nu\over4}\intop_{\Omega_\varepsilon}
{\tilde v_\varphi^4\over r^4}dx\cr
&\quad-{\nu\over4}\intop_{-a}^a
{\tilde v_\varphi^4\over r^2}\bigg|_{r=\varepsilon}dz.\cr}
$$
Using the above considerations in (4.16), employing that
$\tilde v_\varphi|_{r=\varepsilon}=0$ we obtain after integration with respect
to time the equality
$$\eqal{
&{1\over4}\intop_{\Omega_\varepsilon}
{\tilde v_\varphi^4(t)\over r^2}dx+{3\over4}\nu\intop_{\Omega_\varepsilon^t}
\bigg|\nabla{\tilde v_\varphi^2\over r}\bigg|^2dxdt'+{3\over4}\nu
\intop_{\Omega_\varepsilon^t}{\tilde v_\varphi^4\over r^4}dxdt'\cr
&=-{3\over2}\intop_{\Omega_\varepsilon^t}{v_r\over r}
{\tilde v_\varphi^4\over r^2}dxdt'+\intop_{\Omega_\varepsilon^t}v\cdot\nabla
\zeta_1v_\varphi{\tilde v_\varphi|\tilde v_\varphi|^2\over r^2}dxdt'\cr
&\quad-\nu\intop_{\Omega_\varepsilon^t}[2\nabla v_\varphi\nabla\zeta_1+
v_\varphi\Delta\zeta_1]{\tilde v_\varphi|\tilde v_\varphi|^2\over r^2}dxdt'+
\intop_{\Omega_\varepsilon^t}\tilde f_\varphi\tilde v_\varphi
{\tilde v_\varphi^2\over r^2}dxdt'\cr
&\quad+{1\over4}\intop_{\Omega_\varepsilon}{\tilde v_\varphi^4(0)\over r^2}dx.
\cr}
\leqno(4.17)
$$
The second term on the r.h.s. of (4.17) is estimated by
$$
c(1/r_0)\|rv_\varphi\|_{L_\infty(\Omega_\varepsilon^t)}^3
\intop_{\Omega_\varepsilon^t}v^2dxdt'\le c(1/r_0)d_5^2d_6^3,
$$
where $d_5$ appears in (2.2) and $d_6$ in (3.1).
Similarly the third integral is bounded by
$$
c(1/r_0)\|rv_\varphi\|_{L_\infty(\Omega_\varepsilon^t)}^2
\intop_{\Omega_\varepsilon^t}(|\nabla v_\varphi|^2+|v_\varphi|^2)dxdt'\le c
(1/r_0)d_5^2d_6^2.
$$
Applying the H\"older and the Young inequalities to the term with
$\tilde f_\varphi$ we estimate it by
$$\eqal{
&\intop_{\Omega_\varepsilon^t}|\tilde f_\varphi|{|\tilde v_\varphi|^3\over r^2}
dxdt'=\intop_{\Omega_\varepsilon^t}{|\tilde f_\varphi|\over r^{1/2}}
{|\tilde v_\varphi|^3\over r^{3/2}}dxdt'\cr
&\le\bigg\|
{\tilde f_\varphi\over r^{1/2}}\bigg\|_{L_{20/11}(\Omega_\varepsilon^t)}
\bigg\|
{\tilde v_\varphi^2\over r}\bigg\|_{L_{10\over3}(\Omega_\varepsilon^t)}^{3/2}
\cr
&\le{\varepsilon_1^{4/3}\over4/3}
\bigg\|{\tilde v_\varphi^2\over r}\bigg\|_{L_{10/3}(\Omega_\varepsilon^t)}^2+
{1\over4\varepsilon_1^4}\bigg\|
{\tilde f_\varphi\over r^{1/2}}\bigg\|_{L_{20/11}(\Omega_\varepsilon^t)}^4.\cr}
$$
Employing the above estimates in (4.17) and choosing $\varepsilon_1$
sufficiently small we obtain (4.15). This concludes the proof.

\noindent
To examine problem (4.5) it is convenient to introduce new quantities $\eta$
and $\vartheta$ by the relations
$$
\psi=\eta r^2,\quad \chi=\vartheta r.
\leqno(4.18)
$$
Then problem (4.5) assumes the form
$$
\Delta\eta+{2\eta_{,r}\over r}=\vartheta\quad {\rm in}\ \ \Omega,\quad
\eta|_S=0,\ \ \vartheta=0\quad {\rm for}\ \ r\le\varepsilon.
\leqno(4.19)
$$
Introducing the new quantities
$$
\tilde\eta=\eta\zeta_1^2,\quad \tilde\vartheta=\vartheta\zeta_1^2
\leqno(4.20)
$$
we see that problem (4.19) takes the form
$$\eqal{
&\Delta\tilde\eta+{2\tilde\eta_{,r}\over r}=\tilde\vartheta+2\nabla\eta
\nabla\zeta_1^2+\eta\Delta\zeta_1^2+{2\over r}\zeta_{1,r}^2\eta\equiv
\tilde\vartheta+\vartheta_1\equiv\vartheta_2,\cr
&\tilde\eta|_{S_2}=0,\ \ \tilde\eta|_{r=2r_0}=0,\ \ \vartheta_2=0\quad
{\rm for}\ \ r\le\varepsilon.\cr}
\leqno(4.21)
$$

\proclaim Lemma 4.3.
Assume that $\vartheta_{,r}\in L_2(\Omega^T)$ and $v$ is a weak solution
satisfying (2.2). Then
$$
\intop_{\Omega^t}|\nabla\tilde\eta_{,rz}|^2dxdt'+6\intop_{\Omega^t}
{\tilde\eta_{,rz}^2\over r^2}dxdt'\le\intop_{\Omega^t}\vartheta_{2,r}^2dxdt'
+c(1/r_0)d_5^2,
\leqno(4.22)
$$
where $d_5$ is introduced in (2.2).

\Proof
Differentiating $(4.21)_1$ with respect to $r$ and $z$ yields
$$
\Delta\tilde\eta_{,rz}-{3\tilde\eta_{,rz}\over r}+2{\tilde\eta_{,rrz}\over r}=
\vartheta_{2,rz}.
\leqno(4.23)
$$
Multiplying (4.23) by $\tilde\eta_{,rz}$ and integrating over $\Omega$ implies
$$\eqal{
&\intop_\Omega\Delta\tilde\eta_{,rz}\tilde\eta_{,rz}dx-3\intop_\Omega
{\tilde\eta_{,rz}^2\over r^2}dx+2\intop_\Omega\tilde\eta_{,rrz}\tilde\eta_{,rz}
drdz\cr
&=\intop_\Omega\vartheta_{2,rz}\tilde\eta_{,rz}dx.\cr}
$$
Integrating by parts in the first term on the l.h.s. yields
$$\eqal{
&\intop_{S_2}\bar n\cdot\nabla\tilde\eta_{,rz}\tilde\eta_{,rz}dS_2-
\intop_\Omega|\nabla\tilde\eta_{,rz}|^2dx-3\intop_\Omega
{\tilde\eta_{,rz}^2\over r^2}dx\cr
&\quad+\intop_\Omega(\tilde\eta_{,rz}^2)_{,r}drdz=\intop_\Omega
(\vartheta_{2,r}\tilde\eta_{,rz})_{,z}dx-\intop_\Omega\vartheta_{2,r}
\tilde\eta_{,rzz}dx.\cr}
\leqno(4.24)
$$
The first integral on the l.h.s. of (4.24) equals
$$
I_1\equiv\intop_{S_2}\tilde\eta_{,rzz}\tilde\eta_{,rz}dS_2.
$$
Expressing (4.21) in the cylindrical coordinates yields
$$
\tilde\eta_{,rr}+\tilde\eta_{,zz}+3{\tilde\eta_{,r}\over r}=\tilde\vartheta+
\vartheta_1.
\leqno(4.25)
$$
In view of (4.21) we have that $\tilde\eta|_{S_2}=0$ also
$\tilde\eta_{,r}|_{S_2}=0$ and $\tilde\eta_{,rr}|_{S_2}=0$. Therefore
$\vartheta_1|_{S_2}=0$.

\noindent
Projecting (4.25) on $S_2$ and using that $\tilde\vartheta|_{S_2}=0$ and
$\tilde\vartheta_1|_{S_2}=0$ we have that $\tilde\eta_{,zz}|_{S_2}=0$ so also
$\tilde\eta_{,zzr}|_{S_2}=0$. Therefore, $I_1=0$.

\noindent
Since $\vartheta_2|_{S_2}=0$ so also $\vartheta_{2,r}|_{S_2}=0$. Then the
first term on the r.h.s. of (4.24) vanishes. Hence (4.24) takes the form
$$\eqal{
&\intop_\Omega|\nabla\tilde\eta_{,rz}|^2dx+3\intop_\Omega
{\tilde\eta_{,rz}^2\over r^2}dx+\intop_{-a}^a\tilde\eta_{,rz}^2|_{r=0}dz\cr
&=\intop_\Omega\vartheta_{2,r}\tilde\eta_{,rzz}dx.\cr}
\leqno(4.26)
$$
Applying the Cauchy inequality to the r.h.s. of (4.26) we estimate it by
$$
{1\over2}\intop_\Omega\tilde\eta_{,rzz}^2dx+{1\over2}\intop_\Omega
\vartheta_{2,r}^2dx.
$$
Using this in (4.26) implies
$$
\intop_\Omega|\nabla\tilde\eta_{,rz}|^2dx+6\intop_\Omega
{\tilde\eta_{,rz}^2\over r^2}dx\le\intop_\Omega\vartheta_{2,r}^2dx.
\leqno(4.27)
$$
To estimate the r.h.s. of (4.27) we examine
$$\eqal{
&\intop_\Omega\vartheta_{1,r}^2dx\le c(1/r_0)
\intop_{\Omega_{,\zeta_{1,r}}}(\eta_{,rr}^2+\eta_{,r}^2+\eta^2)dx\cr
&\le c(1/r_0)\intop_{\Omega_{,\zeta_{1,r}}}(\psi_{,rr}^2+\psi_{,r}^2+
\psi^2)dx\le c(1/r_0)\intop_\Omega(v_{z,r}^2+v_z^2)dx,\cr}
$$
where we used the equality $\psi(r,z)=\intop_R^r\psi_{,r'}(r',z)dr'$.
Therefore, (4.27) takes the form
$$\eqal{
&\intop_\Omega|\nabla\tilde\eta_{,rz}|^2dx+6\intop_\Omega
{\tilde\eta_{,rz}^2\over r^2}dx\le\intop_\Omega\vartheta_{2,r}^2dx\cr
&\quad+c(1/r_0)\intop_\Omega(v_{z,r}^2+v_z^2)dx.\cr}
\leqno(4.28)
$$
Integrating (4.28) with respect to time and using (2.2) yield (4.22). This
concludes the proof.

\proclaim Corollary 4.4.
Since $\tilde\eta_{,z}={\tilde\psi_{,z}\over r^2}={\tilde v_r\over r}$,
$\tilde\eta_{,zr}=\bigg({\tilde v_r\over r}\bigg)_{,r}$ and
$\tilde\vartheta={\tilde\chi\over r}$ then (4.22) takes the form
$$\eqal{
&\intop_{\Omega^t}\bigg|\nabla\bigg({\tilde v_r\over r}\bigg)_{,r}\bigg|^2
dxdt'+6\intop_{\Omega^t}{1\over r^2}\bigg({\tilde v_r\over r}\bigg)_{,r}^2dxdt'\cr
&\le\intop_{\Omega_\varepsilon^t}\bigg({\tilde\chi\over r}\bigg)_{,r}^2dxdt'
+c(1/r_0)d_5^2.\cr}
\leqno(4.29)
$$

\proclaim Lemma 4.5. \hskip-1pt Assume \hskip-1pt that \hskip-1pt
${\bar\chi(0)\over r}\in L_2(\Omega_\varepsilon)$, \hskip-1pt
${\tilde F\over r}\in L_2(\Omega_\varepsilon^T)$,  \hskip-1pt
${\tilde f_\varphi\over\sqrt{r}}\in L_{20/11}(\Omega_\varepsilon^T)$,\break
${\tilde v_\varphi(0)\over\sqrt{r}}\in L_4(\Omega_\varepsilon)$. Assume that
(2.2) and (3.1) hold, where $d_5$ and $d_6$ are introduced. Assume that
$$
\|u\|_{L_\infty(\Omega_\zeta^T)}\le\root{4}\of{5\over8}\nu,
\leqno(4.30)
$$
which can be satisfied because $u$ vanishes on the axis of symmetry, it is
the H\"older continuous and $\supp\zeta$ is sufficiently small.
Then the following a priori inequality holds
$$
\bigg\|{\tilde\chi\over r}\bigg\|_{V_2^0(\Omega_\varepsilon^t)}^2\le c(1/r_0)
d_5\|\chi\|_{L_{20/7}(\Omega_{\varepsilon,\zeta_{1,r}}^t)}^2+cA^2,\quad t\le T,
\leqno(4.31)
$$
where
$$\eqal{
A^2&=c(1/r_0)d_5^2[1+(1+d_6)d_6^2]+
\bigg\|{\tilde\chi(0)\over r}\bigg\|_{L_2(\Omega_\varepsilon)}^2+
\bigg\|{\tilde F\over r}\bigg\|_{L_2(\Omega_\varepsilon^t)}^2\cr
&\quad+\bigg\|
{\tilde f_\varphi\over\sqrt{r}}\bigg\|_{L_{20/11}(\Omega_\varepsilon^t)}^4+
\bigg\|{\tilde v_\varphi(0)\over\sqrt{r}}\bigg\|_{L_4(\Omega_\varepsilon)}^4.
\cr}
$$

\Proof
From (4.10) and (4.29) we have
$$\eqal{
&\intop_{\Omega_\varepsilon^t}\bigg|\nabla\bigg(
{\tilde v_r\over r}\bigg)_{,r}\bigg|^2dxdt'+6\intop_{\Omega_\varepsilon^t}
{1\over r^2}\bigg({\tilde v_r\over r}\bigg)_{,r}^2dxdt'\cr
&\le{2\over\nu^2}\intop_{\Omega_\varepsilon^t}{\tilde v_\varphi^4\over r^4}
dxdt'+c(1/r_0)d_5\|\chi\|_{L_{20/7}(\Omega_{\varepsilon,\zeta_{1,r}}^t)}+
cA_1^2,\cr}
\leqno(4.32)
$$
where
$$
A_1^2=c(1/r_0)d_5^2+
\bigg\|{\tilde\chi(0)\over r}\bigg\|_{L_2(\Omega_\varepsilon)}^2+
\bigg\|{\tilde F\over r}\bigg\|_{L_2(\Omega_\varepsilon^t)}^2.
$$
To estimate the first term on the r.h.s. of (4.32) we use (4.15) in the form
$$
\intop_{\Omega_\varepsilon^t}{\tilde v_\varphi^4\over r^4}dxdt'\le{2\over\nu}
\intop_{\Omega_\varepsilon^t}{|v_r|\over r}{\tilde v_\varphi^4\over r^2}
dxdt'+cA_2^2,
\leqno(4.34)
$$
where
$$
A_2^2=c(1/r_0)d_5^2(1+d_6)d_6^2+\bigg\|
{\tilde f_\varphi\over\sqrt{r}}\bigg\|_{L_{20/11}(\Omega_\varepsilon^t)}^4+
\bigg\|{\tilde v_\varphi(0)\over\sqrt{r}}\bigg\|_{L_4(\Omega_\varepsilon)}^4.
\leqno(4.35)
$$
Now we examine the first term on the r.h.s. of (4.34) in the following way
$$\eqal{
&{2\over\nu}\intop_{\Omega_\varepsilon^t}{|\tilde v_r|\over r^3}r^2v_\varphi^2
{\tilde v_\varphi^2\over r^2}dxdt'={2\over\nu}\intop_{\Omega_\varepsilon^t}
{|\tilde v_r|^2\over r^3}u^2{\tilde v_\varphi^2\over r^2}dxdt'\cr
&\le{2\over\nu}\bigg[{\varepsilon\over2}\intop_{\Omega_\varepsilon^t}
{\tilde v_\varphi^4\over r^4}dxdt'+{1\over2\varepsilon}
\|u\|_{L_\infty(\Omega_\zeta^t)}^4\intop_{\Omega_\varepsilon^t}
{\tilde v_r^2\over r^6}dxdt'\bigg].\cr}
$$
Setting $\varepsilon={\nu\over2}$ in the above inequality we obtain from
(4.34) the inequality
$$
\intop_{\Omega_\varepsilon^t}{\tilde v_\varphi^4\over r^4}dxdt'\le{4\over\nu^2}
\|u\|_{L_\infty(\Omega_\zeta^t)}^4\intop_{\Omega_\varepsilon^t}
{\tilde v_r^2\over r^6}dxdt'+cA_2^2.
\leqno(4.36)
$$
Using (4.36) in (4.32) and applying the Hardy inequality
$$
\intop_{\Omega_\varepsilon}{1\over r^4}\bigg|{\tilde v_r\over r}\bigg|^2dx\le
\intop_{\Omega_\varepsilon}{1\over r^2}\bigg({\tilde v_r\over r}\bigg)_{,r}^2dx
\leqno(4.37)
$$
we obtain
$$\eqal{
&6\intop_{\Omega_\varepsilon^t}{1\over r^2}
\bigg({\tilde v_r\over r}\bigg)_{,r}^2dxdt'\le{8\over\nu^4}
\|u\|_{L_\infty(\Omega_\zeta^t)}^4\intop_{\Omega_\varepsilon^t}{1\over r^2}
\bigg({\tilde v_r\over r}\bigg)_{,r}^2dx\cr
&\quad+c(1/r_0)d_5\|\chi\|_{L_{20/7}(\Omega_{\varepsilon,\zeta_{1,r})}^t}^2+
c(A_1^2+A_2^2).\cr}
\leqno(4.38)
$$
Assuming that
$$
6-{8\over\nu^4}\|u\|_{L_\infty(\Omega_\zeta^t)}^4\ge1,
$$
which can be expressed in the form
$$
\|u\|_{L_\infty(\Omega_\zeta^t)}\le\root{4}\of{5\over8}\nu,
\leqno(4.39)
$$
we obtain from (4.38) the inequality
$$
\intop_{\Omega_\varepsilon^t}{1\over r^2}
\bigg({\tilde v_r\over r}\bigg)_{,r}^2dxdt'\le c(1/r_0)d_5
\|\chi\|_{L_{20/7}(\Omega_{\varepsilon,\zeta_{1,r}}^t)}^2+c(A_1^2+A_2^2).
\leqno(4.40)
$$
Employing (4.39) and (4.40) in (4.36) implies
$$\eqal{
&\intop_{\Omega_\varepsilon^t}{\tilde v_\varphi^4\over r^4}dxdt'\le
{4\over\nu^2}{5\over8}\nu^4[c(1/r_0)d_5
\|\chi\|_{L_{20/7}(\Omega_{\varepsilon,\zeta_{1,r}}^t)}^2\cr
&\quad+c(A_1^2+A_2^2)]+cA_2^2\cr
&\le c(1/r_0)d_5\|\chi\|_{L_{20/7}(\Omega_{\varepsilon,\zeta_{1,r}}^t)}^2+
c(A_1^2+A_2^2).\cr}
\leqno(4.41)
$$
In view of (4.41) inequality (4.10) takes the form
$$\eqal{
&\bigg\|{\tilde\chi\over r}\bigg\|_{L_\infty(0,T;L_2(\Omega_\varepsilon))}^2+
\bigg\|\nabla{\tilde\chi\over r}\bigg\|_{L_2(\Omega_\varepsilon^t)}^2\cr
&\le c(1/r_0)d_5\|\chi\|_{L_{20/7}(\Omega_{\varepsilon,\zeta_{1,r}}^t)}^2
+c(A_1^2+A_2^2).\cr}
\leqno(4.42)
$$
This concludes the proof.

\proclaim Lemma 4.6.
Let the assumptions of Lemma 4.5 hold. Then
$$
\bigg\|{\chi\over r}\bigg\|_{L_{20/7}(\Omega_\varepsilon^t)}\le cA,
\leqno(4.43)
$$
where $A$ is introduced in Lemma 4.5.

\Proof
From (4.31) we have
$$
\bigg\|{\tilde\chi\over r}\bigg\|_{L_{10/3}(\Omega_\varepsilon^t)}\le c_1
\bigg\|{\chi\over r}\bigg\|_{L_{10/3}(\Omega_{\varepsilon,\zeta_{1,r}}^t)}+cA.
\leqno(4.44)
$$
Let us introduce the sets
$$
\Omega_\varepsilon^{(\lambda)}=\{(r,z)\in\Omega_\varepsilon:\
0<\varepsilon<r\le r_0-\lambda,\ |z|<a\}
$$
and connect with them a set of cut-off functions such that
$$
\zeta^{(\lambda)}=\left\{\eqal{
&1\quad &{\rm for}\ \ (r,z)\in\Omega_\varepsilon^{(\lambda)}\cr
&0\quad &{\rm for}\ \ (r,z)\in\Omega_\varepsilon\setminus
\Omega_\varepsilon^{(\lambda/2)}\cr}\right.
$$
Let $\vartheta={\chi\over r}$. Then (4.44) can be expressed in the form
$$
\|\vartheta\|_{L_{10/3}(\Omega_\varepsilon^{(\lambda)}\times(0,t))}\le c_1
\|\vartheta\|_{L_{10/3}(\Omega_\varepsilon^{(\lambda/2)}\setminus
\Omega_\varepsilon^{(\lambda)}\times(0,t))}+cA.
\leqno(4.45)
$$
From (4.45) we have
$$\eqal{
&\intop_{\Omega_\varepsilon^{(\lambda)}\times(0,t)}|\vartheta|^{10/3}dxdt'\cr
&\le c_1c_1^{10/3}\intop_{\Omega_\varepsilon^{(\lambda/2)}\setminus
\Omega_\varepsilon^{(\lambda)}\times(0,t)}
|\vartheta|^{10/3}dxdt'+c_2c^{10/3}A^{10/3},\cr}
\leqno(4.46)
$$
where $c_2=2^{10/3-1}$. Adding
$$
c_2c_1^{10/3}\intop_{\Omega_\varepsilon^{(\lambda)}\times(0,t)}
|\vartheta|^{10/3}dxdt'
$$
to both sides of (4.46) we obtain
$$\eqal{
&\intop_{\Omega_\varepsilon^{(\lambda)}\times(0,t)}|\vartheta|^{10/3}dxdt'\cr
&\le{c_2c_1^{10/3}\over1+c_2c_1^{10/3}}
\intop_{\Omega_\varepsilon^{(\lambda/2)}\times(0,t)}
|\vartheta|^{10/3}dxdt'+{c_2c^{10/3}\over1+c_2c_1^{10/3}}A^{10/3}.\cr}
\leqno(4.47)
$$
Introducing the notation
$$\eqal{
&f(\lambda)=\intop_{\Omega_\varepsilon^{(\lambda)}\times(0,t)}
|\vartheta|^{10/3}dxdt',\quad
\mu={c_2c_1^{10/3}\over1+c_2c_1^{10/3}}<1\cr
&K={c_2c^{10/3}\over1+c_2c_1^{10/3}}A,\cr}
$$
we obtain from (4.47) the inequality
$$
f(\lambda)\le\mu f(\lambda/2)+K,
$$
which implies the estimate
$$
f(\lambda)\le\sum_{j=0}^\infty\mu^jK={1\over1-\mu}K.
$$
Therefore, Lemma 4.6 is proved.

\section{5. Estimate for $\chi$ in a neighborhood located in a positive
distance from the axis of symmetry}

Let $\zeta_2=\zeta_2(r)$ be a smooth cut-off function such that
$\zeta_2(r)=0$ for $r\le r_0$ and $\zeta_2(r)=1$ for $r\ge2r_0$. Let
$\{\zeta_1(r),\zeta_2(r)\}$ compose a partition of unity in the radial
direction. Let us introduce the notation
$$\eqal{
&\bar\chi=\chi\zeta_2^2,\quad \bar v_\varphi=v_\varphi\zeta_2,\quad
\bar v'=v'\zeta_2^2,\quad \bar F=F\zeta_2^2,\quad \bar f'=f'\zeta_2^2,\cr
&\bar f_\varphi=f_\varphi\zeta_2.\cr}
\leqno(5.1)
$$
In view of (4.6) function $\bar\chi$ is a solution to the problem
$$\eqal{
&\bar\chi_{,t}+v\cdot\nabla\bar\chi-{v_r\over r}\bar\chi-\nu\bigg[\bigg(
r\bigg({\bar\chi\over r}\bigg)_{,r}\bigg)_{,r}+\bar\chi_{,zz}+2
\bigg({\bar\chi\over r}\bigg)_{,r}\bigg]\quad\cr
&=v\cdot\nabla\zeta_2^2\chi-\nu\bigg[(\chi\zeta_{2,r}^2)_{,r}+r
\bigg({\chi\over r}\bigg)_{,r}\zeta_{2,r}^2+2\bigg({\chi\over r}
\zeta_{2,r}^2\bigg)\bigg]\cr
&\quad+{2\bar v_\varphi\bar v_{\varphi,z}\over r}+\bar F\quad &{\rm in}\ \
\Omega^T,\cr
&\bar\chi|_{r=r_0}=0,\ \ \bar\chi|_{S_1\cup S_2}=0,\cr
&\bar\chi|_{t=0}=\bar\chi_0\quad &{\rm in}\ \ \Omega.\cr}
\leqno(5.2)
$$

\proclaim Lemma 5.1.
Assume that $v$ is a weak solution to problem (1.1) satisfying assumptions
of Lemma 2.1. Let the assumptions of Lemma 3.1 hold. Let $d_1$, $d_5$, $d_7$,
$d_6$ be constants introduced by (2.1), (2.2), (2.15), (3.1), respectively.
Let $\bar F\in L_2(0,T;L_{6/5}(\Omega))$, $\bar\chi(0)\in L_2(\Omega)$.
Then solutions to (5.2) satisfy the estimate
$$
\bigg\|{\bar\chi\over r}\bigg\|_{V_2^0(\Omega^t)}\le c
[c(1/r_0)d_1^2+1]A_0,\quad t\le T,
\leqno(5.3)
$$
where
$$
A_0^2=c(1/r_0)(d_5^2+d_6^2d_7^2)+\|\bar F\|_{L_2(0,T;L_{6/5}(\Omega))}^2+
{1\over r_0^2}\|\bar\chi(0)\|_{L_2(\Omega)}^2
\leqno(5.4)
$$
and $r_0$ is introduced by the definition of the cut-off function $\zeta_2(r)$.

\Proof
Multiplying $(5.2)_1$ by ${\bar\chi\over r^2}$, integrating over $\Omega$
and using the boundary conditions yields
$$\eqal{
&{1\over2}{d\over dt}\intop_\Omega{\bar\chi^2\over r^2}dx+\nu\intop_\Omega
\bigg|\nabla{\bar\chi\over r}\bigg|^2dx=\intop_\Omega v\cdot\nabla\zeta_2^2
\chi{\bar\chi\over r^2}dx\cr
&\quad-\nu\intop_\Omega\bigg[(\chi\zeta_{2,r}^2)_{,r}+r
\bigg({\chi\over r}\bigg)_{,r}\zeta_{2,r}^2+2{\chi\over r}\zeta_{2,r}^2\bigg]
{\bar\chi\over r^2}dx\cr
&\quad+2\intop_\Omega{\bar v_\varphi\bar v_{\varphi,z}\over r^2}
{\bar\chi\over r}dx+\intop_\Omega{\bar F\over r}{\bar\chi\over r}dx.\cr}
\leqno(5.5)
$$
Now we estimate the terms from the r.h.s. of (5.5). We estimate the first
term by
$$
\varepsilon_1\bigg\|{\bar\chi\over r}\bigg\|_{L_6(\Omega)}^2+
c(1/r_0,1/\varepsilon_1)\|v\chi\|_{L_{6/5}(\Omega_{\zeta_{2,r}})}^2.
$$
The first term under the square bracket of the second term yields
$$
\intop_\Omega\chi_{,r}\zeta_{2,r}^2\zeta_2^2{\chi\over r}drdz+
\intop_\Omega\chi^2{1\over r^2}(\zeta_{2,r}^2)_{,r}\zeta_2^2dx\equiv I_1,
$$
where the first integral in $I_1$ equals
$$
-{1\over2}\intop_\Omega\chi^2\bigg({1\over r}\zeta_{2,r}^2\zeta_2^2\bigg)drdz.
$$
Hence
$$
|I_1|\le c(1/r_0)\intop_{\Omega_{\zeta_{2,r}}}\chi^2dx.
$$
Similarly, the integral with the second term under the square bracket implies
$$
\intop_\Omega\bigg({\chi\over r}\bigg)_{,r}{\chi\over r}\zeta_{2,r}^2
\zeta_2^2rdrdz=
-{1\over2}\intop_\Omega{\chi^2\over r^2}(\zeta_{2,r}^2\zeta_2^2r)_{,r}drdz
\equiv I_2.
$$
Then
$$
|I_2|\le c(1/r_0)\intop_{\Omega_{\zeta_{2,r}}}\chi^2dx.
$$
Finally, the integral with the last term under the square bracket is bounded by
$$
c(1/r_0)\intop_{\Omega_{\zeta_{2,r}}}\chi^2dx.
$$
Summarizing, the second term on the r.h.s. of (5.5) is bounded by
$$
c(1/r_0)\intop_{\Omega_{\zeta_{2,r}}}\chi^2dx.
$$
The third term on the r.h.s. of (5.5) yields
$$\eqal{
&\bigg|\intop_\Omega{\bar v_\varphi^2\over r^2}
\bigg({\bar\chi\over r}\bigg)_{,z}dx\bigg|\le\varepsilon_2\intop_\Omega
\bigg({\bar\chi\over r}\bigg)_{,z}^2dx+c(1/\varepsilon_2)
\intop_{\Omega_{\zeta_2}}{\bar v_\varphi^4\over r^4}dx\cr
&\le\varepsilon_2\intop_\Omega\bigg({\bar\chi\over r}\bigg)_{,z}^2dx+c
(1/\varepsilon_2,1/r_0)\intop_\Omega v_\varphi^4dx.\cr}
$$
Finally, the last term on the r.h.s. of (5.5) is bounded by
$$
\varepsilon_3\bigg\|{\bar\chi\over r}\bigg\|_{L_6(\Omega)}^2+c(1/\varepsilon_3)
\bigg\|{\bar F\over r}\bigg\|_{L_{6/5}(\Omega)}^2.
$$
Employing the above estimates in (5.5) and choosing
$\varepsilon_1-\varepsilon_3$ sufficiently small yield
$$\eqal{
&{d\over dt}\intop_\Omega\bigg|{\bar\chi\over r}\bigg|^2dx+\nu\intop_\Omega
\bigg|\nabla{\bar\chi\over r}\bigg|^2dx\le c(1/r_0)\|v\|_{L_2(\Omega)}^2
\|\chi\|_{L_3(\Omega_{\zeta_{2,r}})}^2\cr
&\quad+c(1/r_0)\|\chi\|_{L_2(\Omega_{\zeta_{2,r}})}^2+c(1/r_0)\intop_\Omega
v_\varphi^4dx+c\bigg\|{\bar F\over r}\bigg\|_{L_{6/5}(\Omega)}^2.\cr}
\leqno(5.6)
$$
Integrating (5.6) with respect to time and exploiting estimates (2.1) and
(3.6) we obtain
$$\eqal{
&\bigg\|{\bar\chi\over r}\bigg\|_{L_\infty(0,t;L_2(\Omega))}^2+\nu
\bigg\|\nabla{\bar\chi\over r}\bigg\|_{L_2(\Omega^t)}^2\cr
&\le c(1/r_0)d_1^2\|\chi\|_{L_2(0,t;L_3(\Omega_{\zeta_{2,r}}))}^2+
c(1/r_0)(d_5^2+d_6^2d_7^2)\cr
&\quad+c\bigg\|{\bar F\over r}\bigg\|_{L_2(0,t;L_{6/5}(\Omega))}^2+
\bigg\|{\bar\chi(0)\over r}\bigg\|_{L_2(\Omega)}^2,\cr}
\leqno(5.7)
$$
where $t\le T$, $d_1$ is introduced in (2.1), $d_5$ in (2.2), $d_6$ in (3.1)
and $d_7$ in (2.15).
To estimate the first term on the r.h.s. of (5.7) by data we express (5.7)
in the form
$$
\|\bar\chi\|_{L_3(\Omega^t)}\le c_1\|\chi\|_{L_3(\Omega_{\zeta_{2,r}}^t)}+A_0,
\leqno(5.8)
$$
where $c_1=c(1/r_0)d_1$ and $A_0$ is defined by (5.4).

\noindent
To apply the local considerations (see \cite[Ch. 4, Sect. 10]{LSU}) we
introduce the sets
$\Omega^{(\lambda)}=\{(r,z)\in\Omega:\ r\ge r'_0+\lambda\}$ and corresponding
cut-off functions $\zeta^{(\lambda)}(x)$ such that $\zeta^{(\lambda)}(x)=1$
for $x\in\Omega^{(\lambda)}$ and $\zeta^{(\lambda)}(x)=0$ for
$x\in\Omega\setminus\Omega^{(\lambda/2)}$, so
$|\nabla\zeta^{(\lambda)}|\le{c\over\lambda}$. Moreover, we assume that
$r'_0+\lambda=2r_0$.

Then (5.8) can be expressed in the form
$$
\|\chi\|_{L_3(\Omega^{(\lambda)}\times(0,t))}\le c_1
\|\chi\|_{L_3(\Omega^{(\lambda/2)}\setminus\Omega^{(\lambda)}\times(0,t))}+A_0.
\leqno(5.9)
$$
Hence
$$
\|\chi\|_{L_3(\Omega^{(\lambda)}\times(0,t))}^3\le4c_1^3
\|\chi\|_{L_3(\Omega^{(\lambda/2)}\setminus\Omega^{(\lambda)}\times(0,t))}^3
+4A_0^3.
\leqno(5.10)
$$
By the filling-hole argument we have
$$
\intop_{\Omega^{(\lambda)}\times(0,t)}|\chi|^3dxdt'\le{4c_1^3\over4c_1^3+1}
\intop_{\Omega^{(\lambda/2)}\times(0,t)}|\chi|^3dxdt'
+{4\over4c_1^3+1}A_0^3.
\leqno(5.11)
$$
Introducing the notation
$$
f(\lambda)=\intop_{\Omega^{(\lambda)}\times(0,t)}|\chi|^3dxdt',\quad
\mu={4c_1^3\over4c_1^3+1}<1,\quad K={4\over4c_1^3+1}A_0^3
$$
we obtain from (5.11) the inequality
$$
f(\lambda)\le\mu f(\lambda/2)+K
$$
which implies the estimate
$$
f(\lambda)\le\sum_{j=0}^\infty\mu^jK={1\over1-\mu}K.
$$
Employing the estimate in (5.7) we get (5.3). This concludes the proof.

\section{6. Estimate for $v'$}

From Lemmas 4.5 and 4.6 we have
$$
\bigg\|{\tilde\chi\over r}\bigg\|_{V_2^0(\Omega_\varepsilon^t)}\le cA,\quad
t\le T,
\leqno(6.1)
$$
where
$$\eqal{
A^2&=c(1/r_0)d_5^2[1+(1+d_6)d_6^2]+
\bigg\|{\tilde\chi(0)\over r}\bigg\|_{L_2(\Omega_\varepsilon)}^2\cr
&\quad+\bigg\|{\tilde F\over r}\bigg\|_{L_2(\Omega_\varepsilon^t)}^2+
\bigg\|{\tilde f_\varphi\over\sqrt{r}}\bigg\|_{L_{20/11}(\Omega_\varepsilon^t)}^4
+\bigg\|{\tilde v_\varphi(0)\over\sqrt{r}}\bigg\|_{L_4(\Omega_\varphi)}^4,
\quad t\le T.\cr}
\leqno(6.2)
$$
Next, Lemma 5.1 implies
$$
\bigg\|{\bar\chi\over r}\bigg\|_{V_2^0(\Omega^t)}\le c[c(1/r_0)d_1^2+1]A_0,\quad
t\le T,
\leqno(6.3)
$$
where
$$\eqal{
A_0^2&=c(1/r_0)(d_5^2+d_6^2d_7^2)+\|\bar F\|_{L_2(0,t;L_{6/5}(\Omega))}^2\cr
&\quad+{1\over r_0^2}\|\bar\chi(0)\|_{L_2(\Omega)}^2,\quad t\le T.\cr}
\leqno(6.4)
$$
Inequalities (6.1) and (6.3) imply
$$
\bigg\|{\chi\over r}\bigg\|_{V_2^0(\Omega_\varepsilon^t)}\le cA+
c[c(1/r_0)d_1^2+1]A_0\equiv cA_*,\quad t\le T.
\leqno(6.5)
$$
Let us consider problem (4.5).
It is convenient to introduce new quantities $\eta$ and $\vartheta$ by the
relations
$$
\psi=\eta r^2,\quad \chi=\vartheta r.
\leqno(6.6)
$$
Then problem (4.5) assumes the form
$$
\Delta\eta+{2\eta_{,r}\over r}=\vartheta\quad {\rm in}\ \ \Omega,\quad
\eta|_S=0,\ \ \vartheta=0\quad {\rm for}\ \ r\le\varepsilon.
\leqno(6.7)
$$
Since (6.5) holds we have

\proclaim Lemma 6.1.
Assume that $\vartheta\in H^1(\Omega)$.
Then for a sufficiently smooth solution of (6.7) we have
$$\eqal{
&\intop_\Omega(\eta^2+|\nabla\eta|^2+|\nabla\eta_{,r}|^2+|\nabla\eta_{,z}|^2+
|\nabla\eta_{,zr}|^2+|\nabla\eta_{,zz}|^2)dx\cr
&\quad+\intop_\Omega\bigg({\eta_{,r}^2\over r^2}+{\eta_{,zr}^2\over r^2}\bigg)
dx+\intop_{-a}^a(\eta_{,r}^2|_{r=R}+\eta_{,rz}^2|_{r=R})dz\cr
&\quad+\intop_{-a}^a(\eta^2|_{r=0}+\eta_{,r}^2|_{r=0}+\eta_{,z}^2|_{r=0}+
\eta_{,rz}^2|_{r=0}+\eta_{,zz}^2|_{r=0})dz\cr
&\le c\intop_\Omega(\vartheta^2+\vartheta_{,r}^2+\vartheta_{,z}^2)dx.\cr}
\leqno(6.8)
$$

\Proof
Multiplying $(6.7)_1$ by $\eta$, integrating over $\Omega$ and using the
boundary conditions yields
$$
-\intop_\Omega|\nabla\eta|^2dx+2\intop_\Omega{\eta_{,r}\eta\over r}rdrdz=
\intop_\Omega\eta\vartheta dx.
$$
Applying the Cauchy and the Young inequalities to the r.h.s. of the above
equality gives
$$
\intop_\Omega|\nabla\eta|^2dx+\intop_{-a}^a\eta^2|_{r=0}dz\le\varepsilon
\intop_\Omega\eta^2dx+c(1/\varepsilon)\intop_\Omega\vartheta^2dx.
$$
In view of sufficiently small $\varepsilon$ and the Poincare inequality
we get the estimate
$$
\intop_\Omega(\eta^2+|\nabla\eta|^2)dx+\intop_{-a}^a\eta^2|_{r=0}dz\le c
\intop_\Omega\vartheta^2dx.
\leqno(6.9)
$$
Differentiating $(6.7)_1$ with respect to $r$, multiplying the result by
$\eta_{,r}$ and integrating over $\Omega$ yields
$$
\intop_\Omega\Delta\eta_{,r}\eta_{,r}dx-3\intop_\Omega{\eta_{,r}^2\over r^2}dx
+2\intop_\Omega{\eta_{,rr}\eta_{,r}\over r}dx=\intop_\Omega\vartheta_{,r}
\eta_{,r}dx,
\leqno(6.10)
$$
where we used that $(6.7)_1$ takes the form
$$\eqal{
&\eta_{,rr}+\eta_{,zz}+3{\eta_{,r}\over r}=\vartheta,\quad {\rm so}\quad
\eta_{,rrr}+\eta_{,rzz}+3{\eta_{,rr}\over r}-3{\eta_{,r}\over r^2}=
\vartheta_{,r}\cr
&{\rm and}\quad \Delta\eta_{,r}+2{\eta_{,rr}\over r}-3{\eta_{,r}\over r^2}=
\vartheta_{,r}.\cr}
\leqno(6.11)
$$
Integrating by parts in (6.10) and using that $\eta_{,r}|_{S_2}=0$ we obtain
$$\eqal{
&\intop_{S_1}\bar n\cdot\nabla\eta_{,r}\eta_{,r}dS_1-\intop_\Omega
|\nabla\eta_{,r}|^2dx-3\intop_\Omega{\eta_{,r}^2\over r^2}dx+\intop_\Omega
(\eta_{,r}^2)_{,r}drdz\cr
&=\intop_\Omega\vartheta_{,r}\eta_{,r}dx.\cr}
\leqno(6.12)
$$
In view of (6.11) we have
$$
\bar n\cdot\nabla\eta_{,r}|_{S_1}=\eta_{,rr}|_{S_1}=-3{\eta_{,r}\over r}|_{S_1}.
$$
Moreover, integrating by parts in the r.h.s. of (6.12) gives
$$
\intop_\Omega(\vartheta\eta_{,r}r)_{,r}drdz-\intop_\Omega\vartheta
(\eta_{,rr}r+\eta_{,r})drdz,
$$
where the first integral vanishes because $\vartheta|_{r=R}=0$ and
$\vartheta|_{r=0}=0$.

\noindent
Therefore, (6.12) takes the form
$$\eqal{
&\intop_\Omega|\nabla\eta_{,r}|^2dx+3\intop_\Omega{\eta_{,r}^2\over r^2}dx+
3\intop_{-a}^a\eta_{,r}^2|_{r=R}dz-\intop_{-a}^a\eta_{,r}^2|_{r=R}dz\cr
&\quad+\intop_{-a}^a\eta_{,r}^2|_{r=0}dz\le\varepsilon\bigg(\intop_\Omega
\eta_{,rr}^2dx+\intop_\Omega{\eta_{,r}^2\over r^2}dx\bigg)+c(1/\varepsilon)
\intop_\Omega\vartheta^2dx\cr}
\leqno(6.13)
$$
For sufficiently small $\varepsilon$, (6.13) implies the estimate
$$\eqal{
&{1\over2}\intop_\Omega|\nabla\eta_{,r}|^2dx+{5\over2}\intop_\Omega
{\eta_{,r}^2\over r^2}dx+2\intop_{-a}^a\eta_{,r}^2|_{r=R}dz+\intop_{-a}^a
\eta_{,r}^2|_{r=0}dz\cr
&\le c\intop_\Omega\vartheta^2dx.\cr}
\leqno(6.14)
$$
Differentiating $(6.7)_1$ with respect to $z$, multiplying the result by
$\eta_{,z}$ and integrating over $\Omega$ gives
$$
\intop_\Omega\Delta\eta_{,z}\eta_{,z}dx+2\intop_\Omega\eta_{,rz}\eta_{,z}drdz=
\intop_\Omega\vartheta_{,z}\eta_{,z}dx.
$$
Integrating by parts and using that $\eta_{,zz}|_{S_2}=0$,
$\eta_{,z}|_{S_1}=0$, $\vartheta|_{S_2}=0$ we obtain
$$
\intop_\Omega|\nabla\eta_{,z}|^2dx-\intop_\Omega(\eta_{,z}^2)_{,r}drdz=
\intop_\Omega\vartheta\eta_{,zz}dx.
$$
Continuing, we have
$$
\intop_\Omega|\nabla\eta_{,z}|^2dx+\intop_{-a}^a\eta_{,z}^2|_{r=0}dz\le
\varepsilon\intop_\Omega\eta_{,zz}^2dx+c(1/\varepsilon)\intop_\Omega
\vartheta^2dx.
$$
Hence, for sufficiently small $\varepsilon$, we obtain
$$
\intop_\Omega|\nabla\eta_{,z}|^2dx+\intop_{-a}^a\eta_{,z}^2|_{r=0}dz\le
c\intop_\Omega\vartheta^2dx.
\leqno(6.15)
$$
Differentiating $(6.7)_1$ with respect to $r$ and $z$ yields
$$
\Delta\eta_{,rz}-{3\eta_{,zr}\over r}+2{\eta_{,zrr}\over r}=\vartheta_{,zr}.
\leqno(6.16)
$$
Multiplying (6.16) by $\eta_{,rz}$ and integrating over $\Omega$ implies
$$\eqal{
&\intop_\Omega\Delta\eta_{,rz}\eta_{,rz}dx-3\intop_\Omega
{\eta_{,rz}^2\over r^2}dx+2\intop_\Omega\eta_{,zrr}\eta_{,zr}drdz\cr
&=\intop_\Omega\vartheta_{,zr}\eta_{,zr}dx.\cr}
\leqno(6.17)
$$
Integrating by parts in the first integral on the l.h.s. yields
$$\eqal{
&\intop_S\bar n\cdot\nabla\eta_{,rz}\eta_{,rz}dS-\intop_\Omega
|\nabla\eta_{,rz}|^2dx-3\intop_\Omega{\eta_{,zr}^2\over r^2}dx\cr
&\quad+\intop_\Omega(\eta_{,zr}^2)_{,r}drdz=\intop_\Omega(\vartheta_{,r}
\eta_{,rz})_{,z}dx-\intop_\Omega\vartheta_{,r}\eta_{,rzz}dx.\cr}
\leqno(6.18)
$$
The first integral on the r.h.s. vanishes because $\vartheta_{,r}|_{S_2}=0$
and the first integral on the l.h.s. equals
$$
\intop_{S_1}\eta_{,zrr}\eta_{,zr}dS_1+\intop_{S_2}\eta_{,zrz}\eta_{,zr}dS_2
\equiv I_1.
$$
Since $\eta_{,rr}=-{3\eta_{,r}\over r}$ on $S_1$, so
$\eta_{,rrz}=-{3\eta_{,rz}\over r}$ on $S_1$ also. Therefore, the first
integral in $I_1$ takes the form
$$
-3\intop_{S_1}{\eta_{,rz}^2\over r}rdz=-3\intop_{-a}^a\eta_{,rz}^2|_{r=R}dz.
$$
Projecting (6.11) on $S_2$ gives $\eta_{,zz}|_{S_2}=0$, so also
$\eta_{,zzr}|_{S_2}=0$. Therefore, the second term in $I_1$ vanishes.

\noindent
In view of the above considerations, (6.18) takes the form
$$\eqal{
&\intop_\Omega|\nabla\eta_{,zr}|^2dx+3\intop_{-a}^a\eta_{,rz}^2|_{r=R}dz+
3\intop_\Omega{\eta_{,zr}^2\over r^2}dx\cr
&\quad-\intop_\Omega(\eta_{,zr}^2)_{,r}drdz=\intop_\Omega\vartheta_{,r}
\eta_{,zzr}dx.\cr}
\leqno(6.19)
$$
Performing integration by parts in the last term on the l.h.s. and applying
the Cauchy and the Young inequalities to the r.h.s. we derive
$$\eqal{
&\intop_\Omega|\nabla\eta_{,zr}|^2dx+6\intop_\Omega{\eta_{,zr}^2\over r^2}dx+
4\intop_{-a}^a\eta_{,rz}^2|_{r=R}dz\cr
&\quad+2\intop_{-a}^a\eta_{,rz}^2|_{r=0}dz\le\intop_\Omega\vartheta_{,r}^2dx.
\cr}
\leqno(6.20)
$$
Differentiating $(6.7)_1$ twice with respect to $z$, multiplying the result
by $\eta_{,zz}$ and integrating over $\Omega$ we arrive to
$$
\intop_\Omega\Delta\eta_{,zz}\eta_{,zz}dx+2\intop_\Omega\eta_{,zzr}\eta_{,zz}dx
=\intop_\Omega\vartheta_{,zz}\eta_{,zz}dx.
$$
Since $\eta_{,zz}$ vanishes on $S_1\cup S_2$ the above equality takes the form
$$\eqal{
&\intop_\Omega|\nabla\eta|_{,zz}|^2dx+\intop_{-a}^a\eta_{,zz}^2|_{r=0}dz=
-\intop_\Omega(\vartheta_{,z}\eta_{,zz})_{,z}dx\cr
&\quad+\intop_\Omega\vartheta_{,z}\eta_{,zzz}dx.\cr}
\leqno(6.21)
$$
Since $\eta_{,zz}|_{S_2}=0$ the first integral on the r.h.s. of (6.21)
vanishes. Applying the Cauchy and the Young inequalities to the second term
on the r.h.s. of (6.21) we derive
$$
\intop_\Omega|\nabla\eta_{,zz}|^2dx+\intop_{-a}^a\eta_{,zz}^2|_{r=0}dz\le
c\intop_\Omega\vartheta_{,z}^2dx.
\leqno(6.22)
$$
From (6.9), (6.14), (6.15), (6.20) and (6.22) we obtain (6.8). This concludes
the proof.

Estimate (6.8) does not contain the norm $\|\nabla\eta_{,rr}\|_{L_2(\Omega)}$
because to estimate it we need vanishing of $\eta_{,rrr}|_{S_1}$. But the
boundary conditions on $S_1$ do not imply it. Therefore we recall a smooth
cut-off function $\zeta=\zeta_1(r)$ such that $\zeta_1(r)=1$ for $r\le r_0$
and $\zeta_1(r)=0$ for $r\ge2r_0$, $2r_0<R$. Introducing the notation
$$
\tilde\eta=\eta\zeta_1^2,\quad \tilde\vartheta=\vartheta\zeta_1^2
\leqno(6.23)
$$
we see that $\tilde\eta$ is a solution to the problem
$$\eqal{
&\Delta\tilde\eta+{2\tilde\eta_{,r}\over r}=\tilde\vartheta+2\nabla\eta
\nabla\zeta_1^2+\eta\Delta\zeta_1^2+{\eta\zeta_{1,r}^2\over r}\equiv
\tilde\vartheta+\vartheta_1\equiv\vartheta_2,\cr
&\tilde\eta|_{S_2}=0.\cr}
\leqno(6.24)
$$

\proclaim Lemma 6.2.
Let $\tilde\eta$ be a solution to (6.24). Let $\vartheta_{,r}\in L_2(\Omega)$,
$v'\in H^1(\Omega)$.\\
Then the following estimate holds
$$\eqal{
&\intop_\Omega|\nabla\tilde\eta_{,rr}|^2r^2dx+\intop_\Omega|\tilde\eta_{,rr}|^2
dx+\intop_\Omega|\tilde\eta_{,rr}|^2dx+\intop_{-a}^a|\tilde\eta_{,r}|^2dz\cr
&\le c\intop_\Omega\vartheta_{,r}^2dx+c\|v'\|_{H^1(\Omega)}^2.\cr}
\leqno(6.25)
$$

\Proof
Differentiating (6.24) twice with respect to $r$, multiplying the result by
$r^2\tilde\eta$ and integrating over $\Omega$ yields
$$\eqal{
&\intop_\Omega(\Delta\tilde\eta)_{,rr}r^2\tilde\eta_{,rr}dx+2\intop_\Omega
\bigg({\tilde\eta_{,r}\over r}\bigg)_{,rr}r^2\tilde\eta_{,rr}dx\cr
&=\intop_\Omega\vartheta_{2,rr}r^2\tilde\eta_{,rr}dx.\cr}
\leqno(6.26)
$$
Since $\Delta\tilde\eta=\tilde\eta_{,rr}+\tilde\eta_{,zz}+
{\tilde\eta_{,r}\over r}$ we have
$$\eqal{
(\Delta\tilde\eta)_{,rr}&=\tilde\eta_{,rrrr}+\tilde\eta_{,zzrr}+
{\tilde\eta_{,rrr}\over r}-2{\tilde\eta_{,rr}\over r^2}+
2{\tilde\eta_{,r}\over r^3}\cr
&=\Delta\tilde\eta_{,rr}-2{\tilde\eta_{,rr}\over r^2}+
2{\tilde\eta_{,r}\over r^3}.\cr}
$$
Employing the expression in (6.26) implies
$$\eqal{
&\intop_\Omega\bigg(\Delta\tilde\eta_{,rr}-{2\tilde\eta_{,rr}\over r^2}+
2{\tilde\eta_{,r}\over r^3}\bigg)r^2\tilde\eta_{,rr}dx\cr
&\quad+2\intop_\Omega\bigg({\tilde\eta_{,r}\over r}\bigg)_{,rr}r^2
\tilde\eta_{,rr}dx=\intop_\Omega\vartheta_{2,rr}r^2\tilde\eta_{,rr}dx,\cr}
\leqno(6.27)
$$
where $\tilde\eta_{,rr}|_{S_2}=0$ and $\tilde\eta$ vanishes with all
derivatives on $S_1$.

\noindent
Integrating by parts in (6.27) and using the boundary conditions we arrive
to the equality
$$\eqal{
&\intop_\Omega\nabla\tilde\eta_{,rr}\nabla(r^2\tilde\eta_{,rr})dx+
6\intop_\Omega|\tilde\eta_{,rr}|^2dx\cr
&\quad-6\intop_\Omega\tilde\eta_{,rr}\tilde\eta_{,r}drdz-
2\intop_\Omega\tilde\eta_{,rrr}\tilde\eta_{,rr}rdx=\intop_\Omega
\vartheta_{2,rr}r^2\tilde\eta_{,rr}dx.\cr}
\leqno(6.28)
$$
Continuing calculations in (6.28) gives
$$\eqal{
&\intop_\Omega|\nabla\tilde\eta_{,rr}|^2r^2dx+2\intop_\Omega\tilde\eta_{,rrr}
\tilde\eta_{,rr}rdx+6\intop_\Omega|\tilde\eta_{,rr}|^2dx\cr
&\quad-6\intop_\Omega\tilde\eta_{,rr}\tilde\eta_{,r}drdz-2\intop_\Omega
\tilde\eta_{,rrr}\tilde\eta_{,rr}rdx=\intop_\Omega\vartheta_{2,rr}r^2
\tilde\eta_{,rr}dx.\cr}
$$
Next, we have
$$\eqal{
&\intop_\Omega|\nabla\tilde\eta_{,rr}|^2r^2dx+6\intop_\Omega
|\tilde\eta_{,rr}|^2dx-3\intop_\Omega\partial_r|\tilde\eta_{,r}|^2drdz\cr
&=\intop_\Omega(\vartheta_{2,r}r^2\tilde\eta_{,rr})_{,r}dx-\intop_\Omega
\vartheta_{2,r}(r^2\tilde\eta_{,rr})_{,r}dx.\cr}
$$
Since $\tilde\eta_{,r}|_{r=R}=0$, $\tilde\eta_{,rr}|_{r=R}=0$,
$\vartheta|_{r\le\varepsilon}=0$ we obtain
$$\eqal{
&\intop_\Omega|\nabla\tilde\eta_{,rr}|^2r^2dx+6\intop_\Omega
|\tilde\eta_{,rr}|^2dx+3\intop_{-a}^a|\tilde\eta_{,r}|^2|_{r=0}dz\cr
&=-\intop_\Omega(\vartheta_{2,r}\tilde\eta_{,rrr}r^2+2\vartheta_{2,r}
\tilde\eta_{,rr}r)dx.\cr}
\leqno(6.29)
$$
Applying the Cauchy and the Young inequalities to the r.h.s. of (6.29) implies
$$\eqal{
&\intop_\Omega|\nabla\tilde\eta_{,rr}|^2r^2dx+6\intop_\Omega|\tilde\eta_{,rr}|^2
dx+6\intop_{-a}^a|\tilde\eta_{,r}|^2|_{r=0}dz\cr
&\le c\intop_\Omega\vartheta_{2,r}^2dx\le c\intop_\Omega\vartheta_{,r}^2dx+
c\intop_\Omega\vartheta_{1,r}^2dx.\cr}
\leqno(6.30)
$$
The second integral on the r.h.s. of (6.30) will be estimated in the following
way
$$\eqal{
&\intop_\Omega\vartheta_{1,r}^2dx=\intop_\Omega\bigg(2\nabla\eta\nabla\zeta_1+
\eta\Delta\zeta_1+{\eta\zeta_{1,r}\over r}\bigg)_{,r}^2r^2dx\cr
&\le c\intop_{\Omega_{\zeta_{1,r}}}(\eta_{,rr}^2+\eta_{,r}^2+\eta^2)dx
\le c\intop_{\Omega_{\zeta_{1,r}}}(\psi_{,rr}^2+\psi_{,r}^2+\psi^2)dx\cr
&\le c\intop_{\Omega_{\zeta_{1,r}}}(v_{z,r}^2+v_z^2)dx+c
\intop_{\Omega_{\zeta_{1,r}}}\psi^2dx\equiv I.\cr}
$$
Since $\psi|_{S_2}=0$ we have
$$
\psi(r,z,t)=\intop_{-a}^z\psi_{,z}dz
$$
and the second term in $I$ is estimated by
$$
c\intop_\Omega v_r^2dx.
$$
Summarizing,
$$
\intop_\Omega\vartheta_{1,r}^2dx\le c\|v'\|_{H^1(\Omega)}^2.
\leqno(6.31)
$$
From (6.30) and (6.31) we obtain (6.25). This concludes the proof.

\noindent
To find an estimate for $\eta_{,rrr}$ in a neighborhood located in a positive
distance from the axis of symmetry we introduce the notation
$$
\bar\eta=\eta\zeta_2^2,\quad \bar\vartheta=\vartheta\zeta_2^2.
\leqno(6.32)
$$
Then (6.11) takes the form
$$
\bar\eta_{,rr}+\bar\eta_{,zz}+3{\bar\eta_{,r}\over r}=\bar\vartheta+
2\eta_{,r}\zeta_{2,r}^2+\eta\zeta_{2,rr}^2+3{\eta\zeta_{2,r}^2\over r}.
\leqno(6.33)
$$

\proclaim Lemma 6.3.
Let $\vartheta\in H^1(\Omega)$. Then
$$
\|\bar\eta_{,rrr}\|_{L_2(\Omega)}\le c\intop_\Omega(\vartheta^2+
\vartheta_{,r}^2+\vartheta_{,z}^2)dx.
\leqno(6.34)
$$

\Proof
Differentiating (6.33) with respect to $r$ yields
$$\eqal{
&\bar\eta_{,rrr}+\bar\eta_{,rzz}+3{\bar\eta_{,rr}\over r}-
3{\bar\eta_{,r}\over r^2}=\bar\vartheta_{,r}+\bigg(2\eta_{,r}\zeta_{2,r}^2+
\eta\zeta_{2,rr}^2\cr
&\quad+3{\eta\zeta_{2,r}^2\over r}\bigg)_{,r}.\cr}
\leqno(6.34)
$$
In view of (6.8) we obtain from (6.34) the inequality
$$\eqal{
&\|\bar\eta_{,rrr}\|_{L_2(\Omega)}^2\le c\|\bar\eta_{,rzz}\|_{L_2(\Omega)}^2+
c(1/r_0)\|\bar\eta_{,rr}\|_{L_2(\Omega)}^2\cr
&\quad+c(1/r_0)\|\eta_{,r}\|_{L_2(\Omega)}^2+c
\|\bar\vartheta_{,r}\|_{L_2(\Omega)}^2+c(1/r_0)
(\|\eta_{,rr}\|_{L_2(\Omega_{\zeta_2})}^2\cr
&\quad+\|\eta_{,r}\|_{L_2(\Omega_{\zeta_2})}^2+
\|\eta\|_{L_2(\Omega_{\zeta_2})}^2)\le c(1/r_0)\intop_\Omega
(\vartheta^2+\vartheta_{,r}^2+\vartheta_{,z}^2)dx.\cr}
\leqno(6.35)
$$
This concludes the proof.

However, the norm on the l.h.s. of (6.5) is over $\Omega_\varepsilon$ we can
pass to the limit $\varepsilon=0$ because the r.h.s. of (6.5) is independent
of $\varepsilon$.

\proclaim Lemma 6.4.
Assume that ${\chi\over r}\in V_2^0(\Omega^T)$. Then
$$
\bigg\|{v_r\over r}\bigg\|_{V_2^1(\Omega^T)}\le c
\bigg\|{\chi\over r}\bigg\|_{V_2^0(\Omega^T)}.
\leqno(6.36)
$$

\Proof
From (6.8) we have
$$
\intop_\Omega(\eta_{,z}^2+\eta_{,zrr}^2+\eta_{,zzr}^2+\eta_{,zzz}^2)dx\le c
\bigg\|{\chi\over r}\bigg\|_{H^1(\Omega)}^2.
\leqno(6.37)
$$
Using that
$$\eqal{
&\intop_\Omega\eta_{,z}^2dx=\intop_\Omega\bigg|{v_r\over r}\bigg|^2dx,\cr
&\intop_\Omega(|\eta_{,zrr}|^2+|\eta_{,zzr}|^2+|\eta_{,zzz}|^2)dx\cr
&=\intop_\Omega\bigg[\bigg({v_r\over r}\bigg)_{,rr}^2+
\bigg({v_r\over r}\bigg)_{,zr}^2+\bigg({v_r\over r}\bigg)_{,zz}^2\bigg]dx\cr}
$$
we obtain from (6.37) the estimate
$$
\bigg\|{v_r\over r}\bigg\|_{H^2(\Omega)}^2\le c
\bigg\|{\chi\over r}\bigg\|_{H^1(\Omega)}^2.
\leqno(6.38)
$$
Integrating (6.38) with respect to time implies
$$
\bigg\|{v_r\over r}\bigg\|_{L_2(0,T;H^2(\Omega))}^2\le c
\bigg\|{\chi\over r}\bigg\|_{L_2(0,T;H^1(\Omega))}^2.
\leqno(6.39)
$$
From (6.15) we have
$$
\intop_\Omega\bigg|\nabla{v_r\over r}\bigg|^2dx=\intop_\Omega
|\nabla\eta_{,z}|^2dx\le c\bigg\|{\chi\over r}\bigg\|_{L_2(\Omega)}^2
\leqno(6.40)
$$
and (6.9) yields
$$
\intop_\Omega\bigg|{v_r\over r}\bigg|^2dx=\intop_\Omega
\bigg({\psi\over r^2}\bigg)_{,z}^2dx\le\intop_\Omega|\nabla\eta|^2dx\le
c\bigg\|{\chi\over r}\bigg\|_{L_2(\Omega)}^2.
$$
The above two estimates yield
$$
\bigg\|{v_r\over r}\bigg\|_{H^1(\Omega)}^2\le c
\bigg\|{\chi\over r}\bigg\|_{L_2(\Omega)}^2.
\leqno(6.41)
$$
In view of the assumptions of the lemma estimate (6.41) implies
$$
\bigg\|{v_r\over r}\bigg\|_{L_\infty(0,T;H^1(\Omega))}^2\le c
\bigg\|{\chi\over r}\bigg\|_{L_\infty(0,T;L_2(\Omega))}^2.
\leqno(6.42)
$$
From (6.39) and (6.42) we derive (6.36). This concludes the proof.

\noindent
Next we have

\proclaim Lemma 6.5.
Assume that ${\chi\over r}\in V_2^0(\Omega^T)$. Then
$$
\|v'\|_{V_2^1(\Omega^T)}\le c\bigg\|{\chi\over r}\bigg\|_{V_2^0(\Omega^T)}.
\leqno(6.43)
$$

\Proof
Since $v_r={\psi_{,z}\over r}=r\eta_{,z}$, $v_z=-{\psi_{,r}\over r}=
-r\eta_{,r}-2\eta$ we have
$$\eqal{
&\intop_\Omega(|\nabla v_r|^2+|\nabla v_z|^2+|v_r|^2+|v_z|^2)dx\le c
\intop_\Omega(|\eta|^2+|\nabla\eta|^2\cr
&\quad+|\nabla\eta_{,r}|^2+|\nabla\eta_{,z}|^2)dx\le c
\bigg\|{\chi\over r}\bigg\|_{L_2(\Omega)}^2,\cr}
\leqno(6.44)
$$
where (6.8) was used. Similarly,
$$\eqal{
&\intop_\Omega(|\nabla^2(v_r)|^2+|\nabla^2(v_z)|^2+|v_r|^2+|v_z|^2)dx\cr
&\le c\intop_\Omega(|\nabla^2(r\eta_{,z})|^2+|\nabla^2(r\eta_{,r})|^2+
|\nabla\eta|^2+|\eta|^2)dx\cr
&\le c\bigg\|{\chi\over r}\bigg\|_{H^1(\Omega)}^2,\cr}
\leqno(6.45)
$$
where Lemmas 6.1, 6.2, 6.3 were employed.
Taking $L_\infty$ norm with respect to time to (6.44) and $L_2$ norm with
respect to time to (6.45) we obtain (6.43). This concludes the proof.

\section{7. Estimate for the angular component of velocity}

Let us consider the problem
$$\eqal{
&v_{\varphi,t}-\nu\Delta v_\varphi+v'\cdot\nabla v_\varphi+{v_r\over r}
v_\varphi+\nu{v_\varphi\over r^2}=f_\varphi\quad &{\rm in}\ \ \Omega^T,\cr
&v_{\varphi,r}={1\over r}v_\varphi\quad &{\rm on}\ \ S_1^T,\cr
&v_{\varphi,z}=0\quad &{\rm on}\ \ S_2^T,\cr
&v_\varphi|_{t=0}=v_\varphi(0)\quad &{\rm in}\ \ \Omega.\cr}
\leqno(7.1)
$$

\proclaim Lemma 7.1.
Assume that $v_\varphi(0)\in H_0^1(\Omega)$,
$v_\varphi\in L_{5/2}(0,T;W_{5/2}^1(\Omega))$, $f_\varphi\in L_2(\Omega^T)$.
Assume that $A$ defined by (6.2) and $A_0$ defined by (6.4) are finite. Let
$A_*$ be introduced in (6.5). Then
$$\eqal{
&{1\over4}\|v_{\varphi,t}\|_{L_2(\Omega^t)}^2+{\nu\over2}
\|v_\varphi(t)\|_{H_0^1(\Omega)}^2+{\nu\over2}
\|v_\varphi(R,t)\|_{L_2(-a,a)}^2\cr
&\le{\nu\over2}\|v_\varphi(0)\|_{H_0^1(\Omega)}^2+cd_6^2+\varphi(A_*)
\|v_\varphi\|_{L_{5/2}(0,t;W_{5/2}^1(\Omega))}^2\cr
&\quad+\|f_\varphi\|_{L_2(\Omega^t)}^2,\quad t\le T,\cr}
\leqno(7.2)
$$
where
$$
\|u\|_{H_0^1(\Omega)}=\bigg(\intop_\Omega(|\nabla u|^2+{u^2\over r^2}\bigg)
dx\bigg)^{1/2}.
$$

\Proof
Multiplying $(7.1)_1$ by $v_{\varphi,t}$ and integrating the result over
$\Omega$ yields
$$\eqal{
&\intop_\Omega v_{\varphi,t}^2dx+{\nu\over2}{d\over dt}\intop_\Omega
|\nabla v_\varphi|^2dx+{\nu\over2}{d\over dt}\intop_\Omega
{v_\varphi^2\over r^2}dx\cr
&\quad-{\nu\over2}{d\over dt}\intop_{-a}^av_\varphi^2|_{r=R}dz\le
{\varepsilon_1\over2}\intop_\Omega v_{\varphi,t}^2dx+{1\over2\varepsilon_1}
\intop_\Omega|v'\cdot\nabla v_\varphi|^2dx\cr
&\quad+{\varepsilon_2\over2}\intop_\Omega v_{\varphi,t}^2dx+
{1\over2\varepsilon_2}\intop_\Omega{v_r^2\over r^2}v_\varphi^2dx+
{\varepsilon_3\over2}\intop_\Omega v_{\varphi,t}^2dx+{1\over2\varepsilon_3}
\intop_\Omega f_\varphi^2dx.\cr}
\leqno(7.3)
$$
Setting $\varepsilon_1=\varepsilon_2=\varepsilon_3={1\over2}$ we get
$$\eqal{
&{1\over4}\intop_\Omega v_{\varphi,t}^2dx+{\nu\over2}{d\over dt}\intop_\Omega
|\nabla v_\varphi|^2dx+{\nu\over2}{d\over dt}\intop_\Omega
{v_\varphi^2\over r^2}dx\cr
&\quad-{\nu\over2}{d\over dt}\intop_{-a}^av_\varphi^2|_{r=R}dz\le\intop_\Omega
|v'\cdot\nabla v_\varphi|^2dx+\intop_\Omega{v_r^2\over r^2}v_\varphi^2dx\cr
&\quad+\intop_\Omega f_\varphi^2dx\cr}
\leqno(7.4)
$$
Integrating (7.4) with respect to time and using that
$$
\|v'\|_{L_{10}(\Omega^T)}+\bigg\|{v'\over r}\bigg\|_{L_{10}(\Omega^t)}\le
\varphi(A_*),
$$
where Lemmas 6.4 and 6.5 are employed, we obtain from (7.4) the inequality
$$\eqal{
&{1\over4}\intop_{\Omega^t}v_{\varphi,t'}^2dxdt'+{\nu\over2}\intop_\Omega
|\nabla v_\varphi(t)|^2dx+{\nu\over2}\intop_\Omega{v_\varphi^2(t)\over r^2}dx
\cr
&\quad+{\nu\over2}\intop_{-a}^av_\varphi^2|_{r=R,t=0}dz\le{\nu\over2}
\intop_\Omega|\nabla v_\varphi(0)|^2dx+{\nu\over2}\intop_\Omega
{v_\varphi^2(0)\over r^2}dx\cr
&\quad+{\nu\over2}\intop_{-a}^av_\varphi^2|_{r=R,t}dz+\varphi(A_*)
(\|\nabla v_\varphi\|_{L_{5/2}(\Omega^t)}^2+
\|v_\varphi\|_{L_{5/2}(\Omega^t)}^2)\cr
&\quad+\|f_\varphi\|_{L_2(\Omega^t)}^2,\quad t\le T.\cr}
\leqno(7.5)
$$
Using that the third integral on the r.h.s. of (7.5) is bounded by $cd_6$
(see (3.1)) we derive (7.2). This concludes the proof.

To estimate the coefficient near $\varphi(A_*)$ on the r.h.s. of (7.2) we
introduce the Green function to the linear part of problem (7.1). Let us
denote it by $G$.
To obtain an estimate independent of time we consider problem (7.1) in the
intervals $(0,T_0)$ and $(kT_0,(k+1)T_0)$, where $k\in\N$ and $T_0$ is a given
positive number. For $t\in(0,T_0)$ problem (7.1) can be expressed in the
following integral form
$$\eqal{
&v_\varphi(x,t)=\intop_{\Omega^t}\nabla_{y_\alpha}G(x-y,t-\tau)v'_\alpha
(y,\tau)v_\varphi(y,\tau)dyd\tau\cr
&\quad-\intop_{\Omega^t}G(x-y,t-\tau)\bigg({v_r(y,\tau)\over r}
v_\varphi(y,\tau)+\nu{v_\varphi(y,\tau)\over r^2}\bigg)dyd\tau\cr
&\quad+\intop_{\Omega^t}G(x-y,t)v_\varphi(y,0)dy+\intop_{S_1^t}
G(x-z,t-\tau){1\over r}v_\varphi(z,\tau)dzd\tau\cr
&\quad+\intop_{\Omega^t}G(x-y,t-\tau)f_\varphi(y,\tau)dyd\tau,\quad
t\in(0,T_0).\cr}
\leqno(7.6)
$$
To examine problem (7.1) in the interval $(kT_0,(k+1)T_0)$, $k\in\N$, we
introduce a smooth cut-off function $\zeta=\zeta(t)$ such that $\zeta(t)=1$
for $t\in[kT_0,(k+1)T_0]$ and $\zeta(t)=0$ for
$t\not\in[(k-1/2)T_0,(k+3/2)T_0]$.

\noindent
Multiplying (7.1) by $\zeta$ and introducing the notation
$\tilde v_\varphi=v_\varphi\zeta$, $\tilde f_\varphi=f_\varphi\zeta$ we obtain
$$\eqal{
&\tilde v_{\varphi,t}-\nu\Delta\tilde v_\varphi+v'\cdot\nabla\tilde v_\varphi+
{v_r\over r}\tilde v_\varphi\cr
&\quad=\tilde f_\varphi+
v_\varphi\dot\zeta\quad &{\rm in}\ \ \Omega\times((k-1)T_0,(k+2)T_0),\cr
&\tilde v_{\varphi,r}={1\over r}\tilde v_\varphi\quad &{\rm on}\ \
S_1\times((k-1)T_0,(k+2)T_0),\cr
&\tilde v_{\varphi,z}=0\quad &{\rm on}\ \ S_2\times((k-1)T_0,(k+2)T_0),\cr
&\tilde v_\varphi|_{t=(k-1)T_0}=0.\cr}
\leqno(7.7)
$$
Using the Green function we express (7.7) in the following integral equation
$$\eqal{
&\tilde v_\varphi(x,t)=\intop_{(k-1)T_0}^t\intop_\Omega\nabla_{y_\alpha}G
(x-y,t-\tau)v'_\alpha(y,\tau)\tilde v_\varphi(y,\tau)dyd\tau\cr
&\quad-\intop_{(k-1)T_0}\intop_GG(x-y,t-\tau)\bigg({v_r(y,\tau)\over r}
\tilde v_\varphi(y,\tau)+\nu{\tilde v_\varphi(y,\tau)\over r^2}\bigg)dyd\tau\cr
&\quad+\intop_{(k-1)T_0}^t\intop_{S_1}G(x-z,t-\tau){1\over r}\tilde v_\varphi
(z,\tau)dzd\tau\cr
&\quad+\intop_{(k-1)T_0}^t\intop_\Omega G(x-y,t-\tau)(\tilde f_\varphi(y,\tau)
+v_\varphi(y,\tau)\dot\zeta(\tau))dyd\tau,\cr}
\leqno(7.8)
$$
where $t\in((k-1)T_0,(k+2)T_0)$.

\proclaim Lemma 7.2.
Assume that $T_0>0$ is a given positive number. Assume that
$v_\varphi(0)\in W_{5/2}^{1/5}(\Omega)$,
$f_\varphi\in L_{5/3}(\Omega\times((k-1)T_0,(k+1)T_0))$,
${v_\varphi\over r}\in L_{10/3}(\Omega\times((k-1)T_0,(k+1)T_0))$ for
$k\in\N$. Then the following estimates hold
$$\eqal{
&\|v_\varphi\|_{W_{5/2}^{1,1/2}(\Omega^t)}\le c(T_0)\bigg[\varphi(A_*)d_7+
T_0^{2/5}d_6+\|v_\varphi(0)\|_{W_{5/2}^{1/5}(\Omega)}\cr
&\quad\bigg\|{v_\varphi\over r}\bigg\|_{L_{10/3}(\Omega^{T_0})}^2+
\|f_\varphi\|_{L_{5/3}(\Omega^{T_0})}\bigg],\quad t\le T_0,\cr}
\leqno(7.9)
$$
and
$$\eqal{
&\|v_\varphi\|_{W_{5/2}^{1,1/2}(\Omega\times(kT_0,t))}\le c(T_0)
\bigg[\varphi(A_*)d_7+T_0^{2/5}d_6\cr
&\quad+
\bigg\|{v_\varphi\over r}\bigg\|_{L_{10/3}(\Omega\times((k-1)T_0,(k+1)T_0))}^2
+\|f_\varphi\|_{L_{5/3}(\Omega\times((k-1)T_0,(k+1)T_0))}\cr
&\quad+cT_0^{-13/20}d_6^{1/2}d_7^{1/2}\bigg],\quad t\in(kT_0,(k+1)T_0),\cr}
\leqno(7.10)
$$
where $d_6$ is introduced in (3.1), $d_7$ in (2.15) and $A_*$ is defined by
(6.5). Moreover, the last term on the r.h.s. follows from estimate (3.6).

\Proof
First we consider the case $t\in(0,T_0)$. Applying the potential theory
to (7.6) yields
$$\eqal{
&\|v_\varphi\|_{W_\sigma^{1,1/2}(\Omega^t)}\le c(T_0)
[\|v'v_\varphi\|_{L_\sigma(\Omega^{T_0})}\cr
&\quad+
\bigg\|{v_r\over r}v_\varphi\bigg\|_{L_{5\sigma\over5+\sigma}(\Omega^{T_0})}+
\bigg\|{v_\varphi\over r^2}\bigg\|_{L_{5\sigma\over5+\sigma}(\Omega^{T_0})}+
\|v_\varphi\|_{L_\sigma(S_1^{T_0})}\cr
&\quad+\|v_\varphi(0)\|_{W_\sigma^{1-2/\sigma}(\Omega)}+
\|f_\varphi\|_{L_{5\sigma\over5+\sigma}(\Omega^{T_0})}],\quad t\le T_0.\cr}
\leqno(7.11)
$$
To estimate the third term on the r.h.s. of (7.2) we assume that $\sigma=5/2$.
Then (7.11) takes the form
$$\eqal{
&\|v_\varphi\|_{W_{5/2}^{1,1/2}(\Omega^t)}\le c(T_0)
\bigg[\|v'v_\varphi\|_{L_{5/2}(\Omega^{T_0})}+
\bigg\|{v_r\over r}v_\varphi\bigg\|_{L_{5/3}(\Omega^{T_0})}\cr
&\quad+\bigg\|{v_\varphi\over r^2}\bigg\|_{L_{5/3}(\Omega^{T_0})}+
\|v_\varphi\|_{L_{5/2}(S_1^{T_0})}+\|v_\varphi(0)\|_{W_{5/2}^{1/5}(\Omega)}\cr
&\quad+\|f_\varphi\|_{L_{5/3}(\Omega^{T_0})}\bigg],\quad t\le T_0.\cr}
\leqno(7.12)
$$
Now we estimate the particular terms from the r.h.s. of (7.12). Applying the
H\"older inequality the first term on the r.h.s. is bounded by
$$
\|v'v_\varphi\|_{L_{5/2}(\Omega^{T_0})}\le\|v'\|_{L_{10}(\Omega^{T_0})}
\|v_\varphi\|_{L_{10/3}(\Omega^{T_0})}\le\varphi(A_*)d_7,
$$
where we used (2.7), (2.15), (6.5) and (6.43).

\noindent
Similarly, the second term on the r.h.s. of (7.12) is bounded by
$$
\bigg\|v_r{v_\varphi\over r}\bigg\|_{L_{5/3}(\Omega^{T_0})}\le
\|v_r\|_{L_{10}(\Omega^{T_0})}
\bigg\|{v_\varphi\over r}\bigg\|_{L_2(\Omega^{T_0})}\le\varphi(A_*)d_7.
$$
Assuming that $v_\varphi\ge1$ (otherwise we have regularity of axially
symmetric solutions) the third term on the r.h.s. of (7.12) is estimated by
$$
\bigg\|{v_\varphi\over r}\bigg\|_{L_{10/3}(\Omega^{T_0})}^2.
$$
The fourth term is bounded by
$$
cT_0^{2/5}\|u\|_{L_\infty(\Omega^{T_0})}\le cT_0^{2/5}d_6,
$$
where (3.1) is used. Summarizing the above estimates we obtain (7.9).

Now we apply the potential theory to (7.8). Then we obtain
$$\eqal{
&\|v_\varphi\|_{W_{5/2}^{1,1/2}(\Omega\times(kT_0,t))}\le
\|\tilde v_\varphi\|_{W_{5/2}^{1,1/2}(\Omega\times((k-1)T_0,t))}\cr
&\le c(T_0)\bigg[\|v'\tilde v_\varphi\|_{L_{5/2}(\Omega\times
((k-1)T_0,(k+1)T_0))}\cr
&\quad+\bigg\|v_r{\tilde v_\varphi\over r}\bigg\|_{L_{5/3}(\Omega\times
((k-1)T_0,(k+1)T_0))}\cr
&\quad+\bigg\|{\tilde v_\varphi\over r^2}\bigg\|_{L_{5/3}(\Omega\times
((k-1)T_0,(k+1)T_0))}\cr
&\quad+\|\tilde v_\varphi\|_{L_{5/2}(S_1\times((k-1)T_0,(k+1)T_0))}\cr
&+\quad\|\tilde f_\varphi+v_\varphi\dot\zeta\|_{L_{5/3}(\Omega\times
((k-1)T_0,(k+1)T_0))}\bigg]\cr}
\leqno(7.13)
$$
where $t\in[kT_0,(k+1)T_0]$.
Repeating the considerations leading to (7.9) we obtain (7.10). This
concludes the proof.

To estimate the term
$\big\|{v_\varphi\over r}\big\|_{L_{10/3}(\Omega\times((k-1)T_0,(k+1)T_0))'}$
$k\in\N$, from the r.h.s. of (7.9) and (7.10) we introduce the quantity
$w={v_\varphi\over r}$ which is a solution to the problem
$$\eqal{
&w_{,t}+v'\cdot\nabla w+2{v_r\over r}w-\nu\Delta w-{2\nu\over r}w_{,r}=
{f_\varphi\over r}\equiv g_0\quad &{\rm in}\ \ \Omega\times\R_+,\cr
&w_{,r}=0\quad &{\rm on}\ \ S_1\times\R_+,\cr
&w_{,z}=0\quad &{\rm on}\ \ S_2\times\R_+,\cr
&w|_{t=0}=w(0)\quad &{\rm in}\ \ \Omega.\cr}
\leqno(7.14)
$$

\proclaim Lemma 7.3.
Assume that $A_*$ defined by (6.5) is finite. Let $g_0={f_\varphi\over r}$ and
$|||g_0|||_2\equiv\sup_{k\in\N_0}$ $\|g_0\|_{L_2(kT_0,(k+1)T_0;L_2(\Omega))}$
be finite, $k\in\N_0=\N\cup\{0\}$. Let $T_0>0$ be given. Let
$w(0)\in L_2(\Omega)$. Then
$$
\|w\|_{L_{10/3}(\Omega\times(kT_0,t))}\le\varphi(T_0,d_7,A_*,|||g_0|||_2)+
e^{-\nu_*kT_0}\|w(0)\|_{L_2(\Omega)},
\leqno(7.15)
$$
where $\varphi$ is an increasing positive function, $d_7$ is defined by
(2.15), $0<\nu_*<\nu$ and $t\in(kT_0;(k+1)T_0]$.

\Proof
Multiplying (7.14) by $w|w|^{s-2}$ and integrating over $\Omega$ we obtain
$$\eqal{
&{1\over s}{d\over dt}\intop_\Omega|w|^sdx+{4\nu(s-1)\over s^2}\intop_\Omega
|\nabla|w|^{s/2}|^2dx\cr
&\le2\intop_\Omega\bigg|{v_r\over r}\bigg||w|^sdx+{2\nu\over s}\intop_\Omega
\partial_r|w|^sdrdz+\intop_\Omega|g_0|\,|w|^{s-1}dx\cr}
\leqno(7.16)
$$
Performing integration in the second term on the r.h.s. of (7.16) yields
$$\eqal{
&{1\over s}{d\over dt}\intop_\Omega|w|^sdx+{4\nu(s-1)\over s^2}\intop_\Omega
|\nabla|w|^{s/2}|^2dx\cr
&\le2\intop_\Omega\bigg|{v_r\over r}\bigg||w|^sdx+{2\nu\over s}\intop_{-a}^a
|w(R,z,t)|^sdz\cr
&\quad+\|g_0\|_{L_s(\Omega)}\|w\|_{L_s(\Omega)}^{s-1}.\cr}
\leqno(7.17)
$$
In view of the interpolation inequality
$$
\intop_{-a}^a|w(R,z,t)|^sdz\le\varepsilon\intop_\Omega|\nabla|w|^{s/2}|^2dx+
{c_1\over\varepsilon}\intop_\Omega|w|^sdx
$$
with $\varepsilon={s-1\over s}$ we derive from (7.17) the inequality
$$\eqal{
&{d\over dt}\intop_\Omega|w|^sdx+{2\nu(s-1)\over s}\intop_\Omega
|\nabla|w|^{s/2}|^2dx\le2s\intop_\Omega\bigg|{v_r\over r}\bigg||w|^sdx\cr
&\quad+{c_1\over s-1}\intop_\Omega|w|^sdx+s\|g_0\|_{L_s(\Omega)}
\|w\|_{L_s(\Omega)}^{s-1}.\cr}
\leqno(7.18)
$$
We are going to consider (7.18) for $s\in\big[{3\over2},4\big]$. Then (7.18)
takes the form
$$\eqal{
&{d\over dt}\intop_\Omega|w|^sdx+{2\over3}\nu\intop_\Omega|\nabla|w|^{s/2}|^2dx
\le8\intop_\Omega\bigg|{v_r\over r}\bigg||w|^sdx\cr
&\quad+2c_1\intop_\Omega|w|^sdx+4\|g_0\|_{L_s(\Omega)}
\|w\|_{L_s(\Omega)}^{s-1}.\cr}
\leqno(7.19)
$$
Introducing the quantity
$$
\eta=|w|^{s/2}
\leqno(7.20)
$$
inequality (7.19) assumes the form
$$\eqal{
&{d\over dt}\intop_\Omega|\eta|^2dx+{2\over3}\nu\intop_\Omega|\nabla\eta|^2dx
\le8\intop_\Omega\bigg|{v_r\over r}\bigg|\eta^2dx\cr
&\quad+2c_1\intop_\Omega|\eta|^2dx+4\|\eta\|_{L_2(\Omega)}^{s-1}
\|g_0\|_{L_s(\Omega)}.\cr}
\leqno(7.21)
$$
Since $\eta_{,r}|_{S_1}=0$, $\eta_{,z}|_{S_2}=0$ the following Poincare
inequality holds
$$
\intop_\Omega\eta^2dx\le c_p\intop_\Omega|\nabla\eta|^2dx+c_p
\bigg|\intop_\Omega\eta dx\bigg|.
\leqno(7.22)
$$
Introducing the splitting ${2\over3}\nu=\nu_1+\nu_2$, $\nu_i>0$, $i=1,2$,
we obtain from (7.21) the inequality
$$\eqal{
&{d\over dt}\intop_\Omega\eta^2dx+{\nu_1\over c_p}\intop_\Omega\eta^2dx+
\nu_2\intop_\Omega|\nabla\eta|^2dx\le{\nu_1\over c_p}\bigg|\intop_\Omega
\eta dx\bigg|\cr
&\quad+8\intop_\Omega\bigg|{v_r\over r}\bigg|\eta^2dx+2c_1\intop_\Omega\eta^2dx
+4\|\eta\|_{L_2(\Omega)}^{s-1}\|g_0\|_{L_s(\Omega)}.\cr}
\leqno(7.23)
$$
Let $\nu_*={\nu_1\over c_p}$. Then from (7.23) we have
$$\eqal{
&{d\over dt}\bigg(\intop_\Omega\eta^2dxe^{\nu_*t}\bigg)+\nu_2\intop_\Omega
|\nabla\eta|^2dxe^{\nu_*t}\cr
&\le\bigg[\nu_*\intop_\Omega|\eta|dx+2c_1\intop_\Omega\eta^2dx+8\intop_\Omega
\bigg|{v_r\over r}\bigg|\eta^2dx+4\|\eta\|_{L_2(\Omega)}^{s-1}
\|g_0\|_{L_s(\Omega)}\bigg]e^{\nu_*t}.\cr}
\leqno(7.24)
$$
Integrating (7.24) with respect to time from $kT_0$ to $t\in(kT_0,(k+1)T_0]$,
$k\in\N\cup\{0\}\equiv\N_0$, yields
$$\eqal{
&\intop_\Omega\eta^2(t)dxe^{\nu_*t}+\nu_2\intop_{kT_0}^t\intop_\Omega
|\nabla\eta(t')|^2dxe^{\nu_*t'}dt'\cr
&\le\intop_{kT_0}^t\bigg[\nu_*\intop_\Omega|\eta|dx+2c_1\intop_\Omega\eta^2dx+
8\intop_\Omega\bigg|{v_r\over r}\bigg|\eta^2dx\cr
&\quad+4\|\eta\|_{L_2(\Omega)}^{s-1}
\|g_0\|_{L_s(\Omega)}\bigg]e^{\nu_*t'}dt'+\intop_\Omega\eta^2(kT_0)dxe^{\nu_*kT_0}.\cr}
\leqno(7.25)
$$
Continuing, we have
$$\eqal{
&\intop_\Omega\eta^2(t)dx+\nu_2e^{-\nu_*t}\intop_{kT_0}^t\intop_\Omega
|\nabla\eta(t')|^2dxe^{\nu_*t'}dt'\cr
&\le e^{-\nu_*t}\intop_{kT_0}^t\bigg[\nu_*\intop_\Omega|\eta|dx+
2c_1\intop_\Omega\eta^2dx+8\intop_\Omega\bigg|{v_r\over r}\bigg|\eta^2dx\cr
&\quad+4\|\eta\|_{L_2(\Omega)}^{s-1}\|g_0\|_{L_s(\Omega)}\bigg]e^{\nu_*t'}dt+
e^{-\nu_*(t-kT_0)}\intop_\Omega\eta^2(kT_0)dx,\cr
&t\in(kT_0,(k+1)T_0].\cr}
\leqno(7.26)
$$
First we obtain an estimate for $\intop_\Omega\eta^2(kT_0)dx$ for any
$k\in\N_0$. For this purpose we omit the second term on the l.h.s. of (7.26).
Therefore we consider the inequalities
$$\eqal{
&\intop_\Omega\eta^2((k+1)T_0)dx\cr
&\le\intop_{kT_0}^{(k+1)T_0}\bigg[\nu_*\intop_\Omega|\eta|dx+2c_1\intop_\Omega
\eta^2dx+8\intop_\Omega\bigg|{v_r\over r}\bigg|\eta^2dx\cr
&\quad+4\|\eta\|_{L_2(\Omega)}^{s-1}\|g_0\|_{L_s(\Omega)}\bigg]dt+
e^{-\nu_*T_0}\intop_\Omega\eta^2(kT_0)dx.\cr}
\leqno(7.27)
$$
Take $s={9\over5}$. Then the first integral on the r.h.s. of (7.27) is
estimated by
$$\eqal{
&c\intop_{kT_0}^{(k+1)T_0}dt\intop_\Omega(|w|^{9/10}+|w|^{9/5})dx\cr
&\quad+c
\bigg(\intop_{kT_0}^{(k+1)T_0}\intop_\Omega\bigg|{v_r\over r}\bigg|^{10}
dxdt\bigg)^{1/10}\bigg(\intop_{kT_0}^{(k+1)T_0}\intop_\Omega|w|^2
dxdt\bigg)^{9/10}\cr
&\quad+cT_0^{2/45}\bigg(\intop_{kT_0}^{(k+1)T_0}
\intop_\Omega|w|^{9/5}dxdt\bigg)^{2/5}
\bigg(\intop_{kT_0}^{(k+1)T_0}\intop_\Omega|g_0|^{9/5}dxdt\bigg)^{5/9}
\equiv I_1.\cr}
$$
Applying the H\"older inequality and the energy estimate (2.15) with $d_7$
implies
$$\eqal{
I_1&\le\varphi(|\Omega|,A_0)[T_0^{11/20}d_7^{9/10}+T_0^{1/10}d_7^{9/5}+
\varphi(A_*)d_7^{9/10}\cr
&\quad+T_0^ad_7^b\|g_0\|_{L_{9/5}(\Omega\times(kT_0,(k+1)T_0))}]\cr
&\equiv A_1(T_0,d_7,A_*,\sup_k\|g\|_{L_{9/5}(\Omega\times(kT_0,(k+1)T_0))}),
\quad a>0,\ \ b>0.\cr}
$$
Introducing the notation
$$
X_1(t)=\intop_\Omega|w(t)|^{9/5}dx
\leqno(7.28)
$$
we obtain from (7.27) the inequality
$$
X_1((k+1)T_0)\le A_1+e^{-\nu_*T_0}X_1(kT_0).
\leqno(7.29)
$$
Hence, (7.29) implies
$$
X_1(kT_0)\le{A_1\over1-e^{-\nu_*T_0}}+e^{-k\nu_*T_0}X_1(0).
\leqno(7.30)
$$
In view of (7.30) we can consider (7.26) in the any interval
$(kT_0,(k+1)T_0)$, $k\in\N_0$, for
$$
\eta_1=|w|^{9/10}.
\leqno(7.31)
$$
Then we obtain the estimate
$$\eqal{
&\|\eta_1\|_{V_2^0(\Omega\times(kT_0,t))}^2\le(1+e^{\nu_*T_0})
\bigg[A_1(T_0)+{A_1(T_0)\over1-e^{-\nu_*T_0}}\cr
&\quad+e^{-k\nu_*T_0}\|\eta_1(0)\|_{L_2(\Omega)}^2\bigg],\quad
t\in(kT_0,(k+1)T_0].\cr}
\leqno(7.32)
$$
Using that
$$
\intop_\Omega\eta^2(0)dx=\intop|w(0)|^sdx
$$
and
$$
\|\eta\|_{V_2^0(\Omega^T)}^2=\|\,|w|^{s/2}\|_{V_2^0(\Omega^T)}^2\ge
\|w\|_{L_{{5\over3}s}(\Omega^T)}^s
$$
we obtain from (7.32) in the case $s={9\over5}$ the estimate
$$\eqal{
&\|w\|_{L_3(\Omega\times(kT_0,t))}\le(1+e^{\nu_*T_0})^{5/9}
\bigg[{2-e^{-\nu_*T_0}\over1-e^{-\nu_*T_0}}A_1(T_0)\cr
&\quad+e^{-k\nu_*T_0}\|w(0)\|_{L_{9/5}(\Omega)}^{9/5}\bigg]^{5/9}\equiv
A_2(T_0),\quad t\in(kT_0,(k+1)T_0],\cr}
\leqno(7.33)
$$
where $A_2$ depends also on the same quantities as $A_1$.

\noindent
Next we consider the case $s=2$. Then the first term on the r.h.s. of (7.27)
is estimated by
$$\eqal{
&c\bigg[\intop_{kT_0}^{(k+1)T_0}dt\intop_\Omega(|w|+|w|^2)dx\cr
&\quad+\bigg(\intop_{kT_0}^{(k+1)T_0}\intop_\Omega\bigg|{v_r\over r}\bigg|^{10}
dxdt\bigg)^{1/10}\bigg(\intop_{kT_0}^{(k+1)T_0}\intop_\Omega|w|^{20\over9}
dxdt\bigg)^{9/10}\cr
&\quad+\|g\|_{L_2(\Omega\times(kT_0,(k+1)T_0))}
\|w\|_{L_2(\Omega\times(kT_0,(k+1)T_0))}\bigg]\cr
&\le c[T_0d_7+d_7^2+\varphi(A)T_0^{7/30}A_2(T_0)+
\|g_0\|_{L_2(\Omega\times(kT_0,(k+1)T_0))}d_7]\cr
&\equiv A_3(T_0,d_7,A_2;\sup_k\|g_0\|_{L_2(\Omega\times(kT_0,(k+1)T_0))}).\cr}
$$
Setting
$$
X_2(t)=\intop_\Omega|w(t)|^2dx
\leqno(7.34)
$$
we obtain from (7.27) the inequality
$$
X_2((k+1)T_0)\le A_3(T_0)+e^{-\nu_*T_0}X_2(kT_0).
\leqno(7.35)
$$
Hence
$$
X_2(kT_0)\le{A_3(T_0)\over1-e^{-\nu_*T_0}}+e^{-\nu_*kT_0}X_2(0).
\leqno(7.36)
$$
Using the above estimate in (7.26) with $s=2$ we have
$$\eqal{
&\|w\|_{L_{10/3}(\Omega\times(kT_0,t))}^2\le(1+e^{\nu_*T_0})
\bigg({2-e^{-\nu_*T_0}\over1-e^{-\nu_*T_0}}A_3(T_0)\cr
&\quad+e^{-\nu_*kT_0}\|w(0)\|_{L_2(\Omega)}^2\bigg),\quad
t\in(kT_0,(k+1)T_0].\cr}
\leqno(7.37)
$$
From (7.37) it follows (7.15). This concludes the proof.

To prolong in time a local solution to (1.1) we need a version of Lemma 7.1
for solutions to problem (7.7). Therefore, we have

\proclaim Lemma 7.4.
Assume that $g=rf_\varphi\in L_\infty(\Omega\times\R_+)$,
$u(0)\in L_\infty(\Omega)$, $v(0)\in L_2(\Omega)$,
$f\in L_\infty(\R_+;L_{6/5}(\Omega))$,
${v_\varphi(0)\over r}\in L_2(\Omega)$. Let $T_0>0$ be given. Let
$|||f_\varphi|||_2=\sup_k\|f_\varphi\|_{L_2(\Omega\times(kT_0,(k+1)T_0))}<
\infty$,
$\big|\big|\big|{f_\varphi\over r}\big|\big|\big|_{9/5}=
\sup_k\big\|{f_\varphi\over r}\big\|_{L_{9/5}(\Omega\times(kT_0,(k+1)T_0))}
<\infty$. Let $A_*$ (introduced by (6.5)) be finite. Then
$$\eqal{
&\intop_{(k-1)T_0}^{kT_0}\intop_\Omega\tilde v_{\varphi,t}^2dxdt+\nu
\intop_\Omega\bigg(|\nabla\tilde v_\varphi(kT_0)|^2+
\bigg|{\tilde v_\varphi(kT_0)\over r}\bigg|^2\bigg)dx\cr
&\le cd_6^2+c\|\tilde f_\varphi\|_2^2+{c\over T_0^{3/2}}d_6d_7+\varphi
\bigg(A_*,T_0,d_6,d_7,\bigg|\bigg\|{f_\varphi\over r}\bigg\|\bigg|_{9/5},
|||f_\varphi|||_2\bigg)\cr
&\quad+ce^{-\nu_*kT_0}\bigg\|{v_\varphi(0)|\over r}\bigg\|_{L_2(\Omega)}^2
\equiv A_4^2.\cr}
\leqno(7.38)
$$
where $d_6$ is introduced in (3.1), $d_7$ in (2.15) and $A_*$ in (6.5).

\Proof
Multiplying $(7.7)_1$ by $\tilde v_{\varphi,t}$ and integrating the result
over $\Omega$ yield
$$\eqal{
&\intop_\Omega\tilde v_{\varphi,t}^2dx-\nu\intop_\Omega\Delta\tilde v_\varphi
\tilde v_{\varphi,t}dx+\nu\intop_\Omega{\tilde v_\varphi\over r^2}
\tilde v_{\varphi,t}dx\cr
&=\intop_\Omega\tilde f_\varphi\tilde v_{\varphi,t}dx+\intop_\Omega
v_\varphi\dot\zeta\tilde v_{\varphi,t}-\intop_\Omega{v_r\over r}
\tilde v_\varphi\tilde v_{\varphi,t}dx-\intop_\Omega v'\cdot\nabla
\tilde v_\varphi\tilde v_{\varphi,t}dx.\cr}
\leqno(7.39)
$$
Integration by parts implies
$$
-\intop_\Omega\Delta\tilde v_\varphi\tilde v_{\varphi,t}dx=-\intop_{-a}^a
\tilde v_\varphi\tilde v_{\varphi,t}|_{r=R}dz+\intop_\Omega\nabla
\tilde v_\varphi\cdot\nabla\tilde v_{\varphi,t}dx,
$$
where the boundary conditions $(7.7)_{2,3}$ were used.

\noindent
Applying the Cauchy inequality to the r.h.s. terms of (7.39) we get
$$\eqal{
&{1\over2}\intop_\Omega\tilde v_{\varphi,t}^2dx+{\nu\over2}{d\over dt}
\intop_\Omega|\nabla\tilde v_\varphi|^2dx+{\nu\over2}{d\over dt}
\intop_\Omega{\tilde v_\varphi^2\over r^2}dx\cr
&\quad-{\nu\over2}{d\over dt}\intop_{-a}^a\tilde v_\varphi^2|_{r=R}dz\le c
\intop_\Omega\tilde f_\varphi^2dx+c\intop_\Omega v_\varphi^2\dot\zeta^2dx\cr
&\quad+c\intop_\Omega v^{'2}|\nabla\tilde v_\varphi|^2dx+c\intop_\Omega
{v_r^2\over r^2}\tilde v_\varphi^2dx.\cr}
\leqno(7.40)
$$
Integrating (7.40) with respect to time from $t=(k-1)T_0$ to $t=kT_0$
we obtain
$$\eqal{
&{1\over2}\intop_{(k-1)T_0}^{kT_0}\intop_\Omega\tilde v_{\varphi,t}^2dxdt+
{\nu\over2}\intop_\Omega|\nabla\tilde v_\varphi(kT_0)|^2dx+{\nu\over2}
\intop_\Omega\bigg|{\tilde v_\varphi(kT_0)\over r}\bigg|^2dx\cr
&\le{\nu\over2}\intop_{-a}^a\tilde v_\varphi^2(kT_0)|_{r=R}dz+c
\intop_{(k-1)T_0}^{kT_0}\intop_\Omega\tilde f_\varphi^2dxdt\cr
&\quad+c\intop_{(k-1)T_0}^{kT_0}\intop_\Omega v_\varphi^2\dot\zeta^2dxdt+
c\intop_{(k1)T_0}^{kT_0}\bigg(v^{'2}|\nabla\tilde v_\varphi|^2+
{v_r^2\over r^2}\tilde v_\varphi^2\bigg)dxdt.\cr}
\leqno(7.41)
$$
The first term in estimated by $cd_6^2$, the second by the norm
$|||\tilde f_\varphi|||_2=\sup_k
\|\tilde f_\varphi\|_{L_2((k-1)T_0,kT_0;L_2(\Omega))}$, the third by
${c\over T_0^{3/2}}\|v_\varphi\|_{L_4(\Omega\times\R_+)}^2\le{c\over T_0^{3/2}}
d_6d_7$ (see (3.6)) and finally the last is bounded by
$$
\varphi\bigg(A_*,T_0,d_6,d_7,
\bigg|\bigg\|{f_\varphi\over r}\bigg\|\bigg|_{9/5},|||f_\varphi|||_2\bigg)+
ce^{-\nu_*kT_0}\bigg\|{v_\varphi(0)\over r}\bigg\|_{L_2(\Omega)}^2,
$$
where (7.10) and (7.15) were used.
Using the above estimates in (7.41) implies (7.38). This concludes the proof.

Having estimate (7.38) we are able to consider problem (7.1) in the time
interval $(kT_0,(k+1)T_0)$. Then instead of (7.2) we obtain
$$
\|v_{\varphi,t}\|_{L_2(kT_0,t;\Omega)}^2+\|v_\varphi(t)\|_{H_0^1(\Omega)}^2
\le cA_4^2,
\leqno(7.42)
$$
where $t\in(kT_0,(k+1)T_0)$ and $A_4$ is introduced in (7.38).

\section{8. Global existence}

To prove global existence of axially symmetric solutions to problem (1.1)
we first show long time existence (see Theorem 8.1). By the long time
existence we mean that there exists a time $T_e<\infty$, without any
restrictions on $T_e$ from above, such that axially symmetric solutions to
(1.1) satisfy
$$
v\in W_2^{2,1}(\Omega^{T_e}),\quad \nabla p\in L_2(\Omega^{T_e}).
\leqno(8.1)
$$
The existence of such time $T_e$ is connected with the assumption that some
spacial norms of the external force must be integrated in appropriate powers
over time interval $(0,T_e)$. Therefore, we prove global existence, so the
existence on the infinite time interval $(0,\infty)$, by prolonging the long
time solution step by step on intervals $((k-1)T_e,kT_e)$, $k\in\N$. In such
a way we relax restrictions on the external force, because otherwise looking
for solution (8.1) with $T_e=\infty$, we have to assume that $f$ vanishes
sufficiently fast as time goes to infinity. In the step by step approach we
prove existence of such solutions that
$$
v\in\bigcup_{k\in\N}W_2^{2,1}(\Omega\times((k-1)T_e,kT_e)),\quad
\nabla p\in\bigcup_{k\in\N}L_2(\Omega\times((k-1)T_e,kT_e)).
\leqno(8.2)
$$
Looking for solutions (8.2) we have to repeat the long time existence (see
Theorem 8.1) in each step, so we should prove that if there exists a constant
$\alpha$ depending on the data norms such that
$$
\|v(0)\|_{H^1(\Omega)}\le\alpha
\leqno(8.3)
$$
then
$$
\|v(kT_e)\|_{H^1(\Omega)}\le\alpha\quad {\rm for\ any}\ k\in\N.
\leqno(8.4)
$$
The above assertion holds for $v_\varphi$ (see Lemma 7.4 and formula (7.42)).
To prove it for $v'=(v_r,v_z)$ we need first some decay estimate for $\chi$.
This will be shown in Lemma 8.2.

\noindent
To formulate Theorem 8.1. we recall necessary notation. Let the constants
$d_0,d_1,\dots,d_7$ be introduced by the relations (see Lemma 2.1):
$$\eqal{
&\|v(0)\|_{L_2(\Omega)}\le d_0,\quad
\|f\|_{L_\infty(\R_+;L_{6/5}(\Omega))}\le d_2,\cr
&g=f\cdot\eta,\quad
u=v\cdot\eta,\quad \eta=(-x_2,x_1,0),\cr
&\sup_{k\in\N_0}\sup_{t\in[kT,(k+1)T]}\bigg|\intop_{kT}^t\intop_\Omega
gdxdt'+\intop_\Omega u(kT)dx\bigg|\le d_1,\cr
&d_3\ge c(d_1+d_2),\quad d_4\ge c(T)(d_0+d_3),\quad
d_5\ge c(T)(d_0+d_3)\cr
&d_7\ge c(T)(d_4+\sup_k\|f\|_{L_2(\Omega\times(kT,(k+1)T))})\quad
{\rm (see\ (2.15))}.\cr}
$$
$d_6=d_6(\|u(0)\|_{L_\infty(\Omega)},\|g\|_{L_\infty(\Omega\times\R_+)})$
(some function) (see (3.1)).

\noindent
Let us recall the partition of unity $\{\zeta_1(r),\zeta_2(r)\}$ on
the interval $(0,R)$ with the properties: $\zeta_1,\zeta_2$ are smooth
functions such that $\zeta_1(r)=1$ for $r\le r_0$ $\zeta_1(r)=0$ for
$r\ge2r_0$, $\zeta_2(r)=0$ for $r\le r_0$ and $\zeta_2(r)=1$ for $r\ge2r_0$.
Moreover, we have the notation: $\tilde w=w\zeta_1$, $\bar w=w\zeta_2$
for any function $w$.

\noindent
Let us introduce the quantities
$$\eqal{
&A^2(t)=c(1/r_0)d_5^2[1+(1+d_6)d_6^2]+
\bigg\|{\tilde\chi(0)\over r}\bigg\|_{L_2(\Omega)}^2\cr
&\quad+\bigg\|{\tilde v_\varphi(0)\over\sqrt{r}}\bigg\|_{L_4(\Omega)}^4+
\bigg\|{\tilde F\over r}\bigg\|_{L_2(\Omega^t)}^2+
\bigg\|{\tilde f_\varphi\over\sqrt{r}}\bigg\|_{L_{20/11}(\Omega^t)}^2,\cr
&A_0^2(t)=c(1/r_0)(d_5^2+d_6^2d_7^2)+\|\bar F\|_{L_2(0,t;L_{6/5}(\Omega))}^2+
{1\over r_0^2}\|\bar\chi(0)\|_{L_2(\Omega)}^2,\cr
&A_4^2(t)=cd_6^2+|\|\tilde f_\varphi\||_2^2+{c\over T^{3/2}}d_6d_7+
\varphi\bigg(A_*,T,d_6,d_7,
\bigg|\bigg\|{f_\varphi\over r}\bigg\|\bigg|_{9/5},\cr
&\qquad\quad\bigg|\bigg\|f_\varphi\bigg\|\bigg|_2\bigg)+c
\bigg\|{v_\varphi(0)\over r}\bigg\|_{L_2(\Omega)}^2,\cr}
\leqno(8.5)
$$
where $|||w|||_s=\sup_{k\in\N_0}\|w\|_{L_s((k-1)T,kT;L_s(\Omega))}$,
$s\in(1,\infty)$ and $A_*=A+[c(1/r_0)d_1^2+1]A_0$ (see (6.5)).

\proclaim Theorem 8.1. (long time existence).
Let $T_e>0$ be given. Assume that $v(0)\in H^1(\Omega)$,
$f\in L_2(\Omega^{T_e})$ and $A_1(T_e)$, $A_0(T_e)$, $A_4(T_e)$ are finite.
Then there exists axially symmetric solution to problem (1.1) such that
(8.1) holds.

\Proof
Assume that $T_*$ is so small that (2.27) is satisfied for $T=T_*$. Then
Lemma 2.5 implies a local existence of solutions to problem (1.1) such that
$v\in W_2^{2,1}(\Omega^{T_*})$, $\nabla p\in L_2(\Omega^{T_*})$. To apply
Lemma 2.5 we need that $v(0)\in H^1(\Omega)$ and $f\in L_2(\Omega^{T_*})$.
To extend the existence of local solutions on the interval $(T_*,2T_*)$
(or to apply Lemma 2.5 for interval $(T_*,2T_*))$ we have to obtain an
estimate for $\|v(T_*)\|_{H^1(\Omega)}$ such that
$$
\|v(T_*)\|_{H^1(\Omega)}\le\|v(0)\|_{H^1(\Omega)}.
\leqno(8.6)
$$
Moreover, we need that
$$
\|f\|_{L_2(\Omega\times(T_*,2T_*))}\le\|f\|_{L_2(\Omega^{T_*})}.
\leqno(8.7)
$$
Estimate (8.7) holds by assumptions. Therefore, we have to show (8.6) only.
Under assumptions (8.5) and from (6.5) we have the estimate
$$
\bigg\|{\chi\over r}\bigg\|_{V_2^0(\Omega^{T_*})}\le c(1/r_0,d_1)(A(T_e)+
A_0(T_e)).
\leqno(8.8)
$$
Then (6.43) (Lemma (6.5)) implies
$$
\|v'\|_{V_2^1(\Omega^{T_*})}\le c
\bigg\|{\chi\over r}\bigg\|_{V_2^0(\Omega^{T_*})}.
\leqno(8.9)
$$
Hence
$$
\|v'(T_*)\|_{H^1(\Omega)}\le c(1/r_0,d_1)(A(T_e)+A_0(T_e)).
\leqno(8.10)
$$
Similarly, (7.42) yields
$$
\|v_\varphi(T_*)\|_{H^1(\Omega)}\le cA_4(T_e).
\leqno(8.11)
$$
Let
$$
\alpha=c(A(T_e)+A_0(T_e)+A_4(T_e)).
\leqno(8.12)
$$
Since
$$
\|v(0)\|_{H^1(\Omega)}\le\alpha
\leqno(8.13)
$$
we obtain from (8.10) and (8.11) that
$$
\|v(T_*)\|_{H^1(\Omega)}\le\alpha
\leqno(8.14)
$$
and we can apply Lemma 2.5. Assuming that $T_*={T_e\over m}$ we can repeat
the procedure $m$ times and prove the theorem.

We can prove the theorem in a different way. We have existence of weak
solutions to problem (1.1) such that $v\in V_2^0(\Omega^{T_e})$ and the
estimate holds
$$
\|v\|_{V_2^0(\Omega^{T_e})}\le c(\|v(0)\|_{L_2(\Omega)}+
\|f\|_{L_2(\Omega^{T_e})}).
\leqno(8.15)
$$
We restrict an increasing of regularity of the weak solutions by getting an
estimate guaranteeing regularity (8.1). The precise procedure is very
complicated but it could be presented.

\noindent
Under assumptions (8.5) and Lemma 6.5 we have the estimate
$$
\|v'\|_{V_2^1(\Omega^{T_e})}\le c(A(T_e)+A_0(T_e)).
\leqno(8.16)
$$
From the proof of Lemma 3.7 in \cite{Z9} we have
$$
\|v'\|_{L_{10}(\Omega^{T_e})}\le c\|v'\|_{V_2^1(\Omega^{T_e})}.
\leqno(8.17)
$$
Let us consider the Stokes system
$$\eqal{
&v_t-\divv\T(v,p)=-v'\cdot\nabla v+f,\cr
&\divv v=0\cr
&v\cdot\bar n|_S=0,\quad \bar n\cdot\D(v)\cdot\bar\tau_\alpha|_S=0,\quad
\alpha=1,2,\cr
&v|_{t=0}=v(0),\cr}
\leqno(8.18)
$$
where the r.h.s. of $(8.18)_1$ is treated as given. In view of (8.15) and
(8.17)\break
 the r.h.s. of $(8.18)_1$ belongs to $L_{5/3}(\Omega^{T_e})$ under the
assumption that $f\in L_{5/3}(\Omega^{T_e})$ also. The last statement is
satisfied because $f\in L_2(\Omega^{T_e})$. Hence, in view of \cite{Z7, ZZ},
there exists a solution to (8.18) such that
$v\in W_{5/3}^{2,1}(\Omega^{T_e})$, $\nabla p\in L_{5/3}(\Omega^{T_e})$
and the estimate holds
$$
\|v\|_{W_{5/3}^{2,1}(\Omega^{T_e})}+\|\nabla p\|_{L_{5/3}(\Omega^{T_e})}\le
c({\rm data}).
\leqno(8.19)
$$
In view of the imbedding
$$
\|\nabla v\|_{L_{5/2}(\Omega^{T_e})}\le c\|v\|_{W_{5/3}^{2,1}(\Omega^{T_e})}
\leqno(8.20)
$$
we obtain that the r.h.s. of $(8.18)_1$ belongs to $L_2(\Omega^{T_e})$.
Applying again \cite{Z7, ZZ} we get that $v\in W_2^{2,1}(\Omega^{T_e})$,
$\nabla p\in L_2(\Omega^{T_e})$ and the estimate is valid
$$
\|v\|_{W_2^{2,1}(\Omega^{T_e})}+\|\nabla p\|_{L_2(\Omega^{T_e})}\le
c({\rm data}).
\leqno(8.21)
$$
This concludes the proof.

Next we shall show that if time existence $T_e$, appeared in Theorem 8.1, is
sufficiently large then there exists a solution to (1.1) with properties
described by (8.2). The solution behaves similarly in each time interval
$(kT_e,(k+1)T_e)$, $k\in\N_0$. To prove this we have to show that
$\|v(kT_e)\|_{H^1(\Omega)}$ is bounded by the same constant for all $k$.
For this some decay estimates are needed.

\proclaim Lemma 8.2.
Let $T_e>0$ be a sufficiently large given number. Let there exist a weak
solution to problem (1.1). Let assumptions of Lemmas 2.2, 2.3, 3.1, 3.2 be
satisfied. Let the constants $d_4$ (see (2.1)), $d_5$ (see (2.2)), $d_7$
(see (2.15)), $d_6$ (see (3.1)) be finite. Let $r_0>0$ be so small that
$$
\sup_t\|rv_\varphi\|_{L_\infty(\Omega_{r_0})}\le\root{4}\of{a\over3}\nu,\quad
a<1,
$$
holds, where $\Omega_{r_0}=\{x\in\Omega:\ r<r_0\}$. Let
$$
X^2(t)={1\over\nu^2}\bigg\|{v_\varphi^2\over r}\bigg\|_{L_2(\Omega)}^2+
\bigg\|{\chi\over r}\bigg\|_{L_2(\Omega)}^2.
$$
Let $X(0)$ be finite. Let the quantities
$$\eqal{
&d_8=\sup_{k\in\N_0}\|f_\varphi\|_{L_4(\Omega\times(kT_e,(k+1)T_e))},\cr
&d_9=\sup_{k\in\N_0}\bigg\|{F\over r}
\bigg\|_{L_2(\Omega\times(kT_e,(k+1)T_e))},\cr}
$$
be finite. Then there exists a positive increasing function\\
$B_0=B_0(1/r_0,d_4,d_5,d_6,d_7,d_8,d_9)$ such that
$$\eqal{
&X(kT_e)\le B_0+X(0)e^{-\nu_0kT_e},\cr
&\|v(t)\|_{H^1(\Omega)}\le c(B_0+X(kT_e)),\cr}
\leqno(8.22)
$$
hold for any $k\in\N_0$ and $t\in(kT_e,(k+1)T_e)$.

\Proof
From (4.12) we have
$$\eqal{
&{1\over2}{d\over dt}
\bigg\|{\tilde\chi\over r}\bigg\|_{L_2(\Omega_\varepsilon)}^2+{\nu\over2}
\bigg\|\nabla{\tilde\chi\over r}\bigg\|_{L_2(\Omega_\varepsilon)}^2\le
{1\over\nu}\bigg\|{\tilde v_\varphi\over r}\bigg\|_{L_4(\Omega_\varepsilon)}^4
\cr
&\quad+c(1/r_0)\intop_{\Omega_{\varepsilon,\zeta_{1,r}}}\chi^2dx+c(1/r_0)
\intop_{\Omega_{\varepsilon,\zeta_{1,r}}}|v_r|\chi^2dx+
c\intop_{\Omega_\varepsilon}\bigg|{\tilde F\over r}\bigg|^2dx.\cr}
\leqno(8.23)
$$
Moreover, using that $\tilde v_\varphi|_{r=\varepsilon}=0$, (4.16) yields
$$\eqal{
&{1\over4}{d\over dt}\intop_{\Omega_\varepsilon}
{\tilde v_\varphi^4\over r^2}dx+{3\over4}\nu\intop_{\Omega_\varepsilon}
\bigg|\nabla{\tilde v_\varphi^2\over r}\bigg|^2dx\cr
&\quad+{3\over4}\nu\intop_{\Omega_\varepsilon}
{\tilde v_\varphi^4\over r^4}dx=-{3\over2}\intop_{\Omega_\varepsilon}
{v_r\over r}{\tilde v_\varphi^4\over r^2}dx+\intop_{\Omega_\varepsilon}
v\cdot\nabla\zeta_1v_\varphi{\tilde v_\varphi\tilde v_\varphi^2\over r^2}dx\cr
&\quad-\intop_{\Omega_\varepsilon}[2\nabla v_\varphi\nabla\zeta_1+
v_\varphi\Delta\zeta_1]{\tilde v_\varphi\tilde v_\varphi^2\over r^2}dx+
\intop_{\Omega_\varepsilon}\tilde f_\varphi\tilde v_\varphi
{\tilde v_\varphi^2\over r^2}dx.\cr}
\leqno(8.24)
$$
Now we examine the particular terms from the r.h.s. of (8.24). We estimate
the first term by
$$
{3\over2}\|rv_\varphi\|_{L_\infty(\Omega_\zeta)}^2\intop_{\Omega_\varepsilon}
{|\tilde v_r|\over r^3}{\tilde v_\varphi^2\over r^2}dx,
$$
the second by
$$
{\varepsilon_1\over2}\intop_{\Omega_\varepsilon}{\tilde v_\varphi^4\over r^4}
dx+{1\over2\varepsilon_1}c(1/r_0)\intop_{\Omega_{\varepsilon,\zeta_{1,r}}}
v^2v_\varphi^4dx\equiv I_1,
$$
where the second term in $I_1$ is bounded by
$$
{c(1/r_0)\over2\varepsilon_1}\|rv_\varphi\|_{L_\infty(\Omega)}^4
\intop_{\Omega_{\varepsilon,\zeta_{1,r}}}v^2dx.
$$
The third integral on the r.h.s. of (8.24) is divided into two parts. The
first part equals
$$\eqal{
&-2\nu\intop_{\Omega_\varepsilon}\nabla v_\varphi\nabla\zeta_1
{\tilde v_\varphi\tilde v_\varphi^2\over r^2}dx=-2\nu\intop_{\Omega_\varepsilon}
\nabla v_\varphi\nabla\zeta_1{v_\varphi^3\zeta_1^3\over r^2}dx\cr
&=-{\nu\over2}\intop_{\Omega_\varepsilon}\nabla v_\varphi^4\nabla\zeta_1
{\zeta_1^3\over r^2}dx={\nu\over2}\intop_{\Omega_\varepsilon}v_\varphi^4
\nabla\bigg(\nabla\zeta_1{\zeta_1^3\over r^2}\bigg)dx\equiv I_2.\cr}
$$
Hence
$$
|I_2|\le c(1/r_0)\intop_{\Omega_{\varepsilon,\zeta_{1,r}}}v_\varphi^4dx.
$$
Similarly, the second part of the third term on the r.h.s. of (8.24) is
estimated by
$$
c(1/r_0)\intop_{\Omega_{\varepsilon,\zeta_{1,r}}}v_\varphi^4dx.
$$
Finally, the last term on the r.h.s. of (8.24) is bounded by
$$
{\varepsilon_2\over2}\intop_{\Omega_\varepsilon}
{\tilde v_\varphi^4\over r^4}dx+{1\over2\varepsilon_2}
\intop_{\Omega_\varepsilon}\tilde f_\varphi^2\tilde v_\varphi^2dx.
$$
Assuming $\varepsilon_1=\varepsilon_2={\nu\over8}$ and using the above
estimates in (8.24) imply
$$\eqal{
&{1\over4}{d\over dt}\intop_{\Omega_\varepsilon}{\tilde v_\varphi^4\over r^2}dx
+{3\over4}\nu\intop_{\Omega_\varepsilon}
\bigg|\nabla{\tilde v_\varphi^2\over r}\bigg|^2dx+{5\over8}\nu
\intop_{\Omega_\varepsilon}{\tilde v_\varphi^4\over r^4}dx\cr
&\le{3\over2}\|rv_\varphi\|_{L_\infty(\Omega_\zeta)}^2
\intop_{\Omega_\varepsilon}{|\tilde v_r|\over r^3}{\tilde v_\varphi^2\over r^2}
dx+c(1/r_0)\|rv_\varphi\|_{L_\infty(\Omega)}^4
\intop_{\Omega_{\varepsilon,\zeta_{1,r}}}v^2dx\cr
&\quad+c(1/r_0)\intop_{\Omega_{\varepsilon,\zeta_{1,r}}}v_\varphi^4dx+
{4\over\nu}\intop_{\Omega_\varepsilon}\tilde f_\varphi^2\tilde v_\varphi^2dx.
\cr}
\leqno(8.25)
$$
Employing (4.29) in (8.23) gives
$$\eqal{
&{1\over2}{d\over dt}
\bigg\|{\tilde\chi\over r}\bigg\|_{L_2(\Omega_\varepsilon)}^2+{\nu_1\over2}
\bigg\|\nabla{\tilde\chi\over r}\bigg\|_{L_2(\Omega_\varepsilon)}^2+
{\nu_2\over2}\bigg(\intop_{\Omega_\varepsilon}\bigg|\nabla\bigg(
{\tilde v_r\over r}\bigg)_{,r}\bigg|^2dx\cr
&\quad+3\nu_2\intop_{\Omega_\varepsilon}{1\over r^2}
\bigg({\tilde v_r\over r}\bigg)_{,r}^2dx\bigg)\le{1\over\nu}
\intop_{\Omega_\varepsilon}{\tilde v_\varphi^4\over r^4}dx+c(1/r_0)
\intop_{\Omega_{\varepsilon,\zeta_{1,r}}}\chi^2dx\cr
&\quad+c(1/r_0)\intop_{\Omega_{\varepsilon,\zeta_{1,r}}}|v_r|\chi^2dx+c
\intop_{\Omega_\varepsilon}\bigg|{\tilde F\over r}\bigg|^2dx\cr
&\quad+{\nu_2\over2}
c(1/r_0)\intop_{\Omega_\varepsilon}(v_{z,r}^2+v_z^2)dx,\cr}
\leqno(8.26)
$$
where $\nu=\nu_1+\nu_2$, $\nu_i>0$, $i=1,2$.
Using (3.1) and applying the Cauchy inequality to the first term on the r.h.s.
of (8.25) we get
$$\eqal{
&{1\over4}{d\over dt}\intop_{\Omega_\varepsilon}
{\tilde v_\varphi^4\over r^2}dx+{3\over4}\nu\intop_{\Omega_\varepsilon}
\bigg|\nabla{\tilde v_\varphi^2\over r}\bigg|^2dx+\bigg({5\over8}\nu-
{\varepsilon\over2}\bigg)\intop_{\Omega_\varepsilon}
{\tilde v_\varphi^4\over r^4}dx\cr
&\le{9\over8\varepsilon}\|rv_\varphi\|_{L_\infty(\Omega_\zeta)}^4
\intop_{\Omega_\varepsilon}{\tilde v_r^2\over r^6}dx+c(1/r_0)d_6^4
\intop_{\Omega_{\varepsilon,\zeta_{1,r}}}v^2dx\cr
&\quad+c(1/r_0)\intop_{\Omega_{\varepsilon,\zeta_{1,r}}}v_\varphi^4dx+
c\intop_{\Omega_\varepsilon}\tilde f_\varphi^4dx.\cr}
\leqno(8.27)
$$
Setting $\varepsilon={\nu\over4}$ and using the Hardy inequality
$$
\intop_{\Omega_\varepsilon}{\tilde v_r^2\over r^6}dx\le
\intop_{\Omega_\varepsilon}{1\over r^2}\bigg({\tilde v_r\over r}\bigg)_{,r}^2dx
$$
in (8.27) gives
$$\eqal{
&{1\over4}{d\over dt}\intop_{\Omega_\varepsilon}
{\tilde v_\varphi^4\over r^2}dx+{3\over4}\nu\intop_{\Omega_\varepsilon}
\bigg|\nabla{\tilde v_\varphi^2\over r}\bigg|^2dx+{\nu\over2}
\intop_{\Omega_\varepsilon}{\tilde v_\varphi^4\over r^4}dx\cr
&\le{9\over2\nu}\|rv_\varphi\|_{L_\infty(\Omega_\zeta)}^4
\intop_{\Omega_\varepsilon}{1\over r^2}\bigg({\tilde v_r\over r}\bigg)_{,r}^2
dx+c(1/r_0)d_6^4\intop_{\Omega_{\varepsilon,\zeta_{1,r}}}v^2dx\cr
&\quad+c(1/r_0)\intop_{\Omega_{\varepsilon,\zeta_{1,r}}}v_\varphi^4dx+
c\intop_{\Omega_\varepsilon}\tilde f_\varphi^4dx.\cr}
\leqno(8.28)
$$
Multiplying (8.28) by ${2\over\nu^2}$ and adding to (8.26) we obtain
$$\eqal{
&{1\over2\nu^2}{d\over dt}\intop_{\Omega_\varepsilon}
{\tilde v_\varphi^4\over r^2}dx+{1\over2}{d\over dt}\intop_{\Omega_\varepsilon}
{\tilde\chi^2\over r^2}dx+{3\over2\nu}\intop_{\Omega_\varepsilon}
\bigg|\nabla{\tilde v_\varphi^2\over r}\bigg|^2dx\cr
&\quad+{\nu_1\over2}\intop_{\Omega_\varepsilon}\bigg|\nabla
{\tilde\chi\over r}\bigg|^2dx+3\nu_2\intop_{\Omega_\varepsilon}{1\over r^2}
\bigg({\tilde v_r\over r}\bigg)_{,r}^2dx\cr
&\le{9\over\nu^3}\|rv_\varphi\|_{L_\infty(\Omega_\zeta)}^4
\intop_{\Omega_\varepsilon}{1\over r^2}\bigg({\tilde v_r\over r}\bigg)_{,r}^2
dx+c(1/r_0)\intop_{\Omega_\varepsilon}\chi^2dx\cr
&\quad+c(1/r_0)\intop_{\Omega_{\varepsilon,\zeta_{1,r}}}|v_r|\chi^2dx+
c(1/r_0)\intop_{\Omega_\varepsilon}(v_{z,r}^2+v_z^2)dx\cr
&\quad+c(1/r_0)d_6^2\intop_{\Omega_{\varepsilon,\zeta_{1,r}}}v^2dx+c(1/r_0)
\intop_{\Omega_\varepsilon}v_\varphi^4dx+c\intop_{\Omega_\varepsilon}
\bigg|{\tilde F\over r}\bigg|^2dx\cr
&\quad+c\intop_{\Omega_\varepsilon}\tilde f_\varphi^4dx.\cr}
\leqno(8.29)
$$
Setting $\nu_2=a\nu$, $a<1$ and assuming that
$$
{9\over\nu^3}\|rv_\varphi\|_{L_\infty(\Omega\zeta)}^4\le3a\nu\quad {\rm so}
\quad \|rv_\varphi\|_{L_\infty(\Omega_\zeta)}^4\le{a\nu^4\over3}.
\leqno(8.30)
$$
For $a={16\over27}$ the above condition takes the form
$$
\|rv_\varphi\|_{L_\infty(\Omega_\zeta)}\le{2\nu\over3},
$$
we obtain from (8.29) the inequality
$$\eqal{
&{d\over dt}\bigg({1\over\nu^2}
\bigg\|{\tilde v_\varphi^2\over r}\bigg\|_{L_2(\Omega_\varepsilon)}^2+
\bigg\|{\tilde\chi\over r}\bigg\|_{L_2(\Omega_\varepsilon)}^2\bigg)+
\bigg({3\over\nu}\bigg\|\nabla
{\tilde v_\varphi^2\over r}\bigg\|_{L_2(\Omega_\varepsilon)}^2\cr
&\quad+\nu_1\bigg\|\nabla{\tilde\chi\over r}\bigg\|_{L_2(\Omega_\varepsilon)}^2
\bigg)\le c(1/r_0)\intop_{\Omega_\varepsilon}\chi^2dx+c(1/r_0)
\intop_{\Omega_{\varepsilon,\zeta_{1,r}}}|v_r|\chi^2dx\cr
&\quad+c(1/r_0)\intop_{\Omega_\varepsilon}(v_{z,r}^2+v_z^2)dx+c(1/r_0)d_6^2
\intop_{\Omega_\varepsilon}v^2dx+c(1/r_0)\intop_{\Omega_\varepsilon}
v_\varphi^4dx\cr
&\quad+c\intop_{\Omega_\varepsilon}\tilde f_\varphi^4dx+c
\intop_{\Omega_\varepsilon}\bigg|{\tilde F\over r}\bigg|^2dx\cr}
\leqno(8.31)
$$
Let us recall that $\nu_1=(1-a)\nu$. Let $b<1-a$.

\noindent
Then the second expression on the l.h.s. takes the form
$$
\nu\bigg({3\over\nu^2}\bigg\|\nabla
{\tilde v_\varphi^2\over r}\bigg\|_{L_2(\Omega_\varepsilon)}^2+(1-a-b)
\bigg\|\nabla{\tilde\chi\over r}\bigg\|_{L_2(\Omega_\varepsilon)}^2\bigg)+
b\nu\bigg\|\nabla{\tilde\chi\over r}\bigg\|_{L_2(\Omega_\varepsilon)}^2
$$
Let $c_p$ be the constant from the Poincare inequality.\\
Let $\nu_*=\min\big\{{3\nu\over c_p},{(1-a-b)\nu\over c_p}\big\}$ and let
$$
X(t)={1\over\nu^2}
\bigg\|{\tilde v_\varphi^2\over r}\bigg\|_{L_2(\Omega_\varepsilon)}^2+
\bigg\|{\tilde\chi\over r}\bigg\|_{L_2(\Omega_\varepsilon)}^2.
\leqno(8.32)
$$
Then (8.31) takes the form
$$\eqal{
&{d\over dt}X+\nu_*X+b\nu\bigg\|\nabla
{\tilde\chi\over r}\bigg\|_{L_2(\Omega_\varepsilon)}^2\le c(1/r_0)
\intop_{\Omega_{\varepsilon,\zeta_{1,r}}}|v_r|\chi^2dx\cr
&\quad+c(1/r_0)\intop_{\Omega_\varepsilon}(\chi^2+v_{z,r}^2+v_z^2+d_6^2v^2)dx
+c(1/r_0)\intop_{\Omega_\varepsilon}v_\varphi^4dx\cr
&\quad+c\intop_{\Omega_\varepsilon}\tilde f_\varphi^4dx+
c\intop_{\Omega_\varepsilon}\bigg|{\tilde F\over r}\bigg|^2dx.\cr}
\leqno(8.33)
$$
Multiplying (8.33) by $e^{\nu_*t}$ and integrating with respect to time from
$kT$ to $t\in(kT_e,(k+1)T_e]$, $k\in\N_0$, yields
$$\eqal{
&X(t)e^{\nu_*t}+b\nu\intop_{kT}^t\bigg\|\nabla{\tilde\chi\over r}
(t')\bigg\|_{L_2(\Omega_\varepsilon)}^2e^{\nu_*t'}dt'\cr
&\le c(1/r_0)\intop_{kT_e}^t\intop_{\Omega_{\varepsilon,\zeta_{1,r}}}
|v_r|\chi^2dxe^{\nu_*t'}dt'+c(1/r_0)d_5^2(1+d_6^2)e^{\nu_*t}\cr
&\quad+c(1/r_0)d_6^2d_7^2e^{\nu_*t}+
(\|f_\varphi\|_{L_4(\Omega\times(kT_e,(k+1)T_e))}^4\cr
&\quad+\bigg\|{\tilde F\over r}\bigg\|_{L_2(\Omega\times(kT_e,(k+1)T_e))}^2
\bigg)e^{\nu_*t}+X(kT)e^{\nu_*kT_e}.\cr}
\leqno(8.34)
$$
Introducing the notation
$$
d_8=\sup_k\|f_\varphi\|_{L_4(\Omega\times(kT_e,(k+1)T_e))},\quad
d_9=\sup_k\bigg\|{\tilde F\over r}\bigg\|_{L_2(\Omega\times(kT_e,(k+1)T_e))}
$$
and estimating the first term on the r.h.s. of (8.34) by
$$\eqal{
&c(1/r_0)\intop_{kT_e}^t\|v_r\|_{L_2(\Omega_{\varepsilon,\zeta_{1,r}})}
\|\chi\|_{L_4(\Omega_{\varepsilon,\zeta_{1,r}})}^2e^{\nu_*t'}dt'\cr
&\le c(1/r_0)d_4\intop_{kT_e}^t
\|\chi\|_{L_4(\Omega_{\varepsilon,\zeta_{1,r}})}^2e^{\nu_*t'}dt'\cr
&\le\varepsilon\intop_{kT_e}^t\bigg\|\nabla
{\chi\over r}\bigg\|_{L_2(\Omega_{\varepsilon,\zeta_{1,r}})}^2e^{\nu_*t'}dt'+
c(1/r_0,d_4,d_5)e^{\nu_*t}\cr}
$$
we obtain
$$\eqal{
&X(t)e^{\nu_*t}+b\nu\intop_{kT_e}^t\bigg\|\nabla
{\chi\over r}(t')\bigg\|_{L_2(\Omega_\varepsilon)}^2e^{\nu_*t'}dt'\cr
&\le\varepsilon\intop_{kT_e}^t\bigg\|\nabla
{\chi\over r}\bigg\|_{L_2(\Omega_{\varepsilon,\zeta_{1,r}})}^2e^{\nu_*t'}dt'+
B(1/r_0,d_4,d_5,d_6,d_7,d_8,d_9,1/\varepsilon)\cdot\cr
&\quad\cdot e^{\nu_*t}+X(kT)e^{\nu_*kT_e},\cr}
\leqno(8.35)
$$
where $B$ is a positive increasing function of its arguments.

\noindent
By the local iteration technique (see \cite[Ch. 4, Sect. 10]{LSU}) we get
$$
X(t)e^{\nu_*t}+b\nu\intop_{kT_e}^t\bigg\|\nabla
{\chi\over r}\bigg\|_{L_2(\Omega_\varepsilon)}^2e^{\nu_*t'}dt'\le
2Be^{\nu_*t}+2X(kT)e^{\nu_*kT}
\leqno(8.36)
$$
To obtain an estimate for $X(kT_e)$ we skip the second term on the l.h.s. of
(8.36) and set $t=(k+1)T_e$. Then we get
$$
X((k+1)T_e)\le2B+2X(kT)e^{-\nu_*T_e}.
\leqno(8.37)
$$
Let $B_0=2B$ and $2e^{-{\nu_*\over2}T_e}\le1$. Denoting $\nu_0={\nu_*\over2}$
we obtain
$$
X((k+1)T_e)\le B_0+X(kT_e)e^{-\nu_0T_e}.
\leqno(8.38)
$$
Hence
$$
X(kT_e)\le{B_0\over1-e^{-\nu_0T_e}}+X(0)e^{-\nu_0kT_e},\quad k\in\N_0.
\leqno(8.39)
$$
In view of (8.39) we obtain from (8.31) after integration with respect to
time from $kT_e$ to $t\in(kT_e,(k+1)T_e]$ the estimate
$$\eqal{
&\bigg\|{v_\varphi^2\over r}\bigg\|_{V_2^0(\Omega\times(kT_e,t))}^2+
\bigg\|{\chi\over r}\bigg\|_{V_2^0(\Omega\times(kT_e,t))}^2\cr
&\le B_0+X^2(kT_e)e^{-\nu_*(t-kT_e)}\cr}
\leqno(8.40)
$$
The inequality implies the estimate
$$
\|v(t)\|_{H^1(\Omega)}\le c(B_0+\chi(kT_e)),\quad t\in(kT_e,(k+1)T_e],
\leqno(8.41)
$$
for any $k\in\N_0$. This concludes the proof.

\noindent
{\bf Proof of the Main Theorem.}
Let $T=T_e$. In view of Theorem 8.1 and Lemma 8.2 there exists a local
solution to problem (1.1) such that
$$
v\in W_2^{2,1}(\Omega\times(kT,(k+1)T)),\quad
\nabla p\in L_2(\Omega\times(kT,(k+1)T))
$$
and estimate (1.15) is valid. This ends the proof.

\section{References}

\item{BIN} Besov, O. V.; Il'in, V. P.; Nikolskij, S. M.: Integral
representation of functions and theorems of imbedding, Nauka, Moscow 1975
(in Russian).

\item{B} Burczak, J.: Boundary De Giorgi-Ladyzhenskaya classes and their
application to regularity of swirl of the Navier-Stokes equations (to be
published).

\item{CKN} Caffarelli, L.; Kohn. R.; Nirenberg, L.: Partial regularity of
suitable weak solutions of the Navier-Stokes equations, Comm. Pure Appl.
Math. 35 (1982), 771--831.

\item{K} Kochin, N. E.: Vectorial calculations and beginning of tensor
calculations, Akademia Nauk SSSR, Moscow 1951 (in Russian).

\item{LSU} Ladyzhenskaya, O. A.; Solonnikov, V. A.; Uraltseva, N. N.:
Linear and quasilinear equations of parabolic type, Nauka, Moscow 1967
(in Russian).

\item{L1} Ladyzhenskaya, O. A.: On unique solvability of three-dimensional
Cauchy problem for the Navier-Stokes equations under the axial symmetry,
Zap. Nauchn. Sem. LOMI 7 (1968), 155--177 (in Russian).

\item{L2} Ladyzhenskaya, O. A.: Mathematical problem of motion of viscous
incompressible fluid, Nauka, Moscow 1970 (in Russian).

\item{LL} Landau, L.; Lifshitz, E.: Hydrodynamics, Nauka, Moscow 1986 (in
Russian).

\item{S} Solonnikov, V. A.: Estimates of the solutions of a nonstationary
linearized system of the Navier-Stokes equations, Trudy Mat. Inst. Steklov
70 (1964), 213--317; English transl. Amer. Math. Soc. Trans. Ser. 2, 65
(1967), 51--137.

\item{T} Temam, R.: Navier-Stokes equations, theory and numerical analysis,
North-Holland Publishing Company Amsterdam-New York-Oxford, 1979.

\item{UY} Ukhovskij, M. R.; Yudovich, V. I.: Axially symmetric motions of
ideal and viscous fluids filling all space, Prikl. Mat. Mech. 32 (1968),
59--69 (in Russian).

\item{ZZ} Zadrzy\'nska, E.; Zaj\c aczkowski, W. M.: Nonstationary Stokes
system in Sobolev spaces (to be published).

\item{Z1} Zaj\c aczkowski, W. M.: A priori estimate for axially symmetric
solutions to the Navier-Stokes equations near the axis of symmetry, arXiv:
1210.1510, math. AP.

\item{Z2} Zaj\c aczkowski, W. M.: Global regular axially symmetric solutions
to the Navier-Stokes equations in a periodic cylinder, arXiv: 1210. 1513,
math. AP.

\item{Z3} Zaj\c aczkowski, W. M.: The H\"older regularity of swirl, arXiv:\break
1210.1515, math. AP.

\item{Z4} Zaj\c aczkowski, W. M.: Global special regular solutions to the
Navier-Stokes equations in a cylindrical domain under boundary slip
condition, Gakuto Intern. Ser., Math. Sc. Appl. 21 (2004), 1--188.

\item{Z5} Zaj\c aczkowski,W. M.: Global special regular solutions to the
Navier-Stokes equations in axially symmetric domains under boundary slip
conditions, Diss. Math. 432 (2005), pp. 138.

\item{Z6} Zaj\c aczkowski, W. M.: Global regular solutions to the
Navier-Stokes equations in an axially symmetric domain, Topol. Meth. Nonlin.
Anal. 33 (2009), 233--274.

\item{Z7} Zaj\c aczkowski, W. M.: Nonstationary Stokes system in cylindrical
domains and anisotropic Sobolev spaces.

\item{Z8} Zaj\c aczkowski, W. M.: Global axially symmetric solutions with
large swirl to the Navier-Stokes equations, Topol. Meth. Nonlin. Anal. 29
(2007), 295--331.

\item{Z9} Zaj\c aczkowski, W. M.: Global special regular solutions to the
Navier-Stokes equations in a cylindrical domain without the axis of symmetry,
Top. Meth. Nonlin. Anal. 24 (2004), 69--105.

\bye